\documentclass[11pt,twoside]{amsart}
\usepackage{amssymb,rotating}

\newcommand{\QQ}{{\mathbb{Q}}}
\newcommand{\FF}{{\mathbb{F}}}
\newcommand{\CC}{{\mathbb{C}}}
\newcommand{\ZZ}{{\mathbb{Z}}}
\newcommand{\NN}{{\mathbb{N}}}

\newcommand{\Sym}{{\mathfrak{S}}}

\newcommand{\fC}{{\mathfrak{C}}}
\newcommand{\fF}{{\mathfrak{F}}}
\newcommand{\fp}{{\mathfrak{p}}}
\newcommand{\fh}{{\mathfrak{h}}}
\newcommand{\fg}{{\mathfrak{g}}}
\newcommand{\fsl}{{\mathfrak{sl}}}
\newcommand{\cH}{{\mathcal{H}}}
\newcommand{\cB}{{\mathcal{B}}}
\newcommand{\cD}{{\mathcal{D}}}
\newcommand{\cL}{{\mathcal{L}}}
\newcommand{\cS}{{\mathcal{S}}}
\newcommand{\cM}{{\mathcal{M}}}
\newcommand{\cN}{{\mathcal{N}}}
\newcommand{\cW}{{\mathcal{W}}}

\newcommand{\tX}{{\tilde{X}}}

\newcommand{\tT}{{\tilde{T}}}

\newcommand{\ba}{{\mathbf{a}}}
\newcommand{\bc}{{\mathbf{c}}}
\newcommand{\bu}{{\mathbf{u}}}

\newcommand{\bG}{{\mathbf{G}}}
\newcommand{\bA}{{\mathbf{A}}}
\newcommand{\bB}{{\mathbf{B}}}
\newcommand{\bH}{{\mathbf{H}}}
\newcommand{\bJ}{{\mathbf{J}}}
\newcommand{\bLambda}{{\boldsymbol{\Lambda}}}

\newcommand{\ulambda}{{\boldsymbol{\lambda}}}
\newcommand{\umu}{{\boldsymbol{\mu}}}
\newcommand{\uvar}{{\underline{\varnothing}}}
\newcommand{\bN}{{\mathbf{N}}}
\newcommand{\bT}{{\mathbf{T}}}
\newcommand{\bU}{{\mathbf{U}}}

\newcommand{\bW}{{\mathbf{W}}}

\newcommand\modA{\operatorname{-mod}}
\newcommand\ModA{\operatorname{-Mod}}
\newcommand\Irr{\operatorname{Irr}}

\newcommand\Ind{\operatorname{Ind}}

\newcommand\Hom{\operatorname{Hom}}
\newcommand\End{\operatorname{End}}
\newcommand\Fix{\operatorname{Fix}}
\newcommand\Rst{{^*\!R}}
\newcommand\GL{\operatorname{GL}}

\newcommand\GU{\operatorname{GU}}
\newcommand\Uch{\operatorname{Uch}}

\renewcommand{\leq}{\leqslant}
\renewcommand{\geq}{\geqslant}

\renewcommand{\atop}[2]{\genfrac{}{}{0pt}{}{#1}{#2}}
\def\dddots{\mathinner{\mkern1mu\raise1pt
    \vbox{\kern7pt\hbox{.}}\mkern2mu
    \raise4pt\hbox{.}\mkern2mu\raise7pt\hbox{.}\mkern1mu}}

\newtheorem{thm}{Theorem}[section]
\newtheorem{lem}[thm]{Lemma}
\newtheorem{cor}[thm]{Corollary}
\newtheorem{prop}[thm]{Proposition}
\newtheorem{conj}[thm]{Conjecture}

\theoremstyle{definition}
\newtheorem{defn}[thm]{Definition}
\newtheorem{exmp}[thm]{Example}
                                                                            
\theoremstyle{remark}
\newtheorem{rem}[thm]{Remark}
                                                                            
\newenvironment{Par}{\refstepcounter{thm} \medskip \noindent
{\bf \thethm.}}{\par\smallskip}

\begin{document}

\title{Modular representations of Hecke algebras}

\author{Meinolf Geck} 

\email{m.geck@maths.abdn.ac.uk}

\address{Institut Camille Jordan, bat. Jean 
Braconnier, Universit\'e Lyon 1, 21 av Claude Bernard, 
F--69622 Villeurbanne Cedex, France}

\curraddr{Department of Mathematical Sciences, King's College,
Aberdeen University, Aberdeen AB24 3UE, Scotland, U.K.}

\date{}


\begin{abstract} These notes are based on a course given at the EPFL in 
May 2005. It is concerned with the representation theory of Hecke algebras 
in the non-semisimple case. We explain the role that these algebras play 
in the modular representation theory of finite groups of Lie type and 
survey the recent results which complete the classification of the simple 
modules. These results rely on the theory of Kazhdan--Lusztig cells in 
finite Weyl groups (with respect to possibly unequal parameters) and the 
theory of canonical bases for representations of quantum groups.
\end{abstract}

\date{}
\maketitle

\pagestyle{myheadings}

\markboth{Geck}{Modular representations of Hecke algebras}

\section{Introduction} \label{MGsec1}

The theory and the results that we are going to talk about can be seen
as a contribution to the general project of determining the irreducible
representations of all finite simple groups. Recall that such a group is 
either cyclic of prime order, or an alternating group of degree $\geq 5$, 
or a simple group of Lie type,  or one of $26$ sporadic simple groups. 
Here, we will concentrate on the finite groups of Lie type; any such group 
is naturally defined over a finite field $\FF_q$, where $q$ is a prime 
power. We shall consider representations over a field whose characteristic 
is a prime number not dividing $q$, the ``non-defining characteristic case''. 
Iwahori--Hecke algebras associated with finite Coxeter groups naturally 
arise in this context, as endomorphism algebras of certain induced 
representations. This will be explained in Sections~2 and~3, where we give
an introduction to Harish--Chandra series in the ``modular case''. For 
further details and references, we refer to the survey \cite{mycharl}. 

Now the endomorphism algebras arising in this context can be defined
abstractly, in terms of generators and relations. We will study the 
representations of Iwahori--Hecke algebras in this abstract setting, 
without reference to the realization as an endomorphism algebra of some 
induced representation; see Section~4. Our main aim is to explain a natural 
parametrization of the irreducible representations in terms of so-called 
``canonical basic sets'' and the unitriangularity of the decomposition
matrix, following work of Rouquier and the author \cite{mykl}, 
\cite{GeRo2}, \cite{my00}. This heavily relies on the theory of 
Kazhdan--Lusztig cells; see Sections~5 and~6. 

Originally, we only considered ``canonical basic sets'' for 
Iwahori--Hecke algebras with equal parameters. Here, we present this 
theory in the general framework of possibly unequal parameters, following 
Lusztig \cite{Lusztig03}. Some properties of the Kazhdan--Lusztig basis are 
only conjectural in the general case, but the recent work of Bonnaf\'e, 
Iancu and the author \cite{BI}, \cite{BI2}, \cite{GeIa05}, \cite{my05} 
shows that these conjectural properties hold (at least) in the so-called
``asymptotic case'' in type $B_n$. For example, the characterization of 
the canonical basic set in this case in Example~\ref{canbn} is a new result. 

Once the existence of a natural parametrization of the simple modules is
established, it is another question to determine those ``canonical
basic sets'' explicitly. In type $B_n$ (in many ways the hardest case) this
was achieved by Jacon \cite{Jac0}, \cite{Jac2}, \cite{Jac4}, using the 
theory of canonical bases of quantum groups and building on earlier work
of Ariki \cite{Ar}, \cite{Ar3} (the proof of the LLT conjecture), 
Ariki--Mathas \cite{ArMa}, Foda et al.  \cite{FLOTW}, and Uglov \cite{Ug}. 
These results guarantee the existence of a ``canonical basic set'' even 
for those choices of unequal parameters where Lusztig's conjectural 
properties are not (yet) known to hold. Thus, quite remarkably, irreducible 
representations of Iwahori--Hecke algebras of classical type are naturally 
labelled by the crystal bases of certain highest weight modules for the
quantized enveloping algebra $U_v(\hat{\fsl}_l)$. All this will be discussed 
in Sections~7 and~8.

\section{Harish-Chandra series and Hecke algebras} \label{MGsec2}

Let $G$ be a finite group with a split BN-pair of characteristic~$p$ which 
satisfies Chevalley's commutator relations; see \cite[\S 2.5--2.6]{Ca2} 
or \cite[\S 1.6]{myag}. We don't repeat all the axioms here, but just recall 
some basic results. We have $B=UH$ where $U=\operatorname{O}_p(B)$ is the 
largest normal $p$-subgroup of $B$ and $H$ is an abelian $p'$-subgroup such 
that $H=B \cap N$. Let $W=N/H$ be the Weyl group of $G$. For any $w \in W$ 
we denote by $\dot{w}$ a representative in~$N$. The group $W$ is a Coxeter 
group with respect to the set of generators 
\[S=\{ w \in W \mid \mbox{$w \neq 1$ and $B \cup B\dot{w}B$ is a 
subgroup of $G$}\}.\]
For any subset $I \subseteq S$, we have a corresponding parabolic subgroup 
$W_I=\langle I \rangle \subseteq W$.  Moreover, $I$ also defines a parabolic 
subgroup $P_I=BN_IB \subseteq G$, where $N_I=\{H\dot{w} \mid w\in W_I\}$. 
We have a Levi decomposition $P_I=U_IL_I$ where $U_I=\operatorname{O}_p
(P_I)$ is the largest normal $p$-subgroup of $P_I$ and $L_I$ is a 
complementary subgroup which is uniquely determined by the condition 
that $H \subseteq L_I$. The group $L_I$ is called a standard Levi 
subgroup of $G$; it is itself a finite group with a split BN-pair of 
characteristic~$p$ corresponding to the subgroups $B_I=B \cap L_I$ and 
$N_I=\{H\dot{w} \mid w\in W_I\}$. 

\begin{exmp} \label{gln1} Finite groups with a split BN-pair typically 
arise as the fixed point sets of connected reductive algebraic groups 
under a Frobenius map. (A formal definition of such maps is contained
in \cite[\S 4.1]{myag}.) Here is the standard example. Let $p$ be a prime 
and $\overline{\FF}_p$ be an algebraic closure of the finite field with 
$p$ elements. Let $\bG= \GL_n(\overline{\FF}_p)$ be the general linear 
group of $n\times n$-matrices over $\overline{\FF}_p$. Let $\bB\subseteq 
\bG$ be the subgroup consisting of all upper triangular matrices, and 
$\bN\subseteq \bG$ be the subgroup consisting of all monomial matrices. 
Then it is well-known that the groups $\bB$ and $\bN$ form a BN-pair in 
$\bG$, with Weyl group $\bW \cong \Sym_n$, the symmetric group on $n$ 
letters; see \cite[\S 1.6]{myag} for more details. We have a semidirect
product decomposition $\bB=\bU\rtimes \bT_0$ where $\bU$ is the normal
subgroup consisting of all upper triangular matrices with $1$ on the
diagonal, and $\bT_0$ is the subgroup consisting of all invertible
diagonal matrices; then $\bN=\operatorname{N}_{\bG}(\bT_0)$, the normalizer
of $\bT_0$ in $\bG$, and $\bW=\bN/\bT_0$. 

(a) Let $q=p^f$ for some $f\geq 1$. Then we have a unique subfield $\FF_q 
\subseteq \overline{\FF}_p$ with $q$ elements. We have the ``standard''
Frobenius map 
\[ F_q \colon \bG \rightarrow \bG, \qquad (a_{ij}) \mapsto (a_{ij}^q).\]
The group of fixed points is $\bG^{F_q}=\GL_n(\FF_q)$, the general linear 
group over $\FF_q$.  The groups $\bB$ and $\bN$ are $F_q$-stable. By 
taking fixed points under $F_q$, we obtain that $\bB^{F_q}$ and $\bN^{F_q}$
form a split BN-pair of characteristic~$p$ in the finite group $\GL_n(\FF_q)$,
with Weyl group $W \cong \Sym_n$; see \cite[\S 4.2]{myag}.

(b) Now consider the permutation matrix 
\[ Q_n:= \left[\begin{array}{cccc} 0 & \cdots & 0 & 1 \\ \vdots &\dddots &
\dddots & 0 \\ 0 & 1 &\dddots & \vdots \\ 1 & 0 & \cdots & 0 \end{array}
\right] \in M_n(k)\]
and define an automorphism of algebraic groups $\gamma\colon 
\bG\rightarrow \bG$ by 
\[\gamma(A):=Q_n^{-1}\cdot (A^{\text{tr}})^{-1}\cdot Q_n \qquad \mbox{where
$A \in \GL_n(\overline{\FF}_p)$}.\]
Then $\gamma$ commutes with $F_q$ and $\gamma^2$ is the identity. Hence, 
the map $F:=F_q \circ \gamma$ also is a Frobenius map on $\bG$. Since 
$F^2$ is the standard Frobenius map with respect to ${\FF}_{q^2}$, we 
have $\bG^{F} \subseteq \GL_n({\FF}_{q^2})$. Now the restriction of 
$F_q$ to $\GL_n({\FF}_{q^2})$ is an automorphism of order $2$, which we 
denote by $A \mapsto \bar{A}$.  Then we obtain
\[\bG^{F}=\text{GU}_n({\FF}_q):=\{A \in \GL_n ({\FF}_{q^2})\mid 
\bar{A}^{\text{tr}}\cdot Q_n\cdot A=Q_n\},\]
the general unitary group with respect to the hermitian form defined
by~$Q_n$. The groups $\bB$ and $\bN$ are $F$-stable; furthermore,
$F$ induces an automorphism of $\bW \cong \Sym_n$ which we denote by the
same symbol; that automorphism is given by conjugation with the permutation
$w_0 \in \Sym_n$ whose matrix is $Q_n$. Then $\bB^F$ and $\bN^F$ form a 
split BN-pair of characteristic~$p$ in $\text{GU}_n(\FF_q)$, with Weyl group 
$W\cong \bW^F$; see \cite[4.2.6]{myag}. We have that $W$ is a Coxeter group 
of type $B_{n/2}$ (if $n$ is even) or of type $B_{(n-1)/2}$ (if $n$ is odd).
\end{exmp}


Now let $k$ be an algebraically closed field such that $p$ is invertible 
in~$k$. We denote by $kG\modA$ the category of finite-dimensional (left) 
$kG$-modules. Let $\Irr_k(G)$ be the set of simple $kG$-modules (up to 
isomorphism). 
Now let $I \subseteq S$. As we have noted above, the Levi subgroup $L_I$
is like $G$, that is, a group with a split BN-pair of characteristic~$p$.
We define functors
\[ R_I^S \colon kL_I\modA \rightarrow kG\modA \quad\mbox{and}\quad
\Rst_I^S \colon kG\modA \rightarrow kL_I\modA\]
as follows. Let $X \in kL_I\modA$. Then we can regard $X$ as a 
$kP_I$-module via the canonical map $P_I \rightarrow L_I$ with kernel $U_I$;
denote that $kP_I$-module by~$\tX$. Then we set $R_I^S(X)=\Ind_{P_I}^G(\tX)$
(where $\Ind$ denotes the usual induction of modules from subgroups). 

Conversely, let $Y \in kG\modA$. Then, since $U_I$ is normalized by $L_I$, 
the fixed point set $\Fix_{U_I}(Y)$ is naturally a $kL_I$-module, which we 
denote $\Rst_I^S(Y)$. Since $p$ is invertible in $k$, we have the
following alternative description. Let 
\[e_I=\frac{1}{|U_I|}\sum_{u \in U_I} u \in kG;\] 
i.e., $e_I$ is the idempotent corresponding to the trivial $kU_I$-module. 
Then the map $\Fix_{U_I}(Y) \rightarrow e_IY$, $y \mapsto e_Iy$, is an 
isomorphism of $kL_I$-modules. 

The functors $R_I^S$ and $\Rst_I^S$ have functorial properties similar
to the usual induction and restriction: transitivity, adjointness with 
respect to $\Hom$ and a Mackey formula. See the survey in \cite{mycharl}
for further details and references. 

\begin{defn} \label{MGcuspdef} Let $Y \in kG\modA$. We say that $Y$ is 
{\em cuspidal} if $\Rst_I^S(Y)=\{0\}$ for every proper subset $I \subset S$. 
Thus, $Y$ is cuspidal if and only if $\Fix_{U_I}(Y)=\{0\}$ for every proper 
$I \subset S$. We set
\[\Irr_k^\circ(G):=\{ Y \in \Irr_k(G) \mid \mbox{$Y$ is cuspidal}\}.\] 
Similar notations are used for Levi subgroups of $G$.  
\end{defn}

Let $\fC_G$ be the set of all pairs $(I,X)$ where $I\subseteq S$ and
$X \in \Irr_k^\circ(L_I)$. Given two such pairs $(I,X)$ and $(J,Y)$, 
we write $(I,X) \approx (J,Y)$ if there exists some $w \in W$ such 
that $w^{-1}Iw=J$ and ${^{\dot{w}}}X \cong_{kL_J} Y$. (Here,  
${^{\dot{w}}}X$ is the $kL_J$-module with the same underlying vector
space $X$ but where the action is defined by composition with the
group isomorphism $L_I \rightarrow L_J$ given by conjugation with 
$\dot{w}$.) This defines an equivalence relation on $\fC_G$. 

Now let us fix $Y \in \Irr_k(G)$. Then, by the transitivity of
Harish--Chandra induction and restriction, there exists a pair $(I,X) 
\in \fC_G$ such that the following two conditions are satisfied:
\begin{itemize}
\item[(HC1)] $\Rst_I^S(Y) \neq \{0\}$ and $I \subseteq S$ is minimal
with property;
\item[(HC2)] $X$ is a composition factor of $\Rst_I^S(Y)$. 
\end{itemize}
It is known that, if $(I',X')\in \fC_G$ is another pair satisfying (HC1) 
and (HC2), then $(I,X) \approx (I',X')$; see Hiss \cite{Hiss1} for the
case where $k$ has positive characteristic. This leads to the following 
definition. Let $(I,X) \in \fC_G$. Then $\Irr_k(G,(I,X))$ 
is defined to be the set of all $Y \in \Irr_k(G)$ such that (HC1) and (HC2)
hold. The set $\Irr_k(G,(I,X))$ is called the {\em Harish--Chandra 
series} defined by $(I,X)$. We have 
\[ \Irr_k(G)=\coprod_{(I,X) \in \fC_G/\approx} \Irr_k(G,(I,X)).\]
Thus, up to this point, everything formally works as in the classical 
theory of Harish--Chandra series in characteristic zero; see, for example, 
Carter \cite[Chap.~9]{Ca2} or Digne--Michel \cite[Chap.~6]{DiMi}. In 
this  case, by the adjointness between Harish--Chandra induction and 
restriction, the series  $\Irr_k(G,(I,X))$ can simply be characterized 
as the set of all simple $kG$-modules which occur as composition factors 
in $R_I^S(X)$. In the general case where $k$ is no longer assumed to have 
characteristic zero, we have the following result:

\begin{thm}[Hiss \cite{Hiss1}] \label{hiss1} Let $Y \in \Irr_k(G)$
and $(I,X) \in \fC_G$. Then the following three conditions are equivalent:
\begin{itemize}
\item[(a)] $Y \in \Irr_k(G,(I,X))$.
\item[(b)] $Y$ is isomorphic to a submodule of $R_I^S(X)$.
\item[(c)] $Y$ is isomorphic to a quotient of $R_I^S(X)$.
\end{itemize}
\end{thm}

The above formulation of Hiss' result relies on the following fundamental
property of Harish--Chandra induction:
\[ R_I^S(X) \cong_{kG} R_J^S(Y) \qquad \mbox{if $(I,X) \approx (J,Y)$
in $\fC_G$}.\]
This was proved independently by Dipper--Du \cite{DipDu} and 
Howlett--Lehrer \cite{HL2}.

\begin{exmp} \label{MGhcgl2} Let $G=\GL_2(\FF_q)$ where $q$ is a power of the
prime~$p$, with the BN-pair specified in Example~\ref{gln1}. Consider the pair
$(\varnothing,k_H) \in \fC_G$, where $k_H$ denotes the trivial $kH$-module. 
Then $M:=R_{\varnothing}^S(k_H)$ is the permutation module on the cosets of
$B$; in particular, we have 
\[ \dim M=[G:B]=q+1.\]
There are two essentially different cases:
\begin{itemize}
\item The characteristic of~$k$ does not divide $[G:B]$. Then
\[ M \cong_{kG} k_G\oplus Y \quad \mbox{where $Y \in \Irr_k(G)$, 
$\dim Y=q$};\]
in fact, $Y$ is the Steinberg module. Thus we have
\[ \Irr_k(G,(\varnothing,k_H))=\{k_G,Y\}.\]
\item The characteristic of $k$ divides $[G:B]$. Then $M$ is
indecomposable with submodules $\{0\} \subset V \subset V' \subset M$
such that 
\[ V \cong_{kG} M/V' \cong_{kG} k_G \quad \mbox{and} \quad V'/V \in
\Irr_k^\circ(G).\]
Hence, in this case, we have $\Irr_k(G,(\varnothing,k_H))=\{k_G\}$. 
\end{itemize}
\end{exmp}

Now let us return to the general situation and fix a pair $(I,X) \in \fC_G$.
In order to obtain more information about the corresponding Harish--Chandra 
series $\Irr_k(G,(I,X))$, we consider 
\[ \cH:=\End_{kG}\big(R_I^S(X)\big)^\circ,\]
the opposite algebra of the algebra of $kG$-endomorphisms of $R_I^S(X)$.
This is a finite-dimensional algebra over $k$, called a {\em Hecke algebra}.
We denote by $\cH\modA$ the category of finite-dimensional (left) 
$\cH$-modules and by $\Irr(\cH)$ the set of simple $\cH$-modules, up to 
isomorphism. We have a functor 
\[ \fF \colon kG\modA \rightarrow \cH\modA, \qquad
V \mapsto \Hom_{kG}(R_I^S(X),V),\]
where the action of $\varphi \in \cH$ on $\fF(V)$ is given by 
$\varphi.f=f \circ \varphi$ for $f \in \fF(V)$, and where $\fF$ sends a 
homomorphism of $kG$-modules $\rho \colon V \rightarrow V'$ to the 
homomorphism of $\cH$-modules $\fF(\rho)=\rho_* \colon \fF(V) 
\rightarrow \fF(V')$, $f \mapsto \rho \circ f$. 

By general principles (``Fitting correspondence''), $\fF$
induces a bijection between the isomorphism classes of indecomposable 
direct summands of $R_I^S(X)$ and the isomorphism classes of projective 
indecomposable $\cH$-modules; note that the latter are in natural bijection 
with $\Irr(\cH)$. Hence, if $k$ has characteristic zero, then $\fF$ certainly
induces a bijection between $\Irr_k(G,(I,X))$ and $\Irr(\cH)$. This 
statement remains valid in the general case, but the proof is much harder:
one has to establish some properties of the indecomposable direct 
summands of $R_I^S(X)$. The precise result is as follows. 

\begin{thm}[Geck--Hiss--Malle \cite{ghm2}, Geck--Hiss \cite{lymgh}] 
\label{ghm2} In the above setting, $\cH$ is a symmetric algebra
and we have a bijection 
\[ \Irr_k(G,(I,X)) \stackrel{\sim}{\longrightarrow} 
\Irr(\cH),\qquad Y \mapsto \fF(Y).\]
Every indecomposable direct summand of $R_I^S(X)$ has a unique simple 
submodule and a unique simple quotient, and these are isomorphic to 
each other.
\end{thm}

The fact that $\cH$ is symmetric means that there exists a trace function
$\tau \colon \cH \rightarrow k$ such that the associated bilinear form
\[ \cH \times \cH \rightarrow k, \qquad (h,h') \mapsto \tau(hh'),\]
is non-degenerate. The existence of $\tau$ is proved in \cite{lymgh}, 
using the construction of a suitable basis of $\cH$ in \cite{ghm2}; see
the remarks further below. Thus, $R_I^S(X)$ is a module whose endomorphism 
algebra is symmetric and which satisfies the two equivalent conditions 
(b) and (c) in Theorem~\ref{hiss1}. Then the statements in 
Theorem~\ref{ghm2} follow from a general argument combining ideas of 
Green, Cabanes and Linckelmann; see \cite[\S 2 and \S 3]{mycharl} 
for further details. It is interesting to note that these arguments first
appeared in the representation theory of finite groups with a split
BN-pair of characteristic~$p$, where the base field of the representations 
also has characteristic~$p$.

Finally, we describe the structure of $\cH$ in some more detail. For
this purpose, we recall the general definition of an Iwahori--Hecke
algebra associated with a Coxeter group; see \cite{ourbuch} for the 
general theory. Let $W_1$ be a finite Coxeter group. Thus, $W_1$ has a 
presentation with generators $S_1 \subseteq W_1$ and defining relations 
of the form: 
\begin{itemize}
\item $s^2=1$ for all $s\in S_1$;
\item $(st)^{m(s,t)}=1$ for $s\neq t$ in $S_1$, where $m(s,t)$ denotes the
order of $st$.
\end{itemize}
Let $l\colon W_1 \rightarrow \NN$ be the corresponding length function,
where $\NN=\{0,1,2,\ldots\}$. Let $R$ be any commutative ring with $1$,
and let $R^\times$ be the group of multiplicative units. We say that a 
function
\[ \pi \colon W_1\rightarrow R^\times\]
is a parameter function if $\pi(ww')= \pi(w)\cdot \pi(w')$ whenever we have 
$l(ww')= l(w)+l(w')$ for $w,w'\in W_1$. Such a function is uniquely 
determined by the values $\pi(s)$, $s\in S_1$, subject only to the condition
that $\pi(s)=\pi(t)$ if $s,t\in S_1$ are conjugate in $W_1$. Let $H_1=H_R
(W_1,\pi)$ be the corresponding Iwahori--Hecke algebra over $R$ with 
parameters $\{\pi(s)\mid s\in S_1\}$. The algebra $H_1$ is free over $R$ 
with basis $\{T_w\mid w \in W_1\}$, and the multiplication is given by the 
rule
\[T_sT_w=\left\{\begin{array}{cl} T_{sw}&\quad \mbox{if $l(sw)>l(w)$},\\ 
\pi(s)T_{sw}+(\pi(s)-1)T_w & \quad \mbox{if $l(sw)<l(w)$},\end{array} 
\right.\]
where $s\in S_1$ and $w\in W_1$. We define a linear map $\tau \colon H_1
\rightarrow k$ by 
\[ \tau(T_1)=1 \qquad \mbox{and} \qquad \tau(T_w)=0 \quad 
\mbox{for $w\neq 1$}.\]
Then one can show that $\tau$ is a trace function and we have 
\[ \tau(T_wT_{w'})=\left\{\begin{array}{cl} \pi(w) & \qquad 
\mbox{if $w'=w^{-1}$},\\ 0 & \qquad \mbox{otherwise}.\end{array}\right.\]
This implies that $H$ is a symmetric algebra. 

Now let us return to the pair $(I,X) \in \fC_G$ and the corresponding
endomorphism algebra $\cH$. Let $\cN(I,X)$ be the stabilizer of $X$ in 
$\cN(I):=(\operatorname{N}_G(L_I)\cap N)L_I$. Then we have 
\[ \dim \cH=|\cW(I,X)| \qquad \mbox{where}\qquad \cW(I,X):=\cN(I,X)/L.\]
In fact, one can explicitly construct a basis $\{B_w\mid w\in \cW(I,X)\}$ 
of $\cH$, as in \cite{ghm2}; that construction generalizes the one by 
Howlett--Lehrer \cite{HL1} (which is concerned with the case where the
ground field has characteristic zero). Now, there is a semidirect product 
decomposition $\cW(I,X)=W_1 \rtimes \Omega$ where $W_1$ is a finite Coxeter 
group (with generating set $S_1$) and $\Omega$ is a finite group such that 
$\omega S_1 \omega^{-1}=S_1$ for all $\omega \in \Omega$. Let $\cH_1:=
\langle B_w \mid w\in W_1\rangle \subseteq \cH$. The multiplicative properties
of the basis $\{B_w\}$ show that $\cH_1$ is a subalgebra of $\cH$ and that
we have a direct sum decomposition
\[ \cH=\bigoplus_{\omega \in \Omega} \cH_\omega \qquad \mbox{where}\qquad
\cH_\omega:=H_1\cdot B_\omega=B_\omega \cdot \cH_1;\]
furthermore, each $B_\omega$ ($\omega\in \Omega$) is invertible and we 
have $\cH_{\omega \omega'}=\cH_{\omega} \cdot B_{\omega'}$ for all 
$\omega, \omega'\in \Omega$. Hence, as remarked in \cite[Prop.~2.4]{lymgh}, 
the family of subspaces $\{\cH_\omega \mid \omega \in \Omega\}$ is an 
$\Omega$-graded Clifford system in $\cH$, in the sense of 
\cite[Def.~11.12]{CR2}.

\begin{thm}[Howlett--Lehrer \cite{HL1} (\text{char}(k)=0) and 
Geck--Hiss--Malle \cite{ghm2}] \label{howleh}  In the above setting, 
assume that $X$ can be extended to $\cN(I,X)$. Then we have $\cH_1 
\cong H_k(W_1,\pi)$ for a suitable parameter function $\pi \colon W_1
\rightarrow k^\times$.
\end{thm}

In many cases, we have $\Omega=\{1\}$ and the hypothesis of the
above result can be seen to hold. Hence, in these cases, the algebra
$\cH_1$ is an Iwahori--Hecke algebra associated with a finite Coxeter group.
We will see some examples in the following section.

\begin{rem} \label{sigma} Assume that the hypotheses of Theorem~\ref{howleh}
are satisfied, where $k$ is an algebraic closure of the finite field with
$\ell$ elements. Thus, we have $\cH_1\cong H_k(W_1,\pi)$ where $\pi\colon W_1
\rightarrow k^\times$ is a parameter function. Since $W_1$ is finite, there
exists a finite subfield $k_0 \subseteq k$ such that $\pi(W) \subseteq k_0$.
Hence, since $k_0^\times$ is cyclic, there exist some $\xi\in k^\times$ 
and $d\geq 1$ such that 
\[ \xi^d=1 \quad \mbox{and}\quad \pi(s)=\xi^{L(s)} \quad 
\mbox{for all $s\in S_1$},\]
where $L \colon W_1 \rightarrow \NN$ is a {\em weight function} in the 
sense of Lusztig \cite{Lusztig03}, that is, we have $L(ww')=L(w)+L(w')$
whenever we have $l(ww')=l(w)+l(w')$ for $w,w'\in W_1$. 
\end{rem}

\begin{rem} \label{remDJ} Dipper and James extensively studied the case 
where $G=\GL_n(\FF_q)$ and $k$ has positive characteristic; see 
\cite{QuotHom}, \cite{Ja0} and the references there. They obtained a 
complete classification of all cuspidal simple modules and a 
parametrization of the simple modules of the corresponding Hecke algebras 
in this case. An outstanding role in this context is played by the 
introduction of the $q$-Schur algebra \cite{DJ1}.  A special feature of 
$\GL_n(\FF_q)$ is the fact that all cuspidal simple modules in positive 
characteristic have a ``reduction stable'' lift to characteristic zero. 
The fact that this is no longer true for groups of other types is the 
source of substancial, new complications. 
\end{rem}

\section{Unipotent blocks} \label{MGsec3}

We now wish to restrict our attention to the ``unipotent'' modules
of $G$, in the sense of Deligne--Lusztig \cite{DeLu}. For this purpose, 
let us assume from now on that $G$ arises as the fixed point set of a 
connected reductive algebraic group $\bG$ over $\overline{\FF}_p$ under 
a Frobenius map $F \colon \bG \rightarrow \bG$, with respect to some
${\FF}_q$-rational structure on $G$ where $q$ is a power of $p$. 
(See Example~\ref{gln1}.) Thus, we have  
\[G=\bG(\FF_q)=\bG^F=\{g \in \bG \mid F(g)=g\}.\] 
Let $\bB \subseteq \bG$ be an $F$-stable Borel subgroup and $\bT_0 
\subseteq \bG$ be an $F$-stable maximal torus which is contained in~$\bB$. 
Let $\bW$ be the Weyl group of $\bG$ with respect to~$\bT_0$; the Frobenius 
map $F$ induces an automorphism of $\bW$ which we denote by the same symbol. 
Then $G$ is a finite group with a split $BN$-pair of characteristic~$p$ and 
Weyl group~$W$ where 
\[B:=\bB^F,\qquad N:=\operatorname{N}_{\bG}(\bT_0)^F, \qquad W \cong\bW^F;\]
see \cite[\S 2.9]{Ca2}. 

We shall need some results about the modules of $G$ over the field $\CC$. 
Given $V \in \CC G\modA$, the corresponding character is the function
$\chi_V \colon G \rightarrow \CC$ which sends $g \in G$ to the trace of
$g$ acting on $V$. It is well-known that $V \cong_{\CC G} V'$ (where $V,V' 
\in \CC G\modA$) if and only if $\chi_V=\chi_{V'}$. We set 
\[ G^\wedge=\{\chi \colon G \rightarrow \CC \mid \chi=\chi_V \mbox{ for 
some $V \in \Irr_{\CC}(G)$}\}.\]
A function $f \colon G \rightarrow \CC$ will be called a {\em virtual
character} if $f$ is an integral linear combination of $G^\wedge$. In this 
case, we denote by $\langle f:\chi\rangle$ the coefficient of $\chi$ in the
expansion of $f$ in terms of $G^\wedge$, that is, we have 
\[ f=\sum_{\chi \in G^\wedge} \langle f:\chi\rangle\, \chi.\]

Recall that the conjugates of $\bT_0$ in $\bG$ are called the {\em maximal
tori} of $\bG$. Let $g \in \bG$ and consider $\bT:=g\bT_0g^{-1}$. This 
group will be $F$-stable if and only if $g^{-1}F(g) \in
\operatorname{N}_{\bG}(\bT_0)$.  In this case, $F$ restricts to a Frobenius 
map on $\bT$ and, by taking fixed points, we obtain a finite group 
$T:=\bT^F \subseteq G$. Note that $T$ is abelian of order prime to~$p$. 
Consequently, the set $T^\wedge$ also is an abelian group, for the pointwise 
product of characters. 

Let $\cS$ be the set of all pairs $(\bT,\theta)$ where $\bT \subseteq \bG$ 
is an $F$-stable maximal torus and $\theta \in T^\wedge$. Deligne and
Lusztig \cite{DeLu} have associated with each pair $(\bT,\theta) \in \cS$
a virtual character $R_{T,\theta}$ of $G$. The construction uses deep
results from algebraic geometry; see Lusztig \cite{Lusztig77} and 
Carter \cite{Ca2} for more details. Here are some properties:
\begin{itemize}
\item if $\bT=\bT_0$, then $R_{T,\theta}$ is the character of 
$R_\varnothing^S(\theta)$;
\item $R_{T,\theta}(1)=\eta_T \,[G:T]_{p'}$, where $\eta_T=\pm 1$;
\item given any $\chi \in G^\wedge$, we have $\langle R_{T,\theta}:
\chi\rangle \neq 0$ for some $(\bT,\theta)\in \cS$.
\end{itemize}
Now let $\cS_1$ be the set of all pairs $(\bT,\theta) \in \cS$ such that
$\theta$ is the principal character $1_T$ of $T$. Then the characters 
\[\Uch(G):=\{ \chi \in G^\wedge \mid \langle R_{T,\theta}:\chi \rangle 
\neq 0 \mbox{ for some $(T,\theta) \in \cS_1$}\}\]
are called the {\em unipotent characters} of $G$.  Lusztig \cite{Lu4}, 
\cite{Lu5} obtained a classification of the unipotent characters and 
formulae for the corresponding degrees; furthermore, he showed that the 
classification of all irreducible characters can be reduced to the case of 
unipotent characters, in terms of a ``Jordan decomposition'' of characters. 
The unipotent characters only depend on the Weyl group of $\bG$; more 
precicely:

\begin{thm}[Lusztig \cite{Lu4}] \label{unip} There exists a finite 
set $\bLambda$ and polynomials $D_\lambda \in {\QQ}[X]$ ($\lambda
\in \bLambda$) such that there is a bijection 
\[ \bLambda \stackrel{\sim}{\longrightarrow} \Uch(G),\qquad \lambda
\mapsto \chi_\lambda,\]
where $\chi_\lambda(1)=D_\lambda(q)$ for all $\lambda \in \bLambda$.
The set $\bLambda$ and the polynomials $\{D_\lambda\}$ only depend on 
the Coxeter group $\bW$ and the induced map $F \colon \bW \rightarrow \bW$.
\end{thm}

Explicit tables of the polynomials $D_\lambda$ can be found in the
appendix of \cite{Lu4}. Thus, the unipotent characters play an essential 
role in the whole theory. An analogous result in the ``modular situation'' 
has been obtained by Bonnaf\'e--Rouquier \cite{BoRo}. The definition of
unipotent $kG$-modules, where $k$ has positive characteristic $\ell\neq p$,
is based on the following result.

\begin{thm}[Brou\'e--Michel \cite{BrMi}] \label{brmi} Let $\ell$ be a prime
number $\neq p$ and set 
\[ b_\ell(g):=\frac{1}{|G|}\frac{1}{|G|_p} \sum_{(\bT,\theta)}
\eta_T \, R_{T,\theta}(g^{-1}) \in \CC \qquad \mbox{for any $g \in G$},\]
where the sum runs over all pairs $(\bT,\theta) \in \cS$ such that
$\theta$ is of order a power of $\ell$ in $T^\wedge$. Then we have
$b_\ell(g) \in \ZZ_{(\ell)}$ for any $g \in G$, where $\ZZ_{(\ell)}$ 
denotes the localization of $\ZZ$ in the prime ideal generated by $\ell$.
The element 
\[ \beta_\ell:=\sum_{g \in G} \overline{b}_\ell(g)\,g\in {\FF}_\ell G\]
is a central idempotent, where the bar denotes the canonical map
$\ZZ_{(\ell)} \rightarrow \FF_\ell$.
\end{thm}

Let us indicate how the formula for $b_\ell(g)$ comes about. For any 
$\chi\in G^\wedge$, let $e_\chi$ be the corresponding central primitive 
idempotent in $\CC G$.  We have
\[ e_\chi=\frac{1}{|G|} \sum_{g \in G} \chi(1)\chi(g^{-1})\, g.\]
Now let $\cS_\ell$ be the set of all pairs $(\bT,\theta) \in \cS$ such that
$\theta$ is of order a power of $\ell$ in $T^\wedge$. We set 
\[ e_\ell(g):=\frac{1}{|G|} \sum_{\chi \in \cB_\ell} \chi(1)\chi(g^{-1})
\qquad \mbox{for any $g \in G$},\]
where $\cB_\ell$ is the set of all $\chi \in G^\wedge$ such that 
$\langle R_{T,\theta}:\chi\rangle \neq 0$ for some $(\bT,\theta)\in
\cS_\ell$. Hence, $e_\ell:=\sum_{g \in G} e_\ell(g)\,g$ is a central 
idempotent in ${\CC}G$. We want to show that $e_\ell(g)=b_\ell(g)$ for 
all $g \in G$. For this purpose, we consider the character 
$\chi_{\text{reg}}$ of the regular representation of $G$. We have
\[ \chi_{\text{reg}}=\sum_{\chi \in G^\wedge} \chi(1)\chi=
\frac{1}{|G|_p} \sum_{(\bT,\theta) \in \cS} 
\eta_T R_{T,\theta},\]
where the first equality is well-known and the second equality holds by
\cite[Cor.~7.5.6]{Ca2}. Using the partition of $G^\wedge$ into geometric 
conjugacy classes (see \cite[\S 7.3]{Ca2}), we conclude that 
\[ \sum_{\chi \in \cB_\ell} \chi(1)\chi=\frac{1}{|G|_p} 
\sum_{(\bT,\theta) \in \cS_\ell} \eta_T R_{T,\theta},\]
Thus, we have $e_\ell(g)=b_\ell(g)$ for all $g \in G$, as claimed.
In particular, we have shown that $b_\ell:=\sum_{g \in G} b_\ell(g)g$
is a central idempotent in $\CC G$. 

Finally, using the character formula for $R_{T,\theta}$ (see 
\cite[7.2.8]{Ca2}) one easily shows that $b_\ell(g)\in \QQ$ for all $g\in G$. 
Now Brou\'e--Michel prove that the coefficients actually lie in $\ZZ_{(\ell)}$. 
Hence, we can reduce $b_\ell$ modulo $\ell$ and obtain a central 
idempotent $\beta_\ell \in \FF_\ell G$, as desired. \hfill $\Box$

\medskip
Now let $k$ be an algebraically closed field of characteristic $\ell$,
where $\ell$ is a prime $\neq p$. Then $\beta_\ell \in {\FF}_\ell G
\subseteq kG$ is a central idempotent and, hence, we have a direct sum 
decomposition 
\[ kG=\beta_\ell\, kG \oplus (1-\beta_\ell)\,kG\]
where both $\beta_\ell\, kG$ and $(1-\beta_\ell)\,kG$ are two-sided ideals.
Correspondingly, we have a decomposition of $kG\modA$ into those modules 
on which $\beta_\ell$ acts as the identity on the one hand, and those 
modules on which $\beta_\ell$ acts as zero on the other hand. We set 
\[ \Uch_\ell(G):=\{ Y \in \Irr_k(G) \mid \beta_\ell.Y=Y\}\]
and 
\[\Uch_\ell^\circ(G):=\Uch_\ell(G)\cap \Irr_k^\circ(G).\]
Note that, if we formally set $\ell=1$ and take $k=\CC$, then $\Uch_{1}(G)$ 
indeed is the set of simple $\CC G$-modules whose character is unipotent. 
Thus, the above definitions generalize the definition of unipotent characters
to $kG$-modules where the characteristic of $k$ is a prime $\ell\neq p$.
Now we can state:

\begin{cor}[Hiss \protect{\cite{Hiss1}}] \label{hcuni} Let $k$ be an 
algebraically closed field of characteristic $\ell$, where $\ell$ is a 
prime $\neq p$. Then we have 
\[\Uch_{\ell}(G)=\coprod_{(I,X)} \Irr_k(G,(I,X)),\]
where the union runs over all pairs $(I,X) \in \fC_G$ such that 
$X\in\Uch_\ell^\circ(L_I)$.
\end{cor}

This is based on a compatiblity of Harish-Chandra induction with the 
operator $R_{T,\theta}$; see Lusztig \cite[Cor.~6]{Lu1}.

\begin{thm}[Geck--Hiss \cite{bs1}] \label{MGexistbs} Assume that the 
center of $\bG$ is connected and that $\ell$ is good for $\bG$. Then we
have $|\Uch_\ell(G)|=|\Uch(G)|$.
\end{thm}

We say that the prime $\ell$ is {\em good} for $\bG$ if it is good for each 
simple factor involved in~$\bG$; the conditions for the various simple 
types are as follows.
\begin{center} $\begin{array}{rl} A_n: & \mbox{no condition}, \\
B_n, C_n, D_n: & \ell \neq 2, \\
G_2, F_4, E_6, E_7: &  \ell \neq 2,3, \\
E_8: & \ell \neq 2,3,5.  \end{array}$
\end{center}
If $\ell$ is not good for $G$, we have $|\Uch_\ell(G)|\neq |\Uch(G)|$.
For further information on the cardinalities of the sets $|\Uch_\ell(G)|$
in this case, see \cite[6.6]{lymgh} and \cite[\S 4.1]{DGHM}.

\begin{conj}[James for type $\GL_n$] \label{ourconj} Recall the finite
set $\bLambda$ and the polynomials $D_\lambda$ in Theorem~\ref{unip}.
Now let $e\geq 1$. Then there exist polynomials $D_\lambda^{(e)} \in 
{\QQ}[X]$ ($\lambda \in \bLambda$) such that the following holds: for any 
prime $\ell\neq p$ such that $\ell$ does not divide $|\bW|$ and $e$ is the
multiplicative order of $q$ modulo $\ell$, we have a bijection
\[ \bLambda \stackrel{\sim}{\longrightarrow} \Uch_\ell(G), \qquad
\lambda \mapsto X^\lambda,\]
such that $\dim X^\lambda=D_\lambda^{(e)}(q)$ for any $\lambda \in \bLambda$.
The polynomials $D_\lambda$ and $D_\lambda^{(e)}$ have the same degree
and the same leading coefficient.
\end{conj}

The above conjecture is known to hold in some examples of small rank
by explicit verification, most notably for $G=\GL_n({\FF}_q)$ where 
$n\leq 10$; see James \cite{Ja}. For the example $G=\GU_3({\FF}_q)$; 
see Okuyami--Waki \cite{OkWa}.

Let us now turn to the discussion of some examples of cuspidal unipotent 
modules and the corresponding Hecke algebras.

\begin{exmp} \label{iwa} Assume that $G=\bG^F$ and that the induced 
map $F \colon \bW \rightarrow \bW$ is the identity. The pair $(\varnothing,
k_H)$ is cuspidal and $k_H$ is unipotent. Furthermore, 
$R_\varnothing^S(k_H)$ is nothing but the permutation module of $G$ on 
the cosets of $B$. This is the case originally considered by Iwahori;
see \cite[\S 8.4]{ourbuch} and the references there. We have
\[ \cH\cong H_k(W,\bar{q})\]
where the notation indicates that the parameter function $\pi\colon 
W \rightarrow k$ is given by $\pi(s)=\bar{q}:=q\cdot 1_k\in k$ for all
$s\in S$.
\end{exmp}

\begin{exmp} \label{gln2} Let $G=\GL_n(\FF_q)$. Then we have 
\[ |\Uch_\ell(G)|=|\Uch(G)|=\mbox{number of partitions of $n$}.\]
Now let $e$ be the multiplicative order of $q$ modulo~$\ell$. Then we have 
\begin{equation*}
 |\Uch_\ell^\circ(\GL_n(\FF_q))|=\left\{\begin{array}{cl} 1 & \quad
\mbox{if $n=1$ or $n=e\ell^i$ for some $i \geq 0$},\\ 0 & \quad 
\mbox{otherwise}.  \end{array}\right.\tag{$*$}
\end{equation*}
There are several proofs for this result: see James \cite{Ja0}, Dipper
\cite{QuotHom}, Geck--Hiss--Malle \cite[\S 7]{ghm1} or \cite[\S 2]{ghm2}.
The unique $X^\circ \in \Uch_{\ell}^\circ(G)$ (if it exists) can be lifted
to a cuspidal module in characteristic zero (which is not unipotent); we have
\[ \dim X^\circ=(q^{n-1}-1)(q^{n-2}-1) \cdots (q^2-1)(q-1).\]
Let us now consider the Hecke algebras arising in $G$. Let $\lambda$ be a
partition of $n$, with non-zero parts $\lambda_1 \geq \lambda_2 \geq \cdots
\geq \lambda_r>0$. Then we have a corresponding Levi subgroup $L_I\subseteq 
G$ such that 
\[ L_I \cong \GL_{\lambda_1}({\FF}_q) \times  \cdots  \times
\GL_{\lambda_r}({\FF}_q);\] 
furthermore, any Levi subgroup of $G$ arises in this way (up to $\approx$).
Now, a simple module $X \in \Irr_k(L_I)$ is isomorphic to a tensor product 
of simple modules for the various factors in the above direct product, 
and $X$ is cuspidal unipotent if and only if every factor has this property. 
Let us assume that $X \in \Uch_\ell^\circ(L_I)$. Then, by ($*$), $\lambda$ 
can only have parts equal to $1$ or $e\ell^i$ for some $i\geq 0$. Let
\begin{align*}
m_{-1} &:= \mbox{multiplicity of $1$ as a part of $\lambda$},\\
m_{i} &:= \mbox{multiplicity of $e\ell^i$ as a part of $\lambda$} \qquad
(i=0,1,2,\ldots).
\end{align*}
(Here, we understand that $m_{-1}=0$ if $e=1$.) We have $\cN(I,X)=\cN(I)$ 
and $X$ can be extended to $\cN(I)$; furthermore, $\cW(I,X)$ is a Coxeter 
group $W_1$ isomorphic to 
\[ \Sym_{m_{-1}} \times \Sym_{m_0} \times \Sym_{m_1} \times \Sym_{m_2} 
\times \cdots.\]
(Thus, $\Omega=\{1\}$ in this case.) Finally, the parameter function 
$\pi \colon W_1 \rightarrow k^\times$ is given as follows. We have $\pi(s)
=\bar{q}$ for a generator $s$ in the $\Sym_{m_{-1}}$-factor, and 
$\pi(s)=1$ otherwise. See Dipper \cite[Part~II, \S 4]{QuotHom} for more 
details. The simple $\cH$-modules are classified by Dipper--James 
\cite{DJ1}.
\end{exmp}

\begin{exmp} \label{gln3} Let $G=\GU_n(\FF_q)$. Then we have again
\[ |\Uch_\ell(G)|=|\Uch(G)|=\mbox{number of partitions of $n$}.\]
Furthermore, we have
\[ |\Uch_\ell^\circ(\GU_n(\FF_q))| \leq p_2(n),\]
where $p_2(n)$ denotes the number of all partitions of $n$ with distinct 
parts; see Geck--Hiss--Malle \cite[Prop.~6.2 and Prop.~6.8]{ghm1}. There are 
examples where equality holds (see case (b) below); in general, the exact 
number of cuspidal unipotent simple $kG$-modules is not known! Let us now
consider the Hecke algebras arising in $G$.

Let us write $n=m'+2m$ where $m,m' \in \NN$, and let $\lambda$ be a 
partition of $m$. For any $i \in \{1,\ldots,m\}$, let $n_i\geq 0$ be the
multiplicity of $i$ as a part of $\lambda$. Then we have a corresponding 
Levi subgroup $L_I\subseteq G$ such that 
\[ L_I \cong \GU_{m'}({\FF}_q) \times \prod_{i=1}^m 
\bigl(\underbrace{\GL_i({\FF}_{q^2}) \times \cdots \times 
\GL_i({\FF}_{q^2})}_{\text{$n_i$ factors}}\bigr);\]
furthermore, any Levi subgroup of $G$ arises in this way (up to $\approx$).
Let us assume that $X \in \Uch_\ell^\circ(L_I)$. Then, by 
Example~\ref{gln2}, $\lambda$ can only have parts equal to $1$, $e$, 
$e\ell$, $e\ell^2$, $\ldots$, where $e$ is the multiplicative order of 
$q^2$ modulo $\ell$. As in the previous example, let 
\begin{align*}
m_{-1} &:= \mbox{multiplicity of $1$ as a part of $\lambda$},\\
m_{i} &:= \mbox{multiplicity of $e\ell^i$ as a part of $\lambda$} \qquad
(i=0,1,2,\ldots).
\end{align*}
(Here, we understand that $m_{-1}=0$ if $e=1$.) By \cite[Prop.~4.3]{ghm2}, 
we have $\cN(I,X)=\cN(I)$ and $X$ can be extended to $\cN(I)$; 
furthermore, $\cW(I,X)$ is a Coxeter group $W_1$ of type 
\[ B_{m_{-1}} \times B_{m_0} \times B_{m_1} \times B_{m_2} \times \ldots.\]
(Thus, $\Omega=\{1\}$ in this case.) Finally, by \cite[Prop.~4.4]{ghm2},
the parameter function $\pi\colon W_1 \rightarrow k^\times$ is given by 
\begin{center}
\begin{picture}(300,70)
\put(   0,   7){$B_{m_i}$:}
\put(  50,   7){\circle{7}}
\put(  90,   7){\circle{7}}
\put( 130,   7){\circle{7}}
\put( 250,   7){\circle{7}}
\put( 290,   7){\circle{7}}
\put( 53.5,  9){\line(1,0){ 33}}
\put( 53.5,  5){\line(1,0){ 33}}
\put(93.5,   7){\line(1,0){ 33}}
\put(133.5,  7){\line(1,0){ 33}}
\put( 179,   7){$\ldots$}
\put(213.5,  7){\line(1,0){ 33}}
\put(253.5,  7){\line(1,0){ 33}}
\put(   40, 17){$p_1'(i)$}
\put(   87, 17){$1$}
\put(  127, 17){$1$}
\put(  247, 17){$1$}
\put(  287, 17){$1$}
\put(   0,  43){$B_{m_{-1}}$:}
\put(  50,  43){\circle{7}}
\put(  90,  43){\circle{7}}
\put( 130,  43){\circle{7}}
\put( 250,  43){\circle{7}}
\put( 290,  43){\circle{7}}
\put( 53.5, 45){\line(1,0){ 33}}
\put( 53.5, 41){\line(1,0){ 33}}
\put(93.5,  43){\line(1,0){ 33}}
\put(133.5, 43){\line(1,0){ 33}}
\put( 179,  43){$\ldots$}
\put(213.5, 43){\line(1,0){ 33}}
\put(253.5, 43){\line(1,0){ 33}}
\put(   44, 54){$p_1$}
\put(   87, 54){$\bar{q}^2$}
\put(  127, 54){$\bar{q}^2$}
\put(  247, 54){$\bar{q}^2$}
\put(  287, 54){$\bar{q}^2$}
\end{picture}
\end{center}
for $i=0,1,2,\ldots$, where $p_1,p_1(i)'\in k^\times$ and $\bar{q}$ denotes
the image of $q$ in $k$. The parameters $p_1,p_1'(i)$ are only known in 
special cases. Let $d$ be the multiplicative order of $-q$ modulo $\ell$. 
Then
\[ e=\left\{\begin{array}{cl} d/2 & \quad \mbox{if $d$ is even},\\
d & \quad \mbox{if $d$ is odd}.\end{array}\right.\]
The following distinction between odd and even values for $d$ already occurs
in the work of Fong--Srinivasan \cite{fs1}.

\medskip
(a) {\em Assume that $d=2e$}. Then $m'=t(t+1)/2$, $p_1=\bar{q}^{2t+1}$ 
and $p_1'(i)\neq -1$; see \cite[Lemma~4.9]{ghm2}. In 
this situation, the simple $\cH$-modules are classified by Dipper--James 
\cite{DJ2} (see Example~\ref{canbn} below). A counting argument then 
yields that 
\[ |\Uch_\ell^\circ(G)|=\left\{\begin{array}{cl} 1 & \quad \mbox{if 
$n=s(s+1)/2$ for some $s\geq 1$},\\ 0 & \quad \mbox{otherwise};
\end{array}\right.\]
see \cite[Theorem~4.11]{ghm2}. The unique $X^\circ \in \Uch_\ell^\circ(G)$ 
(if it exists) can even be lifted to a cuspidal unipotent character in 
characteristic zero. The dimension of $X^\circ$ is given by the polynomial
in \cite[p.~457]{Ca2}. This case is studied further by Gruber--Hiss 
\cite{GrHi}.

\medskip
(b) {\em Assume that $d=1$}, that is, $\bar{q}=-1$. Then $p_1=p_1'(i)=-1$; 
see \cite[Lemma~4.6]{ghm2} and \cite[Prop.~2.3.5]{Gruber}. In this 
situation, the simple $\cH$-modules are again classified by Dipper--James 
\cite{DJ2} (see Example~\ref{canbn} below). A counting argument then 
yields that 
\[ |\Uch_\ell^\circ(G)|=p_2(n);\]
see \cite[Theorem~4.12]{ghm2}.  The dimensions of these modules are 
not known.

\medskip
(c) {\em Assume that $d>1$ is odd}. Then $p_1'(i)=-1$ for all $i\geq 0$;
see \cite[Prop.~2.3.5]{Gruber}. As before, the simple modules of the
Hecke algebra corresponding to the $B_{m_i}$-factor ($i\geq 0$)
are classified by Dipper--James \cite{DJ2}. The analogous problem for 
the factor of type $B_{m_{-1}}$ has been solved by Ariki--Mathas \cite{ArMa} 
and Ariki \cite{Ar2}, in terms of so-called Kleshchev bipartitions. These 
bipartitions are defined in a recursive way. Recently, Jacon \cite{Jac0} 
obtained a whole family of different parametrizations by so-called 
FLOTW bipartitions, which have the advantage of being defined in a 
non-recursive way; furthermore, they can be adapted more accurately to
the values of the parameter function $\pi$. (We will discuss all this
in more detail in Sections~7 and~8.)

\medskip
\noindent \underline{\bf Open problem:} {\em Determine $p_1$ for $d>1$ 
odd.  Or even better, find a general practical method for determining 
the function $\pi\colon W_1 \rightarrow k^\times$}. 

\medskip
Some partial results on the computation of $\pi(s)$ are contained in 
\cite[\S 3]{ghm2}; Hiss (unpublished) actually has a conjecture about
$p_1$. Note that, in characteristic zero, the parameters $\pi(s)$ are 
completely known; see Lusztig \cite{Lusztig77}, Table~II, p.~35. In 
this case, we have $\pi(s)=q^{L(s)}$ for $s\in S_1$, where  $L\colon 
W_1\rightarrow \NN$ is a weight function such that $L(s)>0$ for 
all $s\in S_1$.
\end{exmp}

To summarize, the above results show that it would be very interesting to 
know a parametrization of the simple modules of an Iwahori--Hecke algebra
$H_k(W_1,\pi)$ where $\pi \colon W_1 \rightarrow k^\times$ has values in 
an algebraically closed field $k$ of characteristic $\ell>0$. In particular, 
it would be interesting to know in which way the parametrization depends on 
the function $\pi$. In this context, we may assume that $\pi(s)=\xi^{L(s)}$ 
for all $s\in S_1$, where $\xi \in k^\times$ has finite order and $L\colon 
W_1\rightarrow\NN$ is a weight function (see Remark~\ref{sigma}).

\section{Generic Iwahori--Hecke algebras and specializations} 
\label{MGsec4}

We now consider a finite Coxeter group $W$ with generating set $S$ and a 
weight function $L \colon W\rightarrow \NN$, without reference to a 
realization of $W$ in the framework of groups with a BN-pair. Our aim is 
to develop the representation theory of the associated Iwahori--Hecke 
algebras over various fields $k$. One of first decisive observations is 
the fact that these algebras can be defined ``generically'' over a 
polynomial ring, where the parameters are powers of the indeterminate. 
The precise definitions are as follows. (A general reference is 
\cite{ourbuch}.)

Let $A={\ZZ}[v,v^{-1}]$ be the ring of Laurent polynomials in an 
indeterminate~$v$. Then there exists an associative algebra $\bH=\bH_A(W,L)$
over $A$, which is free as an $A$-module with basis $\{T_w\mid w\in W\}$
such that the multiplication is given by 
\[T_sT_w=\left\{\begin{array}{cl} T_{sw}&\quad \mbox{if $l(sw)>l(w)$},\\ 
v^{2L(s)}T_{sw}+(v^{2L(s)}-1)T_w & \quad \mbox{if $l(sw)<l(w)$},\end{array} 
\right.\]
where $s\in S$ and $w\in W$. The fact that the parameters are even powers
of $v$ will play a role in connection with the construction of the
Kazhdan--Lusztig basis of $W$ (to be discussed in Section~5) and 
also in connection with the question of splitting fields. 

\begin{Par} {\bf Specialization.} \label{spec} Let $k$ be a field and 
$\xi \in k^\times$ be an element which has a square root in $k^\times$. 
Then there is a ring homomorphism $\theta \colon A \rightarrow k$ such 
that $\theta(v^2)=\xi$. Considering $k$ as an $A$-module via $\theta$, we 
set 
\[ \bH_{k,\xi}:=k \otimes_A \bH.\]
Thus, $\bH_{k,\xi}$ is an associative $k$-algebra with a basis $\{T_w \mid
w \in W\}$ such that 
\[T_sT_w=\left\{\begin{array}{cl} T_{sw}&\quad \mbox{if $l(sw)>l(w)$},\\ 
\xi^{L(s)}T_{sw}+(\xi^{L(s)}-1)T_w & \quad \mbox{if $l(sw)<l(w)$},
\end{array} \right.\]
where $s\in S$ and $w\in W$. This shows that the endomorphism algebras
arising in the context of Harish-Chandra series as in the previous sections
are obtained via specialization from generic Iwahori--Hecke algebras.
\end{Par}

In order to deal with non-semisimple specializations of $\bH$, we shall
need some fundamental results on the structure of $\bH_{K,v}$ where $K$ 
is the field of fractions of $A$ and $\theta \colon A \hookrightarrow K$ 
is the inclusion. For technical simplicity, we will assume from now on that 
\begin{center}
{\em $W$ is a finite Weyl group},
\end{center}
that is, the product of two generators of $W$ is $2$, $3$, $4$ or $6$. Then
it is known that every complex representation of $W$ can be realized over 
$\QQ$, that is, the group algebra $\QQ W$ is split semisimple 
(see \cite[6.3.8]{ourbuch}). Using the specialization $v \mapsto 1$ and
Tits' deformation argument, it is not hard to show that $\bH_{K',v}$ is 
split semisimple and abstractly isomorphic to the group algebra $K'W$, 
where $K'$ is a sufficiently large finite extension of $K$; see 
\cite[8.1.7]{ourbuch}. We have the following more precise result, which 
combines work of Benson--Curtis, Lusztig, Digne--Michel; see 
\cite[9.3.5]{ourbuch}:

\begin{thm} \label{gen1} The algebra $\bH_{k,v}$ is split semisimple and
abstractly isomorphic to the group algebra $KW$. 
\end{thm}

Now, an isomorphism $\bH_{K,v}\cong KW$ certainly induces a bijection between
$\Irr(\bH_{K,v})$ and $\Irr_{\QQ}(W)$. To describe this, it will be 
convenient to work with characters. As in the case of finite groups, the 
character of an $\bH_{K,v}$-module $V$ is the function $\chi_V \colon 
\bH_{K,v} \rightarrow K$ sending $h \in \bH_{K,v}$ to the trace of $h$ 
acting on $V$. Since $A$ is integrally closed in $K$, a general 
argument (see \cite[7.3.9]{ourbuch}) shows that $\chi_V(T_w) \in A$ for 
all $w \in W$. Let $\bH_{K,v}^\wedge$ be the set of irreducible characters 
of $\bH_{K,v}$. Once Theorem~\ref{gen1} is establised, {\em Tits' 
Deformation Theorem} (see \cite[7.4.6]{ourbuch}) yields the following result.

\begin{thm} \label{gen2} There is a bijection $W^\wedge 
\stackrel{\sim}{\rightarrow} \bH_{K,v}^\wedge$, denoted $\chi \mapsto
\chi_v$, which is uniquely determined by the condition that 
\[ \chi(w)=\chi_v(T_w)|_{v=1} \qquad \mbox{for all $w\in W$}.\]
\end{thm}

Explicit tables or combinatorial algorithms for the values $\chi_v(T_w)$ 
(where $w$ has minimal length in its conjugacy class in $W$) are known 
for all types of $W$ and all weight functions $L$; this is one of the 
main themes of \cite{ourbuch}.

As before, $\bH$ is a symmetric algebra with trace function $\tau\colon\bH 
\rightarrow A$ given by $\tau(T_1)=1$ and $\tau(T_w)=0$ for $w\neq 1$.
By extension of scalars, we obtain a trace function $\tau_K \colon 
\bH_{K,v}\rightarrow K$. Since $\bH_{K,v}$ is split semisimple, we can 
write $\tau_K$ as a linear combination of $\bH_{K,v}^\wedge$, where all 
irreducible characters appear with a non-zero coefficient. (This is a 
general result about split semisimple symmetric algebras; see 
\cite[7.2.6]{ourbuch}.) Thus, we can write
\[ \tau_K=\sum_{\chi \in W^\wedge} \bc_\chi^{-1} \, \chi_v
\qquad \mbox{where} \qquad \bc_\chi \in K^\times.\]
Since $A$ is integrally closed in $K$, we have $\bc_\chi \in A$ for all 
$\chi \in W^\wedge$. (Again, this follows by a general argument on 
symmetric algebras; see \cite[7.3.9]{ourbuch}.) The
constants $\bc_{\chi}$ appear in the {\em orthogonality relations} for 
the irreducible characters of $\bH_{K,v}$; see \cite[7.2.4]{ourbuch}. 
Given $\chi,\chi'\in W^\wedge$, we have 
\[ \sum_{w \in W} v^{-L(w)} \chi_v(T_w)\,\chi_v(T_{w^{-1}})= \left
\{\begin{array}{cl} \chi(1)\, \bc_\chi &\qquad \mbox{if $\chi=\chi'$},\\
0 & \qquad \mbox{if $\chi\neq \chi'$}.\end{array}\right.\] 
We have the following more precise statement about the form of 
$\bc_\chi$; see the historical remarks in \cite[\S 10.7 and 
\S 11.6]{ourbuch} for the origins of this result.

\begin{thm}\label{schur} For any $\chi \in W^\wedge$, we have 
\[\bc_\chi=f_\chi v^{-2\alpha_\chi}\times \mbox{a product of cyclotomic
polynomials in $v$},\]
where $f_\chi$ is a positive integer and $\alpha_\chi \in \NN$. If 
$L(s)>0$ for all $s\in S$, then $f_\chi$ is divisible only by primes 
which are not good for $W$;  here, good primes are defined as in the 
remarks following Theorem~\ref{MGexistbs}.
\end{thm}

The polynomials $\bc_\chi$ are explicitly known for all $W,L$; see the 
appendices of \cite{Ca2} or \cite{ourbuch}. We write $\alpha_E=\alpha_\chi$ 
and $f_E=f_\chi$ if $E \in \Irr_{\QQ}(W)$ affords $\chi$.

\begin{defn} \label{Lgood} Let $p$ be a prime number. We say that
$p$ is $L$-good if $p$ does not divide any of the numbers $f_E$ for 
$E \in \Irr_{\QQ}(W)$. By Theorem~\ref{schur}, a good prime for $W$
is $L$-good if $L(s)>0$ for all $s\in S$. 
\end{defn}

\begin{rem}\label{extreme} Consider the extreme case where $L(s)=0$ for 
all $s\in S$. Then $\bH=AW$ and $\bc_\chi=|W|/\chi(1)$ for all $\chi 
\in W^\wedge$. Consequently, we have $f_\chi=1$ and $\alpha_\chi=0$ for 
all $\chi$. The $L$-good primes are prime numbers which do not divide
the order of $W$. Of course, this case does not give any new information. 
However, for $W$ of type $B_n$, it is interesting to look at cases 
where $L(s)=0$ for some $s\in S$; see Example~\ref{schurbn} below.
\end{rem}

\begin{exmp} \label{g2} Let $W=W(G_2)$ be the Weyl group of type $G_2$;
that is, we have $W=\langle s,t \mid s^2=t^2=(st)^6=1\rangle$. Since $s,t$
are not conjugate, we can take any $a,b \in \NN$ and obtain a weight 
function $L \colon W\rightarrow \NN$ such that $L(s)=a$ and $L(t)=b$.
We have 
\[\Irr_{\QQ}(W)=\{{\bf 1},\varepsilon,\varepsilon_1,\varepsilon_2,
E_{\pm}\}\]
where ${\bf 1}$ is the unit representation, $\varepsilon$ is the
sign representation, $\varepsilon_1$, $\varepsilon_2$ are two further
$1$-dimensional representations, and $E_{\pm 1}$ are two $2$-dimensional
representations. 
By \cite[8.3.4]{ourbuch}, the polynomials $\bc_E$ are given by 
\begin{align*}
\bc_{\bf 1} &= (v^{2a}{+}1)(v^{2b}{+}1)(v^{4a+4b}{+}v^{2a+2b}{+}1),\qquad 
\bc_{\varepsilon} =v^{-6a-6b}\, \bc_{\text{ind}},\\
\bc_{\varepsilon_1}&=v^{-6b}(v^{2a}{+}1)(v^{2b}{+}1)(v^{4a}{+}v^{2a+2b}
{+}v^{4b}), \qquad \bc_{\varepsilon_2}=v^{6b-6a}\bc_{\varepsilon_1},\\
\bc_{E_{\pm}}&=2v^{-2a-2b}(v^{2a+2b}{\pm} v^{a+b}{+}1)(v^{2a}{\mp} 
v^{a+b}{+} v^{2b}).
\end{align*}
This yields the following table:
\[\begin{array}{|c|cc|cc|cc|} \hline & \multicolumn{2}{c|}{b>a>0} & 
\multicolumn{2}{c|}{b=a>0} & \multicolumn{2}{c|}{b>a=0} \\ E & 
f_E&\alpha_E&f_E&\alpha_E&f_E&\alpha_E\\\hline
{\bf 1}       & 1&    0 & 1 &  0 & 2 & 0 \\
\varepsilon   & 1&3b+3a & 1 & 6a & 2 &3b \\
\varepsilon_1 & 1&3b-2a & 3 &  a & 2 &3b \\
\varepsilon_2 & 1&    a & 3 &  a & 2 &0 \\
E_{+}      & 2&    b & 6 &  a & 2 &b \\
E_{-}      & 2&    b & 2 &  a & 2 &b \\
\hline \end{array}\]
\end{exmp}


\begin{exmp} \label{f4} Let $W=W(F_4)$ be the Weyl group of type $F_4$,
with generators and diagram given by:
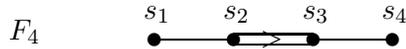
\begin{center}
\begin{picture}(200,20)
\put( 10, 5){$F_4$}
\put( 61,13){$s_1$}
\put( 91,13){$s_2$}
\put(121,13){$s_3$}
\put(151,13){$s_4$}
\put( 65, 5){\circle*{5}}
\put( 95, 5){\circle*{5}}
\put(125, 5){\circle*{5}}
\put(155, 5){\circle*{5}}
\put(105,2.5){$>$}
\put( 65, 5){\line(1,0){30}}
\put( 95, 7){\line(1,0){30}}
\put( 95, 3){\line(1,0){30}}
\put(125, 5){\line(1,0){30}}
\end{picture}
\end{center}
A weight function $L$ is specified by two integers $a:=L(s_1)=L(s_2)\geq 0$ 
and $b:=L(s_3)=L(s_4)\geq 0$. There are $25$ irreducible characters of 
$W$. The polynomials $\bc_E$ are listed in \cite[p.~379]{ourbuch}. One has
to distinguish a number of cases in order to obtain the invariants
$\alpha_E$ and $f_E$; see Table~\ref{invF4}. (Note that these are the
same cases that we found in \cite{my04}.)
\end{exmp}

\begin{table}[htbp] \caption{The invariants $f_\chi$ and $\alpha_\chi$
for type $F_4$} \label{invF4}
\begin{center}
$\begin{array}{|c|cc|cr|cc|cc|cc|} \hline & \multicolumn{2}{c|}{b{>}2a{>}0} 
& \multicolumn{2}{c|}{b{=}2a{>}0} & \multicolumn{2}{c|}{2a{>}b{>}a{>}0} & 
\multicolumn{2}{c|}{b{=}a{>}0} &\multicolumn{2}{c|}{b{>}a{=}0} \\ 
E& f_E & \alpha_E & f_E & \alpha_E & f_E 
& \alpha_E & f_E & \alpha_E & f_E  &\alpha_E \\ \hline
1_1 &1&0        &1&   0 & 1 & 0        & 1 & 0  & 6 & 0  \\
1_2 &1&12b{-}9a &2& 15a & 2 & 11b{-}7a & 8 & 4a & 6 & 12b  \\
1_3 &1&3a       &2&  3a & 2 & {-}b{+}5a& 8 & 4a & 6 & 0  \\
1_4 &1&12b{+}12a&1& 36a & 1 & 12b{+}12a& 1 &24a & 6 & 12b  \\
2_1 &1&3b{-}3a  &2&  3a & 2 & 2b{-}a   & 2 & a  &12 & 3b  \\
2_2 &1&3b{+}9a  &2& 15a & 2 & 2b{+}11a & 2 &13a &12 & 3b  \\
2_3 &1&a        &1&   a & 1 & a        & 2 & a  & 3 & 0  \\
2_4 &1&12b{+}a  &1& 25a & 1 & 12b{+}a  & 2 &13a & 3 & 12b  \\
4_1 &2&3b{+}a   &2&  7a & 2 & 3b{+}a   & 8 & 4a & 6 & 3b  \\
9_1 &1&2b{-}a   &2&  3a & 2 & b{+}a    & 1 & 2a & 2 & 2b \\
9_2 &1&6b{-}2a  &1& 10a & 1 & 6b{-}2a  & 8 & 4a & 2 & 6b  \\
9_3 &1&2b{+}2a  &1&  6a & 1 & 2b{+}2a  & 8 & 4a & 2 & 2b  \\
9_4 &1&6b{+}3a  &2& 15a & 2 & 5b{+}5a  & 1 &10a & 2 & 6b  \\
6_1 &3&3b{+}a   &3&  7a & 3 & 3b{+}a   & 3 & 4a &12 & 3b  \\
6_2 &3&3b{+}a   &3&  7a & 3 & 3b{+}a   &12 & 4a &12 & 3b  \\
12_1&3&3b{+}a   &6&  7a & 6 & 3b{+}a   &24 & 4a & 6 & 3b  \\
4_2 &1&b        &1&  2a & 1 & b        & 2 &  a & 6 & b  \\
4_3 &1&7b{-}3a  &1& 11a & 1 & 7b{-}3a  & 4 & 4a & 6 & 7b  \\
4_4 &1&b{+}3a   &1&  5a & 1 & b{+}3a   & 4 & 4a & 6 & b  \\
4_5 &1&7b{+}6a  &1& 20a & 1 & 7b{+}6a  & 2 &13a & 6 & 7b  \\
8_1 &1&3b       &1&  6a & 1 & 3b       & 1 & 3a &12 & 3b  \\
8_2 &1&3b{+}6a  &1& 12a & 1 & 3b{+}6a  & 1 & 9a &12 & 3b  \\
8_3 &1&b{+}a    &2&  3a & 2 & 3a       & 1 & 3a & 3 & b  \\
8_4 &1&7b{+}a   &2& 15a & 2 & 6b{+}3a  & 1 & 9a & 3 & 7b  \\
16_1&2&3b{+}a   &2&  7a & 2 & 3b{+}a   & 4 & 4a & 6 & 3b \\ \hline
\end{array}$
\end{center}
\end{table}

\begin{exmp} \label{schurbn} Let $n\geq 2$ and $W_n=W(B_n)$ be a Coxeter
group of type $B_n$ with generators and diagram given by:
\begin{center}
\begin{picture}(250,20)
\put( 10, 5){$B_n$}
\put( 61,13){$t$}
\put( 91,13){$s_1$}
\put(121,13){$s_2$}
\put(205,13){$s_{n-1}$}
\put( 65, 5){\circle*{5}}
\put( 95, 5){\circle*{5}}
\put(125, 5){\circle*{5}}
\put(215, 5){\circle*{5}}
\put( 65, 7){\line(1,0){30}}
\put( 65, 3){\line(1,0){30}}
\put( 95, 5){\line(1,0){50}}
\put(160, 5){\circle*{1}}
\put(170, 5){\circle*{1}}
\put(180, 5){\circle*{1}}
\put(195, 5){\line(1,0){20}}
\end{picture}
\end{center}
A weight function $L$ is specified by two integers $b:=L(t)\geq 0$
and $a:=L(s_1)=\cdots =L(s_{n-1})\geq 0$. It is well-known that we have a
parametrization
\[ \Irr_{\QQ}(W_n)=\{E^{\ulambda} \mid \ulambda
\vdash n\}.\]
Here, the notation $\ulambda\vdash n$ means that $\ulambda$ is a 
bipartition of $n$, that is, a pair of partitions $\ulambda=(\lambda_{(1)},
\lambda_{(2)})$ such that $|\lambda_{(1)}|+|\lambda_{(2)}|=n$. For 
example, the trivial representation is labelled by the pair $((n), 
\varnothing)$ and the sign representation is labelled by $(\varnothing, 
(1^n))$; see \cite[\S 5.5]{ourbuch}. The polynomials $\bc_{E^\ulambda}$
are determined as follows. Let $\ulambda=(\lambda_{(1)},\lambda_{(2)})
\vdash n$. By adding zeros if necessary, we write $\lambda_{(1)}$ and
$\lambda_{(2)}$ in the form 
\begin{align*}
\lambda_{(1)}&=\bigl(\lambda_{(1),1}\geq \lambda_{(1),2} \geq \cdots \geq
\lambda_{(1),m} \geq \lambda_{(1),m+1}\geq 0\bigr),\\
\lambda_{(2)}&=\bigl(\lambda_{(2),1}\geq \lambda_{(2),2} \geq \cdots \geq 
\lambda_{(2),m}\bigr),
\end{align*}
for some $m\geq 0$. Then we define the following symbol:
\[ \Lambda=\left(\begin{array}{c} \alpha_1 \qquad\alpha_2 \qquad\cdots 
\qquad \alpha_m\qquad \alpha_{m+1} \\ \beta_1\qquad  \beta_2\qquad  \cdots 
\qquad \beta_m \end{array}\right)\]
where $\alpha_i=i-1+\lambda_{(1),m+2-i}$ and $\beta_i=i-1+
\lambda_{(2),m+1-i}$ for all $i\geq 1$. Then
\begin{align*} \bc_{E^{\ulambda}}&=
\frac{v^{2am(2m+1)(m-2)/3} (v^{2a}+v^{2b})^m}{(v^{2a}-1)^n 
\displaystyle\prod_{i=1}^{m+1}\prod_{j=1}^m 
(v^{2a(\alpha_i-1)+2b}+ v^{2a\beta_j})}\\
&\qquad\times\frac{\displaystyle\prod_{i=1}^{m+1}
\prod_{k=1}^{\alpha_i}(v^{2ak}-1) (v^{2a(k-1)+2b}+1)}{\displaystyle
\prod_{1\leq i'<i\leq m+1} (v^{2a\alpha_i}-v^{2a\alpha_{i'}})}\\
&\qquad\qquad\times \frac{\displaystyle\prod_{j=1}^m\prod_{k=1}^{\beta_j} 
(v^{2ak}-1)(v^{2a(k+1)-2b}+1)}{\displaystyle\prod_{1\leq j'<j\leq m}
(v^{2a\beta_j}-v^{2a\beta_{j'}})} \cdot 
\end{align*}
see \cite[p.~447]{Ca2}. Lusztig \cite[Prop.~22.14]{Lusztig03} has obtained
explicit combinatorial formulae for the invariants $f_\chi$ and 
$\alpha_\chi$ (as a function of the parameters $a,b$). The following 
special cases are worth mentioning.

\medskip
$\bullet$ {\bf The asymptotic case\footnote{This term comes from the work
of Bonnaf\'e--Iancu \cite{BI}, where the Kazhdan--Lusztig cells in
this case are described.} in type $B_n$}.  Assume that 
$b>(n-1)a>0$. Then it is not hard to check (directly using the above 
formula) that 
\[ f_{E^\ulambda}=1 \quad \mbox{and}\quad \alpha_{E^\ulambda}=
b\,|\lambda_{(2)}|+a\,(n(\lambda_{(1)})+2n(\lambda_{(2)})- 
n(\lambda_{(2)}^*))\]
where the star denotes the conjugate partition and where, for any
partition $\nu$, we write $n(\nu):=\sum_{i=1}^r(i-1)\nu_i$ for
$\nu=(\nu_1 \geq \nu_2 \geq \cdots \geq \nu_r>0)$; see 
\cite[Remark~5.1]{my02}.

\medskip
$\bullet$ {\bf Type $A_{n-1}$.} The parabolic subgroup generated by 
$s_1,\ldots, s_{n-1}$ is isomorphic to the symmetric group $\Sym_n$, where 
$s_i$ corresponds to the transposition $(i,i+1)$. For a partition
$\nu$ of $n$, denote by $E^\nu$ the restriction of $E^{(\nu,\varnothing)}$ 
to $\Sym_n$. Then it is well-known that 
\[ \Irr_{\QQ}(\Sym_n)=\{E^\nu \mid \nu \vdash n\}.\]
For example, the partition $(n)$ labels the trivial representation and 
$(1^n)$ labels the sign representation. Furthermore, the formulae from 
the ``asymptotic case'' give the correct values for $f_E$ and 
$\alpha_E$, that is, for $a>0$, we have $f_{E^\nu}=1$ and 
$\alpha_{E^\nu}=n(\nu)a$ for all $\nu\vdash n$. See \cite[\S 10.5]{ourbuch}.

\medskip
$\bullet$ {\bf Type $D_n$.} Assume that $a>0$ and $L(t)=b=0$. We set 
$s_0:=ts_1t$. Then $W_n':=\langle s_0,s_1,\ldots,s_{n-1}\rangle$ is a 
Coxeter group of type $D_n$. For partitions $\lambda,\mu$ such that 
$(\lambda,\mu)\vdash n$,
we denote by $E^{[\lambda,\mu]}$ the restriction of $E^{(\lambda,\mu)}$ 
to $W_n'$. Then 
\begin{alignat*}{2}
E^{[\lambda,\mu]}&\cong E^{[\mu,\lambda]} \in \Irr(W_n') && \qquad 
\mbox{if $\lambda\neq \mu$},\\
E^{[\lambda,\lambda]}&\cong E^{[\lambda,+]}\oplus E^{[\lambda,-]} && 
\qquad \mbox{if $\lambda=\mu$ (and $n$ even)},
\end{alignat*}
where $E^{[\lambda,\pm ]}\in \Irr(W_n')$ and $E^{[\lambda,+]} \not\cong
E^{[\lambda,-]}$. This yields (see \cite[Chap.~5]{ourbuch} for more details):
\[ \Irr(W_n')=\{E^{[\lambda,\mu]} \mid (\lambda,\mu)\vdash n, \lambda 
\neq \mu\}\cup \{E^{[\lambda,\pm]} \mid \mbox{$n$ even}, \lambda \vdash 
n/2\}.\]
Let $L'$ be the restriction of $L$ to $W_n'$. Then $L'(s_i)=a$ for
$0 \leq i \leq n-1$, that is, $L'$ is just $a$ times the length function
on $W_n'$; see, for example, \cite[Lemma~1.4.12]{ourbuch}. By
\cite[Prop.~10.5.6]{ourbuch}, we have
\begin{alignat*}{2}
\alpha_{E^{[\lambda,\mu]}}&=\alpha_{E^{(\lambda,\mu)}} \quad \mbox{and}
\qquad f_{E^{[\lambda,\mu]}}&=f_{E^{(\lambda,\mu)}}, \quad 
&\mbox{if $\lambda\neq \mu$},\\ \alpha_{E^{[\lambda,\pm]}}&=
\alpha_{E^{(\lambda,\lambda)}} \quad \mbox{and} \quad 
f_{E^{[\lambda,\pm]}}&=2f_{E^{(\lambda,\lambda)}}, \quad 
&\mbox{if $n$ even and $\lambda\vdash n/2$}.
\end{alignat*}

\medskip
$\bullet$ Another extreme case is given by $a=0$ and $b>0$. Using the above 
formula for $\bc_{E^\ulambda}$,  one easily checks that 
$\alpha_{E^\ulambda}=|\lambda_{(2)}|\,b$ for all $\ulambda \vdash n$.
\end{exmp}

\begin{Par} {\bf Modular decomposition numbers.} \label{moddec}
Let $k$ be a field and $\theta \colon A \rightarrow k$ be a ring 
homomorphism. Let $\xi=\theta(v^2)$ and consider the specialized algebra 
$\bH_{k,\xi}$. Given $E\in \Irr_{\QQ}(W)$ and $M\in \Irr(\bH_{k,\xi})$,
we would like to define a {\em decomposition nunber} $[E:M]$, as in 
Brauer's modular representation theory of finite groups. Here, some care 
is needed since we are not necessarily working with discrete valuation 
rings. So let $\chi\in W^\wedge$ be the character of $E$ and $\chi_v$ be 
the corresponding irreducible character of $\bH_{K,v}$. Let 
\[ \rho_v \colon \bH_{K,v} \rightarrow M_d(K), \qquad T_w \mapsto 
(a_{ij}(T_w)),\]
be a matrix representation affording $\chi_v$, where $d=\dim E$.
Now, $\fp=\ker(\theta)$ is a prime ideal in $A$ and the localization 
$A_{\fp}$ is a regular local ring of Krull dimension $\leq 2$; see
Matsumura \cite{Matsu} for these notions. Hence, by Du--Parshall--Scott 
\cite[1.1.1]{DPS1}, we can assume that $\rho_v$ satisfies the condition
\[ \rho_v(T_w) \in M_d(A_{\fp}) \qquad \mbox{for all $w\in W$}.\]
Now, $\theta$ certainly extends to a ring homomorphism $\theta_{\fp}
\colon A_{\fp} \rightarrow k$. Applying $\theta_{\fp}$, we obtain a 
representation
\[ \rho_{k,\xi} \colon \bH_{k,\xi}\rightarrow M_d(k), \qquad
T_w \mapsto \bigl(\theta_{\fp}(a_{ij}(T_w))\bigr).\]
This representation may no longer be irreducible. For any $M \in
\Irr(\bH_{k,\xi})$, let $[E:M]$ be the multiplicity of $M$ as a
composition factor of the $\bH_{k,\xi}$-module affording $\rho_{k,\xi}$.
There are some choices involved in this process, but one can show
that $[E:M]$ is independent of these choices; see \cite[1.1.2]{DPS1}. 
Thus, we obtain a well-defined matrix
\[ D=\bigl([E:M]\bigr)_{E \in \Irr_{\QQ}(W),M \in \Irr(\bH_{k,\xi})}\]
which is called the {\em decomposition matrix} associated with $\theta$. 
(One can also define $D$ without using properties of regular local rings,
but then some mild hypotheses on the ground field $k$ are required;
see \cite[\S 7.4]{ourbuch}.)
\end{Par}

\begin{thm} \label{semispec} In the above setting, assume that
$\theta(\bc_\chi)\neq 0$ for all $\chi \in W^\wedge$. Then $\bH_{k,\xi}$
is split semisimple and, up to reordering the rows and columns, $D$
is the identity matrix.
\end{thm}

The fact that $\bH_{k,\xi}$ is split semisimple is proved by the argument
in \cite[Cor.~9.3.9]{ourbuch}. Once this is established, we can apply
{\em Tits' Deformation Theorem} (see \cite[7.4.6]{ourbuch}) and this
yields the statement concerning $D$. 

If the hypotheses of Theorem~\ref{semispec} are not satisfied, then our 
aim is to find a good parametrization of $\Irr(\bH_{k,\xi})$ using 
properties of $D$. The following example provides a model of what we are
looking for. 

\begin{exmp} \label{expsym} Let $W=\Sym_n$ be the symmetric group, with
generators $\{s_1,\ldots,s_{n-1}\}$ as usual. Since all $s_i$ are 
conjugate, we are in the ``equal parameter case'' and so $L(s_i)=a$
for all $i$. Let us assume that $a>0$. We set
\[ e:=\min \{i \geq 2 \mid 1+\xi^a+\xi^{2a}+\cdots + \xi^{(i-1)a}=0\}.\]
(If no such $i$ exists, we set $e=\infty$.) The following results are 
due to Dipper and James \cite{DJ0}. For any partition 
$\lambda\vdash n$, we have a corresponding {\em Specht module} 
$S_k^{\lambda}\in \bH_{k,\xi}\modA$. This module has the property that
\[ [E^\lambda:M]=\mbox{multiplicity of $M$ as a composition factor of
$S_k^{\lambda}$},\]
for any $M \in \Irr(\bH_{k,\xi})$. (Actually, we shall denote by 
$S_k^{\lambda}$ the module that is labelled by $\lambda^*$ in 
\cite{DJ0}, where the star denotes the conjugate partition.) Furthermore, 
there is an $\bH_{k,\xi}$-equivariant symmetric bilinear form on 
$S_k^{\lambda}$; then $\mbox{rad}(S_k^{\lambda})$, the radical of 
that form, is an $\bH_{k,\xi}$-submodule of $S_k^{\lambda}$. We set 
\[ D^{\lambda}:=S_k^{\lambda}/\mbox{rad}(S_k^{\lambda}) \in 
\bH_{k,\xi}\modA.\]
Let $\Lambda_n^{\circ}:=\{\lambda\vdash n\mid D^\lambda\neq \{0\}\}$. Then 
we have  
\[\Irr(\bH_{k,\xi})=\{D^{\lambda}\mid \lambda \in \Lambda_n^{\circ}\}\quad 
\mbox{and} \quad \Lambda_n^{\circ}=\{\lambda \vdash n \mid \mbox{$\lambda$ 
is $e$-regular}\}.\]
(A partition if $e$-regular if no part is repeated $e$ or more times.)
Thus, we have a ``natural'' parametrization of the simple 
$\bH_{k,\xi}$-modules by a subset of the indexing set for 
$\Irr_{\QQ}(\Sym_n)$. We wish to recover that parametrization from  
properties of the decomposition matrix $D$. For this purpose, we need to 
introduce the dominance order on partitions. Let $\lambda,\mu$ be 
partitions of $n$, with parts
\[\lambda=(\lambda_{1} \geq \lambda_{2} \geq \cdots \geq 0),\qquad
\mu=(\mu_1 \geq \mu_2 \geq \cdots \geq 0).\]
We write $\lambda \trianglelefteq \mu$ if $\sum_{i=1}^j \lambda_i \leq
\sum_{i=1}^j \mu_i$ for all $j$. Then we have:
\begin{equation*}
\left\{\begin{array}{l} \quad [E^\mu:D^\mu]=1 \quad \mbox{for any 
$\mu\in \Lambda_n^\circ$},\\ \quad [E^\lambda:D^\mu] \neq 0 \quad 
\Rightarrow \quad \lambda \trianglelefteq \mu,\end{array}\right. \tag{$*$}
\end{equation*}
see \cite[Theorem~7.6]{DJ0}. 
Note that the above conditions uniquely determine the set $\Lambda_n^\circ$. 
Indeed, let $M \in \Irr(\bH_{k,\xi})$. Then ($*$) shows that the set 
\[ \{\lambda \vdash n \mid [E^\lambda :M]\neq 0\}\]
has a unique maximal element with respect to $\trianglelefteq$, namely,
the unique $\mu \in \Lambda_n^{\circ}$ such that $M=D^\mu$. 
Finally, let us consider the invariants $\alpha_{E^\lambda}$. By
Example~\ref{schurbn}, we have $\alpha_{E^\lambda}=n(\lambda)a$.
Now it is known that, for any $\nu,\nu'\vdash n$, we have 
$\nu\trianglelefteq \nu'\Rightarrow n(\nu')\leq n(\nu)$, with equality 
only if $\nu=\nu'$ (see, for example, \cite[Exc.~5.6]{ourbuch}). Thus, 
given $M \in \Irr(\bH_{k,\xi})$ and setting
\[ \breve{\alpha}_M:=\min\{\alpha_{E^\lambda} \mid \lambda 
\vdash n \mbox{ and } [E^\lambda :M ]\neq 0\},\]
there is a unique $\mu\in \Lambda_n^{\circ}$ such that $\breve{\alpha}_M=
\alpha_{E^\mu}$ and $[E^\mu:M]\neq 0$. Hence, the set $\Lambda_n^{\circ}$ 
can be characterized as follows:
\[ \Lambda_n^{\circ}=\{ \lambda \vdash n\mid \exists M \in \Irr(\bH_{k,\xi})
\mbox{ such that } [E^\lambda:M]\neq 0, \breve{\alpha}_M=
\alpha_{E^\lambda}\}.\]
Note that this characterization does not require the explicit knowledge 
of all decomposition numbers.
\end{exmp}

We can now formalize the above discussion as follows.

\begin{defn} \label{canbase} Let $k$ be a field and $\theta \colon A 
\rightarrow k$ be a ring homomorphism; let $\xi=\theta(v^2)$. For any 
$M \in \Irr(\bH_{k,\xi})$, we define
\[ \breve{\alpha}_M:=\min\{\alpha_E\mid E \in \Irr_{\QQ}(W) \mbox{ and } 
[E:M ]\neq 0\}.\]
We say that a subset $\cB_{k,\xi} \subseteq \Irr_{\QQ}(W)$ is a 
``canonical basic set'' for $\bH_{k,\xi}$ if the following conditions are 
satisfied:
\begin{itemize}
\item[(a)] There is a bijection $\Irr(\bH_{k,\xi})\stackrel{\sim}{\rightarrow}
\cB_{k,\xi}$, denoted $M \mapsto E(M)$, such that $[E(M):M]=1$ and
$\alpha_{E(M)}=\breve{\alpha}_M$.
\item[(b)] Given $E\in \Irr_{\QQ}(W)$ and $M \in \Irr(\bH_{k,\xi})$, we have
\[ [E:M]\neq 0 \qquad \Rightarrow \qquad \breve{\alpha}_M<\alpha_E \quad 
\mbox{or} \quad E=E(M).\]
\end{itemize}
Note that, if a canonical basic set exists, then (a) and (b) uniquely 
determine the set $\cB_{k,\xi}\subseteq \Irr_{\QQ}(W)$ and the 
bijection $M \mapsto E(M)$. If $\cB_{k,\xi}$ exists, then the submatrix 
\[ D^\circ=\bigl([E(M):M']\bigr)_{M,M' \in \Irr(\bH_{k,\xi})}\]
is square and lower triangular with $1$ on the diagonal, when we order 
the modules in $\Irr(\bH_{k,\xi})$ according to increasing values of 
$\breve{\alpha}_M$. More precisely, we have a block lower triangular shape
\[ D^\circ=\left(\begin{array}{cccc} D_0^\circ &&&0 \\
& D_1^{\circ} &&\\ && \ddots & \\ *&&& D_N^\circ \end{array}\right),\]
where the block $D_i^\circ$ has rows and columns labelled by those $E(M)$
and $M'$, respectively, where $\breve{\alpha}_M=\breve{\alpha}_{M'}=i$, 
and each $D_i^{\circ}$ is the identity matrix.
\end{defn}

\begin{rem} \label{semispec1} Assume that $\theta(\bc_\chi)\neq 0$ for 
all $\chi \in W^\wedge$. Then, by Theorem~\ref{semispec}, $D$ is
the identity matrix and so $\cB_{k,\xi}=\Irr_{\QQ}(W)$ is the unique
canonical basic set. Hence, the interesting cases are those where
$\theta(\bc_\chi)=0$ for some $\chi$. By Theorem~\ref{schur}, this implies 
that the characteristic of $k$ must be a prime which is not $L$-good, 
or that $\theta(v)$ must be a root of unity in~$k$. 
\end{rem}

The following example shows that canonical basic sets do not always exist.

\begin{exmp} \label{decg2} Let us consider the Iwahori--Hecke algebra of 
type $G_2$ as in Example~\ref{g2}. Let $L$ be the weight function given
by length, let $k$ be a field of characteristic~$2$ and $\theta \colon 
A \rightarrow k$ be a specialization such that $\xi=\theta(v^2)=1$. Using 
the explicit matrix representations in \cite[\S 8.6]{ourbuch}, we find 
that the decomposition matrix is given by 
\[ \begin{array}{|cc|cc|} \hline E & \alpha_E & 
\multicolumn{2}{c|}{[E:M]} \\ \hline
{\bf 1}       & 0 & 1 & 0 \\
\varepsilon_1 & 1 & 1 & 0 \\
\varepsilon_2 & 1 & 1 & 0 \\
\varepsilon   & 6 & 1 & 0 \\
E_+           & 1 & 0 & 1 \\
E_-           & 1 & 0 & 1 \\
\hline \end{array}\]
Thus, there is no subset $\cB_{k,1}$ satisfying the conditions in 
Definition~\ref{canbase}. Note that $2$ is not $L$-good in this case. 
\end{exmp}

In Section~6, we will show by a general method, following Geck 
\cite{mykl}, \cite{my00} and Geck--Rouquier \cite{GeRo2}, that 
canonical basic sets do exist when the characteristic of $k$ is zero 
or a prime which is $L$-good. This general method relies on some
properties of the Kazhdan-Lusztig basis of $\bH$, which are not yet
known to hold in general if $L$ is not constant on $S$. But note
that, as far as unequal parameters are concerned, we only have to deal
with groups of type $G_2$, $F_4$ and $B_n$ (any $n\geq 2$). By an 
explicit (and easy) computation, one can determine the decomposition 
matrices in type $G_2$ for all choices of $L$ and $\theta \colon A 
\rightarrow k$ (using the matrix representations in \cite[\S 8.6]{ourbuch}).
By inspection, one finds that there is a canonical basic set $\cB_{k,\xi}$ 
if the characteristic of $k$ is zero or an $L$-good prime. A similar 
statement can be verified in type $F_4$, using the decomposition matrices 
computed by Geck--Lux \cite{GeLu}, Bremke \cite{Br} and 
McDonough--Pallikaros \cite{McPa1}. Thus, the remaining case is type
$B_n$, and this will be discussed in Example~\ref{canbn} and Section~8.

\section{The Kazhdan--Lusztig basis and the $\ba$-function} \label{MGsec4a}

We keep the setting of the previous section, where $W$ is a finite Weyl 
group and $L \colon W \rightarrow \NN$ is a weight function. Now we
introduce the Kazhdan--Lusztig basis of $\bH$ and explain some theoretical 
constructions arising from it, most notably Lusztig's ring $\bJ$. We will
illustrate the power of these methods by giving new and conceptual proofs 
for some of the results in Section~4 on the structure of $\bH_{K,v}$.
In the following section, we shall discuss applications to non-semisimple 
specializations of $\bH$. 

Let us assume that $L(s)>0$ for all $s\in S$. It will be convenient to 
rescale the basis elements of $\bH$ as follows:
\[ \tT_w:=v^{-L(w)}T_w \qquad \mbox{for any $w\in W$}.\]
Then the multiplication formulae read:
\[\tT_s\tT_w=\left\{\begin{array}{cl}\tT_{sw}&\quad 
\mbox{if $l(sw)>l(w)$},\\ \tT_{sw}+(v^{L(s)}-v^{-L(s)})\tT_w &
\quad \mbox{if $l(sw)<l(w)$},\end{array} \right.\]
where $s \in S$ and $w \in W$. To define the Kazhdan--Lusztig basis of 
$\bH$, we need some ring homomorphisms of $\bH$. First of all, we have a 
ring homomorphism $A \rightarrow A$, $f \mapsto \bar{f}$, where $\bar{f}$ 
is obtained from $f$ by substituting $v \mapsto v^{-1}$. This extends to a 
ring homomorphism $j \colon \bH \rightarrow \bH$ such that 
\[ j\colon \sum_{w \in W} a_w \tT_w\mapsto \sum_{w\in W} 
\bar{a}_w (-1)^{l(w)} \tT_w;\]
we have $j^2=\mbox{id}_{\bH}$. Next, define an $A$-linear map
$\bH\rightarrow \bH$, $h \mapsto h^\dagger$, by 
\[ \tT^\dagger:=(-1)^{l(w)} \tT_{w^{-1}}^{-1}\qquad (w \in W).\]
Then $h \mapsto h^\dagger$ is an $A$-algebra automorphism whose square is 
the identity; furthermore, $\dagger$ and $j$ commute with each other. 
Hence we obtain a ring involution $h \mapsto \bar{h}$ of $\bH$ by composing
$j$ and $\dagger$, that is, we have 
\[ \overline{\sum_{w \in W} a_w \tT_w}:=\sum_{w\in W} \bar{a}_w 
\tT_{w^{-1}}^{-1}.\]

\begin{thm}[Kazhdan--Lusztig; see Lusztig \cite{Lusztig03}] \label{klbase}
For each $w \in W$, there exists a unique $c_w \in \bH$ such that 
\[ \overline{c}_w=c_w \qquad \mbox{and} \qquad  c_w \equiv \tT_w
\bmod \bH_{<0},\]
where $\bH_{<0}$ denotes the set of all $v^{-1}{\ZZ}[v^{-1}]$-linear 
combinations of basis elements $\bT_y$ ($y \in W$). The elements 
$\{c_w \mid w \in W\}$ form an $A$-basis of $\bH$.
\end{thm}

For example, we have $c_1=\tT_1$ and $c_s=\tT_s+v^{-L(s)}\tT_1$ for $s\in S$.

Now the key to understanding the representations of a specialized algebra 
$\bH_{k,\xi}$ is the construction of Lusztig's ring $\bJ$.  This involves 
the following ingredients. Since $\{c_w\}$ is an $A$-basis of $\bH$, we can 
write
\[ c_xc_y=\sum_{z\in W} h_{x,y,z} c_z \qquad \mbox{where $h_{x,y,z}
\in A$};\]
we have $\bar{h}_{x,y,z}=h_{x,y,z}$ for all $x,y,z \in W$. We define a 
function $\ba \colon W \rightarrow \NN$ as follows. Let $z \in W$.
Then we set 
\[ \ba(z):=\min \{i \geq 0 \mid v^i h_{x,y,z} \in {\ZZ}[v]\mbox{ for
all $x,y \in W$}\}.\]
It is easy to see that $\ba(1)=0$ and that $\ba(z)=\ba(z^{-1})$ for all
$z \in W$. For $x,y,z\in W$, we set 
\[ \gamma_{x,y,z}:=\mbox{ constant term of $v^{\ba(z)}h_{x,y,z^{-1}}
\in {\ZZ}[v]$}.\]
Thus, we obtain a family of integers $\{\gamma_{x,y,z} \mid x,y,z\in W\}
\subseteq \ZZ$ and we can try to use them to define a ring. Let $\bJ$ be 
the free abelian group with basis $\{t_w \mid w \in W\}$. We define a 
bilinear product on $\bJ$ by 
\[ t_x \cdot t_y:=\sum_{z \in W} \gamma_{x,y,z} \,t_{z^{-1}} \qquad
(x,y \in W).\]
We would like to show that $\bJ$ is an associative ring with an identity 
element.  For this purpose, we need some further (and rather subtle) 
properties of the Kazhdan--Lusztig basis of $\bH$. 

\begin{Par} {\bf Lusztig's conjectures.} \label{III} Let $z \in W$ and
consider $\tau(c_z) \in A$. By Theorem~\ref{klbase} and 
\cite[Prop.~5.4]{Lusztig03}, we can write
\[ \tau(c_z)=n_z \,v^{-\Delta(z)}+\mbox{combination of smaller powers 
of $v$},\]
where $n_z \in \ZZ\setminus \{0\}$ and $\Delta(z)\in \NN$. We set 
\[ \cD:=\{z \in W \mid \ba(z)=\Delta(z)\}.\]
One easily checks that $1 \in \cD$ and that $\cD^{-1}=\cD$. Now
Lusztig \cite[Chap.~14]{Lusztig03} has formulated $15$ properties P1--P15
of the above objets ($h_{x,y,z}$, $\gamma_{x,y,z}$, $\cD$ etc.) and 
conjectured that they always hold. For our purposes here, we only need 
the following ones:
\begin{itemize}
\item[\bf P2.] If $d \in \cD$ and $x,y\in W$ satisfy $\gamma_{x,y,d}\neq 0$,
then $x=y^{-1}$.
\item[\bf P3.] If $y\in W$, there exists a unique $d\in \cD$ such that
$\gamma_{y^{-1},y,d}\neq 0$.
\item[\bf P4.] If $x,y,z \in W$ satisfy $h_{x,y,z} \neq 0$, then 
$\ba(z)\geq \ba(x)$ and $\ba(z)\geq \ba(y)$.
\item[\bf P5.] If $d\in \cD$, $y\in W$, $\gamma_{y^{-1},y,d}\neq 0$, then
$\gamma_{y^{-1},y,d}=n_d=\pm 1$.
\item[\bf P6.] If $d\in \cD$, then $d^2=1$.
\item[\bf P7.] For any $x,y,z\in W$, we have $\gamma_{x,y,z}=\gamma_{y,z,x}$.
\item[\bf P8.] Let $x,y,z\in W$ be such that $\gamma_{x,y,z}\neq 0$. Then
$\ba(x)=\ba(y)=\ba(z)$. 
\item[\bf P15'.] If $x,x',y,w\in W$ satisfy $\ba(w)=\ba(y)$, then
\[\sum_{u \in W} h_{x,u,y}\,\gamma_{w,x',u^{-1}}=\sum_{u\in W} h_{x,w,u}\,
\gamma_{u,x',y^{-1}}.\]
\end{itemize}
The above properties are known to hold in the following situations:
\begin{itemize}
\item[(I)] Recall that we are assuming that $W$ is a finite Weyl group. 
Let us also assume that the weight function $L$ is constant on $S$ (the 
``equal parameter case''). Then, thanks to a deep geometric interpretation 
of the basis $\{c_w\}$, one can show that P1--P15 hold; see Lusztig 
\cite[Chap.~15]{Lusztig03} and the references there. In the case where 
$W=\Sym_n$ is the symmetric group, elementary proofs for P1--P15 (that 
is, without reference to a geometric interpretation) can be found in 
\cite{my05a}.
\item[(II)] Assume that $W$ is of type $B_n$ with diagram and weight
function given as follows:
\begin{center}
\begin{picture}(250,40)
\put( 10,25){$B_n$}
\put( 10,05){$L$ :}
\put( 65,05){$b$}
\put( 95,05){$a$}
\put(125,05){$a$}
\put(215,05){$a$}
\put( 61,33){$t$}
\put( 91,33){$s_1$}
\put(121,33){$s_2$}
\put(205,33){$s_{n-1}$}
\put( 65,25){\circle*{5}}
\put( 95,25){\circle*{5}}
\put(125,25){\circle*{5}}
\put(215,25){\circle*{5}}
\put( 65,27){\line(1,0){30}}
\put( 65,23){\line(1,0){30}}
\put( 95,25){\line(1,0){50}}
\put(160,25){\circle*{1}}
\put(170,25){\circle*{1}}
\put(180,25){\circle*{1}}
\put(195,25){\line(1,0){20}}
\end{picture}
\end{center}
where $a,b \in \NN$. If $b > (n-1)a>0$ (the ``asymptotic case'' as 
in Example~\ref{schurbn}), then P1--P14 and P15' are known to hold by 
Bonnaf\'e--Iancu \cite{BI}, Bonnaf\'e \cite{BI2}, Geck \cite{my05} and 
Geck--Iancu \cite{GeIa05}. 

If $b=0$ and $a>0$ (the case relevant for type $D_n$), then P1--P15 are 
known to hold by a reduction argument to case~(I); this is already due to
Lusztig (see \cite[\S 2]{my00} and the references there).
\end{itemize}
Partial results for type $F_4$ with unequal parameters are contained in
\cite{my04}. 
\end{Par}

From now on, we will assume that P2--P8 and P15' hold for $W, L$.

\begin{thm}[Lusztig \protect{\cite[18.3, 18.9]{Lusztig03}}] \label{phi} 
$\bJ$ is an associative ring with identity  $1_{\bJ}=\sum_{d \in \cD} 
n_dt_d$.  Furthermore, setting $\bJ_A= A\otimes_{\ZZ} \bJ$, the $A$-linear 
map $\phi \colon \bH \rightarrow \bJ_A$ defined by 
\[ \phi(c_w^\dagger)=\sum_{\atop{z\in W, d \in \cD}{\ba(z)=\ba(d)}}
h_{w,d,z} \, \hat{n}_z \, t_z \qquad (w \in W)\]
is a homomorphism preserving the identity elements. Here, we set 
$\hat{n}_z:=n_d$ where $d \in \cD$ is the unique element such that
$\gamma_{z,z^{-1},d}\neq 0$.
\end{thm}

Now let $a \geq 0$ and consider the following $A$-submodules of $\bH$:
\[ \bH^{\geq a}:=\langle c_w^\dagger \mid \ba(w)\geq a\rangle_A\quad
\mbox{and}\quad \bH^{>a}:=\langle c_w^\dagger \mid \ba(w)>a\rangle_A.\]
Then, by P4, both $\bH^{\geq a}$ and $\bH^{>a}$ are two-sided ideals of
$\bH$. Hence 
\[ \bH^a:=\bH^{\geq a}/\bH^{>a}=\langle [c_w^\dagger] \mid \ba(w)=a
\rangle_A\]
is an $(\bH,\bH)$-bimodule. We define an $A$-bilinear map $\bJ_A \times 
\bH^a \rightarrow \bH^a$ by
\[ t_x \star [c_w^\dagger]:=\sum_{z \in W} \gamma_{x,w,z^{-1}}\, \hat{n}_w
\hat{n}_z\, [c_z^\dagger] \]
where $x \in W$ and $w \in W$ is such that $\ba(w)=a$.

\begin{thm}[Lusztig \protect{\cite[18.10]{Lusztig03}}] \label{bimod} 
$\bH^a$ is a $(\bJ_A,\bH)$-bimodule and we have 
\[ h.[c_w^\dagger]=\phi(h) \star [c_w^\dagger] \qquad\mbox{for all
$h \in \bH$ and $w \in W$, $\ba(w)=a$}.\]
\end{thm}

We invite the reader to check that the above two theorems indeed are proved 
by purely algebraic arguments, using only P2--P8 and P15'. 

Now consider a ring homomorphism $\theta \colon A \rightarrow k$ where 
$k$ is a field. Since all of the above constructions are defined over $A$, 
we can extend scalars from $A$ to $k$ and obtain:
\begin{itemize}
\item the specialized algebras $\bH_{k,\xi}$ (where $\xi=\theta(v^2)$) and 
$\bJ_k=k\otimes_A J_A$; 
\item the induced $k$-algebra homomorphism $\phi_{k,\xi}\colon 
\bH_{k,\xi} \rightarrow\bJ_k$;
\item a $(\bJ_k,\bH_{k,\xi})$-bimodule structure on $\bH^a_{k,\xi}$ 
such that
\[ h.[c_w^\dagger]=\phi_{k,\xi}(h)\star [c_w^\dagger] \]
for all $h\in \bH_{k,\xi}$ and all $w \in W$ such that $\ba(w)=a$.
\end{itemize}
We are now in a position to obtain a number of applications.

\begin{prop}[Lusztig \protect{\cite[18.12]{Lusztig03}}] \label{radkl}  In
the above setting, the kernel of $\phi_{k,\xi}$ is a nilpotent ideal and, 
hence, contained in the Jacobson radical of $\bH_{k,\xi}$. Consequently, 
$\phi_{k,\xi}$ is an isomorphism if $\bH_{k,\xi}$ is semisimple.
\end{prop}

\begin{proof} To illustrate the use of the above constructions, we repeat
Lusztig's proof here. Let $h\in \ker(\phi_{k,\xi})$ and $a\geq 0$. Then 
$h.[c_w^\dagger]=\phi_{k,\xi}(h)\star [c_w^\dagger]=0$ for all $w \in W$ 
such that $\ba(w)=a$. In other words, this means that 
\[ h\bH_{k,\xi}^{\geq a}\subseteq \bH_{k,\xi}^{\geq a+1} \qquad 
\mbox{for $a=0,1,2, \ldots$}.\]
Now let $N:=\max\{\ba(z)\mid z\in W\}$. Repeating the above argument $N+1$
times, we obtain 
\[ h_1 \cdots h_{N+1}=h_1\cdots h_{N+1}c_1^\dagger \in 
\bH_k^{\geq N+1}=\{0\}\]
for any $h_1,\ldots,h_{N+1} \in \ker(\phi_{k,\xi})$. Hence, 
$\ker(\phi_{k,\xi})$ is a nilpotent ideal. 
\end{proof}

The following examples are taken from Lusztig \cite[Chap.~20]{Lusztig03}.

\begin{exmp} \label{exp1} Consider the specialization $\theta\colon A 
\rightarrow \QQ$ such that $\theta(v)=1$. Then $\bH_{\QQ,1}={\QQ}W$ is the 
group algebra of $W$ and $J_{\QQ}={\QQ} \otimes_{\ZZ} J$ is the algebra 
obtained by extension of scalars from $\ZZ$ to $\QQ$. Since $\QQ W$ is 
semisimple, Proposition~\ref{radkl} shows that  
\[\phi_{\QQ,1} \colon {\QQ}W \rightarrow J_{\QQ}\]
is an isomorphism. In particular, $J_{\QQ}$ is split semisimple. 

More generally, let $k$ be a field whose characteristic is either zero or a 
prime which does not divide the order of $W$. Consider the specialization
$\theta \colon A \rightarrow k$ such that $\theta(v)=1$. Then, again, we have
$\bH_{k,1}=kW$ and this is a split semisimple algebra. Hence, as above, 
$\phi_{k,1} \colon kW \rightarrow \bJ_k$ is an isomorphism. 
\end{exmp}

\begin{exmp} \label{exp2} Let $K=\QQ(v)$, the field of fractions of $A$,
and consider the specialization $\theta\colon A \rightarrow K$ given
by inclusion. Now write 
\[ \phi_{K,v}(\tT_y)=\sum_{x\in W} b_{x,y}\,t_x \qquad \mbox{where}
\qquad b_{x,y} \in A.\]
Explicitly, the coefficients $b_{x,y}$ are obtained by writing $\tT_y$
as an $A$-linear combination of $c_w^\dagger$ and then to use the defining
formula of $\phi$. Thus, $B=(b_{x,y})_{x,y\in W}$ is a matrix with entries 
in $A$. If we set $v=1$, we obtain the matrix of the homomorphism 
$\phi_{\QQ,1} \colon {\QQ}W \rightarrow J_{\QQ}$ considered in the previous 
example. Since $\phi_{\QQ,1}$ is an isomorphism, the determinant of the 
matrix $B$ is a Laurent polynomial whose value at $v=1$ is non-zero. In 
particular, $\det B \neq 0$ and so $\phi_{K,v}$ is an isomorphism.  Hence
we obtain isomorphisms 
\[ \bH_{K,v} \stackrel{\phi_{K,v}}{\longrightarrow}\bJ_K 
\stackrel{\phi_{K,1}^{-1}}{\longrightarrow} KW;\]
in particular, $\bH_{K,v}$ is split semisimple and isomorphic to $KW$. 
Thus, using Lusztig's ring $J$, we have recovered Theorem~\ref{gen1} by 
a general argument (assuming that P2--P8 and P15' hold).
\end{exmp}

Now, via the algebra isomorphisms constructed above, we may identify 
the following sets of simple modules:
\[ \fbox{$\Irr_{\QQ}(W)=\Irr(\bJ_{\QQ})=\Irr(\bJ_K)=\Irr(\bH_{K,v})$}\]
where the second equality is given by extension of scalars from $\QQ$ to
$K$. We shall denote these correspondences as follows. Let $E \in 
\Irr_{\QQ}(W)$. Composing the action of $W$ on $E$ with the  inverse of
the isomorphism ${\QQ}W \stackrel{\sim}{\rightarrow} \bJ_{\QQ}$ in 
Example~\ref{exp1}, we obtain a simple $\bJ_{\QQ}$-module denoted by 
$E_\spadesuit$. Extending scalars from $\QQ$ to $K$, we obtain
$E_{\spadesuit,K} \in \Irr(\bJ_K)$. Finally, composing the action of
$\bJ_K$ on $E_{\spadesuit,K}$ with the isomorphism $\bH_{K,v}
\stackrel{\sim}{\rightarrow} \bJ_K$ in Example~\ref{exp2}, we obtain
a simple $\bH_{K,v}$-module denoted by $E_v$. Thus, we have the
correspondences:
\[ E \in \Irr_{\QQ}(W) \quad \leftrightarrow \quad 
E_{\spadesuit} \in \Irr(\bJ_{\QQ}) \quad \leftrightarrow \quad
E_{v} \in \Irr(\bH_{K,v}).\]
We can also express this in terms of characters. As in the previous
section, denote by $\bH_{K,v}^\wedge$ the set of irreducible characters 
of $\bH_{K,v}$. The set $\bJ_{\QQ}^\wedge$ is defined similarly. Then we 
also have bijective correspondences:
\[ \chi \in W^\wedge\quad \leftrightarrow \quad 
\chi_{\spadesuit} \in \bJ_{\QQ}^\wedge \quad \leftrightarrow \quad
\chi_v \in \bH_{K,v}^\wedge.\]
(If $\chi$ is afforded by $E\in \Irr_{\QQ}(W)$, then $\chi_\spadesuit$ is
the character afforded by $E_\spadesuit$ and $\chi_v$ is the character 
afforded by $E_v$.) We have $\chi(w) \in \ZZ$ for all $w \in W$. Similarly, 
we have $\chi_\spadesuit(t_w) \in \ZZ$ and $\chi_v(\tT_w) \in A$ for all
$w\in W$. This follows by a general argument, since $\ZZ$ is integrally 
closed in $\QQ$ and $A$ is integrally closed in $K$; see 
\cite[7.3.9]{ourbuch}. The discussion in Example~\ref{exp2} easily implies
that 
\[ \chi(w)=\chi_v(T_w)|_{v=1} \qquad \mbox{for all $w\in W$};\]
see \cite[20.3]{Lusztig03}. Thus, we have recovered Theorem~\ref{gen2} by
a general argument (assuming that P2--P8 and P15' hold).

Now the general philosophy will be that, although we don't really know
the algebra $\bJ$ explicitly, it will serve as a theoretical tool for
various module-theoretic constructions. We will see the full power of
this in the following section, where we consider non-semisimple
specialisations. Here, on the ``generic level'', we now use this idea to
give a new interpretation to the invariants $\alpha_\chi$ and $f_\chi$
defined by the formula in Theorem~\ref{schur}. This is done as follows. 
We have a direct sum decomposition 
\[ \bJ=\bigoplus_{a\geq 0} \bJ^a \quad \mbox{where} \quad \bJ^a:=
\langle t_w \mid w \in W \mbox{ such that } \ba(w)=a\rangle_{\ZZ}.\]
By P7, each $\bJ^a$ is a two-sided ideal in $\bJ$. In fact, one easily 
checks that 
\[ t_a:=\sum_{\atop{d \in W}{\ba(d)=a}} n_dt_d \in \bJ^a\]
is a central idempotent in $\bJ$; furthermore, we have $1_J=\sum_{a\geq 0}
t_a$ and $t_at_{a'}=0$ for $a\neq a'$. Hence, for any simple $M\in 
\Irr(\bJ_{\QQ})$, there exists a unique $a\geq 0$ such that $t_a.M=M$
and $t_{a'}.M=\{0\}$ for $a'\neq a$. Note that we have 
\[ \forall z \in W : \quad t_z.M \neq \{0\} \quad
\Rightarrow \quad \ba(z)=a.\]

\begin{prop}[Lusztig] \label{ainv} Let $E \in \Irr_{\QQ}(W)$ and consider 
the corresponding $\bJ_{\QQ}$-module $E_\spadesuit$. Then we have 
$\alpha_E=a$ where $a \geq 0$ is uniquely determined by the condition
that $t_a.E_\spadesuit=E_\spadesuit$ (and $t_{a'}.E_\spadesuit=\{0\}$
for $a'\neq a$). 
\end{prop}

Let us sketch the main arguments of the proof, following Lusztig 
\cite[Chap.~20]{Lusztig03}. Let $\chi$ be the character afforded by
$E\in \Irr_{\QQ}(W)$. Let $a\geq 0$ be such that $t_a.E_\spadesuit=
E_\spadesuit$, as above. First note that, by definition, we have
\[ \chi_v(c_w^\dagger)=\chi_\spadesuit(\phi_{K,v}(c_w^\dagger))=
\sum_{\atop{z\in W,d \in \cD}{\ba(z)=\ba(d)}} h_{w,d,z}\hat{n}_z \,
\chi_\spadesuit(t_z),\]
where the sum need only be extended over all $z \in W$ such that $\ba(z)=a$.
Using the properties P2--P7, one shows that 
\[ v^{\ba(w)}\, \chi_v(c_w^\dagger) \in {\ZZ}[v] \qquad \mbox{and}
\qquad v^{\ba(w)}\, \chi_v(c_w^\dagger) \equiv \chi_\spadesuit(t_w) \;
\bmod v{\ZZ}[v]\]
for all $w \in W$. Now we claim that we have 
\begin{equation*}
a=\max\{\ba(w) \mid w \in W \mbox{ and } \chi_v(c_w^\dagger) 
\neq 0\}.\tag{$*$}
\end{equation*}
First note that $\chi_{\spadesuit}(t_w)\neq 0$ for some $w\in W$ and, 
hence, $\chi_v(c_w^\dagger)\neq 0$. Since $\ba(w)=a$, we have 
the inequality ``$\leq$''. On the other hand, assume that $\chi_{v}
(c_w^\dagger) \neq 0$. Then there exist some $z \in W$, $d\in \cD$ such 
that $h_{w,d,z}\neq 0$ and $\ba(z)=a$. Now P4 shows that 
$\ba(w)\leq \ba(z)$. Thus, ($*$) holds. 

Consequently, we have 
\[ v^{a}\, \chi_v(c_w^\dagger) \in {\ZZ}[v] \quad \mbox{and}
\quad v^{a}\, \chi_v(c_w^\dagger) \equiv \chi_\spadesuit(t_w) \;
\bmod v{\ZZ}[v]\]
for all $w \in W$. Now recall that $c_w \equiv \tT_w \bmod \bH_{<0}$. 
Since $\overline{c}_w=c_w$, we have $c_w^\dagger=j(c_w)$ and so 
$c_w^\dagger\equiv (-1)^{l(w)}\tT_w \bmod \bH_{>0}$. These relations imply
that we can write 
\[ \tT_w=(-1)^{l(w)}c_w^\dagger + \mbox{$v{\ZZ}[v]$-combination of various 
$c_y^\dagger$ ($y \in W$)}.\]
Hence we also have 
\[ v^{a}\, \chi_v(\tT_w) \in {\ZZ}[v] \quad \mbox{and}
\quad v^{a}\, \chi_v(\tT_w) \equiv (-1)^{l(w)}
\chi_\spadesuit(t_w) \; \bmod v{\ZZ}[v]\]
for all $w \in W$. Inserting the above congruence conditions into the 
orthogonality relations for the irreducible characters of $\bH_{k,v}$ 
from the previous section, we deduce that 
\[ v^{2a}\,\chi(1)\,\bc_\chi \equiv \Bigl(\sum_{w \in W} 
\chi_{\spadesuit}(t_w) \chi_{\spadesuit}(t_{w^{-1}})\Bigr) \bmod 
v{\ZZ}[v].\]
Finally, by \cite[20.1]{Lusztig03}, $\bJ$ also is a symmetric algeba, 
with trace function $\mu\colon \bJ \rightarrow  \ZZ$ given by 
\[ \mu(t_z)=\left\{\begin{array}{cl} n_z & \qquad \mbox{if $z \in \cD$},\\
0 & \qquad \mbox{otherwise}.\end{array}\right.\]
We have $\mu(t_xt_y)=1$ if $x=y^{-1}$ and $\mu(t_xt_y)=0$ otherwise. (To
check this, use P2, P3, P5, P6). Hence, we can write 
\[ \mu_{\QQ}=\sum_{\psi \in W^\wedge} d_\psi^{-1}\, \psi_{\spadesuit}
\qquad \mbox{where}\qquad 0\neq d_\psi \in \ZZ.\]
By the orthogonality relations for the irreducible characters of $\bJ_{\QQ}$,
we have 
\[ \sum_{w \in W} \chi_\spadesuit(t_w)\,\chi_\spadesuit(t_{w^{-1}})=
\chi(1) d_\chi.\] 
A comparison with the above formula shows that
\[\bc_\chi=d_\chi v^{-2a}+\mbox{combination of higher 
powers of $v$}.\]
Hence, we must have $a=\alpha_\chi$ and $d_\chi=f_\chi$, as desired. 
Thus, both $\alpha_\chi$ and $f_\chi$ can be interpreted in terms of 
the algebra $\bJ$. \hfill $\Box$

\begin{rem} \label{weight0} One can actually check that all the arguments
in this section go through in the case where we consider a weight function
$L \colon W \rightarrow \NN$ and allow the possibility that $L(s)=0$ for
some $s\in S$. (In fact, Lusztig \cite[5.2]{Lusztig03} proves 
Theorem~\ref{klbase} in this more general set-up. Some care is needed in
the definition of $\cD$ since it may happen that $\tau(c_z)=0$ for 
$z\in W$. Everything works out well by setting $\Delta(z)=\infty$ if
$\tau(c_z)=0$.) 
\end{rem}

\section{Canonical basic sets and Lusztig's ring $J$} \label{MGsec5}

We keep the setting of the previous section. Now let $\theta \colon A 
\rightarrow k$ be a ring homomorphism into a field $k$ whose 
characteristic is either zero or a prime which is $L$-good; see 
Definition~\ref{Lgood}. As before, let 
\[ \bH_{k,\xi}=k \otimes_A \bH \qquad \mbox{where}\qquad 
\xi=\theta(v^2).\] 
Our aim is to establish a general existence result for ``canonical basic 
sets'' $\cB_{k,\xi}$ as in Definition~\ref{canbase}. In doing so, we will
also have to deal with the question of splitting fields for the algebra 
$\bH_{k,\xi}$. By Remark~\ref{semispec1}, the critical case that we have
to study is the case where $\xi$ is a root of unity in $k$. 

We assume throughout that P2--P8 and P15' hold for $W,L$; see 
(\ref{III}). The starting point are the following two results.

\begin{prop} \label{semi} Recall our assumptions on the characteristic of $k$.
Then $\bJ_k$ is semisimple and we have a unique bijection 
\[ \Irr(\bJ_{\QQ}) \stackrel{\sim}{\rightarrow} \Irr(\bJ_k),\qquad
M\mapsto M^k\]
such that $\dim M=\dim M^k$ and $\operatorname{trace}(t_w, M^k)=\theta 
\bigl(\operatorname{trace}(t_w,M)\bigr)$ for all $w \in W$.
\end{prop}

\begin{proof} (See \cite[2.5]{mykl}.) Our assumption on the characteristic
implies that $\theta(f_E)\neq 0$ for all $E \in \Irr_{\QQ}(W)$. The 
relevance of the numbers $f_E$ here is the fact that they are the 
coefficients in the expansion of the trace function $\mu_{\QQ}$ in 
terms of irreducible characters (see the discussion at the end of 
the previous section). Hence, by the same argument as in 
Theorem~\ref{semispec}, we conclude that $\bJ_k$ is split semisimple. 
The fact that we have a unique bijection between the simple modules 
with the above properties is a consequence of {\em Tits' Deformation 
Theorem}; see \cite[7.4.6]{ourbuch}.
\end{proof}

\begin{cor} \label{finalcor} We have a ``canonical'' bijection
\[ \Irr_{\QQ}(W) \stackrel{\sim}{\rightarrow} \Irr(\bJ_k), \quad 
\mbox{denoted} \quad E \mapsto E_\spadesuit^k.\]
For $E \in \Irr_{\QQ}(W)$, we have $\dim E=\dim E_\spadesuit^k$; 
furthermore, $\alpha_E=\ba(z)$ for any $z\in W$ such that 
$t_z.E_\spadesuit^k \neq \{0\}$.
\end{cor}

\begin{proof} We have a bijection $\Irr_{\QQ}(W) 
\stackrel{\sim}{\rightarrow} \Irr(\bJ_{\QQ})$ induced by the algebra 
isomorphism in  Example~\ref{exp2}. The correspondence $E \mapsto 
E_\spadesuit^k$ is obtained by composing that bijection with the one 
in Proposition~\ref{semi}. Now consider the statement concerning the
$a$-invariants. As in the previous section, 
\[ t_a:=\sum_{\atop{d \in W}{\ba(d)=a}} n_dt_d \in \bJ_k^a \qquad 
(a=0,1,2,\ldots)\] 
are central idempotents; furthermore, we have $1_J=\sum_{a\geq 0} t_a$ 
and $t_at_{a'}=0$ for $a\neq a'$. Hence, there exists a unique $a_0$ 
such that $t_{a_0}.E_\spadesuit^k=E_\spadesuit^k$. In particular, we have 
$\ba(z)=a_0$ for all $z \in W$ such that $t_z.E_\spadesuit^k \neq \{0\}$. 
In order to show that $a_0=\alpha_E$, it will now be sufficient to find 
some element $z_0\in W$ such that $\ba(z_0)=\alpha_E$ and 
$t_{z_0}.E_\spadesuit^k\neq \{0\}$. This is seen as follows. Since 
$\bJ_k$ is split semisimple, the character of a simple module cannot be 
identically zero. So there exists some $z_0\in W$ such that $\mbox{trace}
(t_{z_0}, E_\spadesuit^k)\neq 0$. If we had $\ba(z_0)\neq \alpha_E$, then 
$t_{z_0}.E_\spadesuit=\{0\}$ and so Proposition~\ref{semi} would yield
\[ \mbox{trace}(t_{z_0},E_\spadesuit)=\theta(\mbox{trace}(t_{z_0},
E_\spadesuit)) =0,\]
a contradiction. Hence, we must have $\ba(z_0)=\alpha_E$.
\end{proof}

In order to obtain a canonical basic set for $\bH_{k,\xi}$, we need 
to construct a map $\Irr(\bH_{k,\xi})\rightarrow \Irr_{\QQ}(W)$. By
Corollary~\ref{finalcor}, this is equivalent to constructing a map
$\Irr(\bH_{k,\xi})\rightarrow \Irr(\bJ_k)$. 

Following Lusztig (see the proof of \cite[Lemma~1.9]{cells3}), we 
attach an integer $\alpha_M$ to any $M \in \Irr(\bH_{k,\xi})$ by the 
requirement that
\begin{align*} 
c_w^\dagger.M &= \{0\} \quad \mbox{for all $w \in W$ with $\ba(w)>
\alpha_M$},\\ c_w^\dagger.M &\neq \{0\} \quad \mbox{for some 
$w \in W$ with $\ba(w)=\alpha_M$}.
\end{align*}
Now consider the $(\bJ_k, \bH_{k,\xi})$-bimodule $\bH_{k,\xi}^a$ 
where $a:=\alpha_M$. Then we obtain a (left)  $\bJ_k$-module 
\[ \tilde{M}:=\bH_{k,\xi}^a \otimes_{\bH_{k,\xi}} M.\]

\begin{lem} \label{lem1} Let $E \in \Irr_{\QQ}(W)$ and assume that
$E_\spadesuit^k\in \Irr(\bJ_k)$ is a composition factor of $\tilde{M}$.
Then $\alpha_E=\alpha_M$.
\end{lem}

\begin{proof} (See the proof of \cite[Cor.~3.6]{cells3}.) Let $z \in W$ be 
such that $\ba(z)=\alpha_E$ and $t_z.E_\spadesuit^k \neq 0$. Then we also 
have $t_z.\tilde{M} \neq \{0\}$ and so $t_z \star [c_w^\dagger] \neq 0$ for 
some $w \in W$ such that $\ba(w)=a$. The defining formula for ``$\star$'' 
and P8 now imply that $\alpha_E=\ba(z)=\ba(w)=\alpha_M$.
\end{proof}

Now there is a well-defined $k$-linear map $\pi\colon \tilde{M} 
\rightarrow M$ such that 
\[ \pi( [f] \otimes m )=f.m \qquad \mbox{for any $f \in \bH_{k,
\xi}^{\geq a}$ and $m \in M$}.\]
Note that, if $f',f \in \bH_{k,\xi}^{\geq a}$ are such that $[f]=[f']$
in $\bH_{k,\xi}^a$, then $f-f' \in \bH_{k,\xi}^{\geq a+1}$ acts as
zero on $M$ by the definition of $a=\alpha_M$.  

Now, we may also also consider $\tilde{M}$ as an $\bH_{k,\xi}$-module, 
using the map $\bJ_k\modA \rightarrow \bH_{k,\xi}\modA$ defined as follows. 
If $V$ is any $\bJ_k$-module, we can regard $V$ as an $\bH_{k,\xi}$ by 
composing the action of $\bJ_k$ on $V$ with the algebra homomorphism 
$\phi_{k,\xi} \colon \bH_{k,\xi} \rightarrow \bJ_k$. We denote that 
$\bH_{k,\xi}$-module by ${^*V}$. Note that, since $\bH_{k,\xi}$ is not 
necessarily semisimple, $\phi_{k,\xi}$ may not be an isomorphism and, 
consequently, ${^*V}$ may not be simple. 

Using this notation, let us consider the $\bH_{k,\xi}$-module
${^*\!\tilde{M}}$. Note that, by the compatibility in 
Theorem~\ref{bimod}, the left $\bH_{k,\xi}$-module structure on 
${^*\!\tilde{M}}$ is just the one coming from the natural left action 
of $\bH_{k,\xi}$ on $\bH_{k,\xi}^a$.

\begin{lem} \label{lem2} $\pi\colon {^*\!\tilde{M}} \rightarrow M$ is 
a surjective homomorphism of $\bH_{k,\xi}$-modules. If $M'\in 
\Irr(\bH_{k,\xi})$ is a composition factor of $\ker(\pi)$, then
$\alpha_{M'}<a_M$. In particular, $M$ occurs with multiplicity~$1$ as
a composition factor of ${^*\!\tilde{M}}$.
\end{lem}

\begin{proof} (See \cite[Lemma~1.9]{cells3}.)  Let $h \in \bH_{k,\xi}$,
$f \in \bH_{k,\xi}^{\geq a}$ and $m \in M$. Then  $hf \in 
\bH_{k,\xi}^{\geq a}$ by P4 and so 
\[h.\pi([f]\otimes m)=h.(f.m)=(hf).m=\pi([hf] \otimes m)=
\pi(h.([f]\otimes m)).\] 
Thus, $\pi$ is a homomorphism of $\bH_{k,\xi}$-modules. 
By the definition of $a=\alpha_M$, there exists some $m \in M$ and some
$w \in W$ such that $\ba(w)=a$ and $\pi([c_w^\dagger] \otimes m)=
c_w^\dagger.m \neq 0$.  Thus, $\pi\neq 0$. Since $M$ is simple, we conclude
that $\pi$ is surjective. 
As far as the composition factors of $\ker(\pi)$ are concerned,
it will be sufficient to show that $c_w^\dagger.\tilde{m}=0$ for any 
$\tilde{m}\in \ker(\pi)$ and $w \in W$ such that $\ba(w)\geq a$. So let 
$\tilde{m}:=\sum_i [f_i]\otimes m_i \in \ker(\pi)$ where 
$f_i \in \bH_{k,\xi}^{\geq a}$ and $m_i \in M$. Then
\begin{align*}
c_w^\dagger.\tilde{m}&=\Bigl(\sum_i [f_i]\otimes m_i\Bigr)=
\sum_i \Bigl(c_w^\dagger.[f_i]\Bigr) \otimes m_i\\
&=\sum_i [c_w^\dagger f_i] \otimes m_i=\sum_i [c_w^\dagger].f_i 
\otimes m_i \qquad \mbox{since $\ba(w)\geq a$}\\
&=c_w^\dagger.\Bigl(\sum_i f_i.m_i\Bigr) =c_w^\dagger.\pi(\tilde{m})=0,
\end{align*}
as desired.
\end{proof}

\begin{cor} \label{bs} Let $M\in \Irr(\bH_{k,\xi})$. Then there
is a unique $E=E(M)\in \Irr_{\QQ}(W)$ such that 
\begin{itemize}
\item[(a)] $E_\spadesuit^k$ is a composition factor of $\tilde{M}$ and 
\item[(b)] $M$ is a composition factor of ${^*\!E_\spadesuit^k}$.
\end{itemize}
The multiplicity of $M$ as a composition factor of  ${^*\!E_\spadesuit^k}$
is $1$ and $\alpha_M=\alpha_E$. For $M,M'\in \Irr_{\QQ}(W)$, we have 
$E(M) \cong E(M')$ if and only if $M \cong M'$. 
\end{cor}

\begin{proof} Consider a composition series of $\tilde{M}$ as a 
$\bJ_k$-module. Now, composing the action of $\bJ_k$ with the 
homomorphism $\phi_{k,\xi} \colon \bH_{k,\xi} \rightarrow \bJ_k$, 
the composition series becomes a filtration  of ${^*\!\tilde{M}}$ as 
an $\bH_{k,\xi}$-module, where the factors are not necessarily simple. 
Note, however, that all of these factors are of the form 
${^*\!E_\spadesuit^k}$, where $E\in \Irr_{\QQ}(W)$ and $E_\spadesuit^k$
is a composition factor of $\tilde{M}$. By Lemma~\ref{lem2}, there is a 
unique such factor ${^*\!E_\spadesuit^k}$ which has $M$ as a composition 
factor. This defines $E=E(M)$. 

It remains to show that $E(M)\not\cong E(M')$ if $M \not\cong M'$. Now,
by Lemma~\ref{lem1}, we have $\alpha_{E(M)}=\alpha_M$ and $\alpha_{E(M')}=
\alpha_{M'}$. Hence, if $\alpha_M \neq \alpha_{M'}$, then we certainly
have $E(M) \not\cong E(M')$. Now suppose that $\alpha_{M}=\alpha_{M'}$.
Since $E(M)_\spadesuit^k$ is a composition factor of $\tilde{M}$, 
Lemma~\ref{lem2} shows that $M$ is the unique composition factor of
$E(M)_\spadesuit^k$ which has $a$-invariant $\alpha_M$. Thus, $M'$ cannot
be a composition factor of $E(M)_\spadesuit^k$ and so $E(M)_\spadesuit^k
\not\cong E(M')_\spadesuit^k$. By Corollary~\ref{finalcor}, this shows
that $E(M)\not\cong E(M')$, as required.
\end{proof}

The above result defines an injective map
\[\Irr(\bH_{k,\xi})\hookrightarrow \Irr_{\QQ}(W),\qquad M\mapsto E(M).\]
We have $\alpha_M= \alpha_{E(M)}$ and $[E(M):M]=1$ for all $M\in 
\Irr(\bH_{k,\xi})$. Let 
\[ \cB_{k,\xi}:=\{E\in \Irr_{\QQ}(W) \mid E=E(M) \mbox{ for some
$M \in \Irr(\bH_{k,\xi})$}\}.\]
Now we can state the main result of this section.

\begin{thm}[Geck \protect{\cite{mykl}, \cite{my00} and 
Geck--Rouquier \cite{GeRo2}}] \label{basic} Recall that the characteristic
of $k$ is assumed to be zero or a prime which is $L$--good. Then
$\bH_{k,\xi}$ is split and $\cB_{k,\xi}$ is a canonical 
basic set in the sense of Definition~\ref{canbase}.
\end{thm}

First let us sketch the proof that $\bH_{k,\xi}$ is split. Let $k_0
\subseteq k$ be the field of fractions of the image of $\theta$. In order 
to prove that $\bH_{k,\xi}$ is split, it is certainly enough to
show that $\bH_{k_0,\xi}$ is split. By Remark~\ref{semispec1}, we may
assume that $\xi$ is a root of unity in $k_0$. Hence, $k_0$ is a finite
extension of the prime field of $k_0$. Consequently, $k_0$ is either a 
finite field or a finite extension of $\QQ$; in particular, $k_0$ is
perfect. Now Lemma~\ref{lem1}, Lemma~\ref{lem2} and Corollary~\ref{bs} 
also hold with $k$ replaced by $k_0$. Then we can argue as in
\cite[Theorem~3.5]{my00} to prove that $\bH_{k_0,\xi}$ is split, where
standard results on ``Schur indices'' and the ``multiplicity~$1$'' 
statement in Corollary~\ref{bs} play a crucial role. 

In order to show that $\cB_{k,\xi}$ is a canonical basic set, we 
can now assume that $k=k_0$.

We will need some standard results on projective indecomposable 
$\bH_{k,\xi}$-modules (PIM's for short). Let $M\in \Irr(\bH_{k,\xi})$. 
Then there is a PIM $P=P(M)$ (unique up to isomorphism) such that $M$ is 
the unique simple quotient of $P$. We can assume that $P=\bH_{k,\xi}e$ 
where $e\in \bH_{k,\xi}$ is a primitive idempotent. Let us consider the 
filtration of $\bH_{k,\xi}$ defined by the $\ba$-function:
\[\{0\}=\bH_{k,\xi}^{\geq N+1}\subseteq \bH_{k,\xi}^{\geq N} \subseteq 
\bH_{k,\xi}^{\geq N-1} \subseteq \cdots \subseteq \bH_{K,\xi}^{\geq 0}=
\bH_{k,\xi}\]
where $N=\max\{\ba(w) \mid w \in W\}$. Multiplying on the right by $e$
and setting $P^{\geq i}:=\bH_{k,\xi}^{\geq i}e$, we obtain a filtration
\[ \{0\}=P^{\geq N+1}\subseteq P^{\geq N} \subseteq P^{\geq N-1}
\subseteq \cdots \subseteq P^{\geq 0}=P.\]
Now fix $i$ and consider the canonical exact sequence
\[ \{0\} \rightarrow \bH_{k,\xi}^{\geq i+1} \rightarrow 
\bH_{k,\xi}^{\geq i} \rightarrow \bH_{k,\xi}^i \rightarrow \{0\}.\]
Multiplying with $e$ yields a sequence
\[ \{0\} \rightarrow \bH_{k,\xi}^{\geq i+1}\,e \rightarrow 
\bH_{k,\xi}^{\geq i}\,e\rightarrow\bH_{k,\xi}^i\,e\rightarrow \{0\}.\]
which is also exact. (To see this, also work with the idempotent $1-e$.)
It follows that we have a canonical isomorphism of $\bH_{k,\xi}$-modules
\[ P^i \cong \bH_{k,\xi}^i\,e \qquad \mbox{for $i=0,1,2,\ldots$}.\]
Now consider the $(\bJ_k,\bH_{k,\xi})$-bimodule structure of 
$\bH_{k,\xi}^i$; see Theorem~\ref{bimod}. Since $\bH_{k,\zeta}^i\,e$ is 
obtained by right multiplication with some element of $\bH_{k,\xi}$, it is 
clear that $\bH_{k,\xi}^i\,e$ is a $\bJ_k$-submodule of $\bH_{k,\xi}^i$. 
Hence, we can write 
\[ P^i \cong \bH_{k,\xi}^i\,e=V_1\oplus \cdots \oplus V_r \qquad 
\mbox{where} \quad V_1,\ldots,V_r \in  \Irr(\bJ_k).\]
Note that, since $P^i$ is isomorphic to a submodule of $\bH_{k,\xi}^i$, 
this implies that $\alpha_{V_1}=\cdots =\alpha_{V_r}=i$. This leads 
us to the following result.

\begin{prop}[Brauer reciprocity] \label{lem3} Let $M\in 
\Irr(\bH_{k,\xi})$ and $P=P(M)$ be a corresponding {\em PIM}. For any 
$E \in \Irr(\bJ_k)$ such that $\alpha_E=i$, we have 
\begin{align*}
[E:M]&=\mbox{multiplicity of $M$ as a composition factor of 
${^*\!E_{\spadesuit}^k}$}\\
&=\mbox{multiplicity of ${E_\spadesuit^k}$ as a direct summand
in $P^i$},
\end{align*}
where $P^i$ is considered as a left $\bJ_{k}$-module as explained above.
\end{prop}

The proof of the first equality uses an abstract characterization of 
the multiplicities $[E:M]$ in terms of decomposition maps between the 
appropriate Grothendieck groups; see \cite[\S 2]{mybourb} and 
\cite[2.5]{mykl}. Then the second statement can be reduced to the 
classical Brauer reciprocity in the modular representation theory of 
finite groups and associative algebras; see \cite{GeRo2} for the 
details. 

Now we can complete the proof of Theorem~\ref{basic}. Let $E \in 
\Irr_{\QQ}(W)$ and $M \in \Irr(\bH_{k,\xi})$. First we show:
\begin{equation*}
[E:M] \neq 0 \qquad \Rightarrow \qquad \alpha_M \leq \alpha_E.\tag{$*$}
\end{equation*}
This is seen as follows. Let $w \in W$ be such that $\ba(w)=\alpha_M$ and 
$c_w^\dagger.M\neq \{0\}$. Then $c_w^\dagger$ also acts non-zero
on $E_\spadesuit^k$. By the definition of that action, this means
that $\phi_{k,\xi}(c_w^\dagger).E_\spadesuit^k \neq \{0\}$. So there 
exist some $z \in W$ and $d \in \cD$ such that $h_{w,d,z}\neq 0$ and 
$t_z.E_\spadesuit^k \neq \{0\}$. Thus, we necessarily have $\alpha_E=
\ba(z)$ by Corollary~\ref{finalcor}. Now P4 implies that $\ba(w)\leq 
\ba(z)$, as required.

Next we show that 
\[ \breve{\alpha}_M=\alpha_M \qquad \mbox{for any $M\in 
\Irr(\bH_{k,\xi})$}.\]
Indeed, let $E \in \Irr_{\QQ}(W)$ be such that $[E:M]\neq 0$ and
$\breve{\alpha}_M=\alpha_E$. Then ($*$) shows that $\alpha_M \leq 
\breve{\alpha}_M$. On the other hand, $[E(M):M]=1$ by 
Corollary~\ref{bs}. Hence $\breve{\alpha}_M \leq \alpha_{E(M)}=
\alpha_M$, where the last equality holds by Lemma~\ref{lem1}.

Finally, we must show that, for any $E \in \Irr_{\QQ}(W)$ and
$M\in \Irr(\bH_{k,\xi})$, we have the implication:
\[ [E:M]\neq 0 \quad \mbox{and} \quad \alpha_M=\alpha_E \quad
\Rightarrow \quad E=E(M).\]
To prove this, we consider a PIM $P=P(M)$. Let $a:=\min \{i \geq 0 
\mid P^i \neq \{0\}\}$. Then $P=P^{a}$ and there is a natural
surjective homomorphism $P \rightarrow P^{a}$ of $\bH_{k,\xi}$-modules. 
As before, considering $P^a$ as a left $\bJ_k$-module, we can write 
\[ P^a \cong \bH_{k,\xi}^a\,e=V_1\oplus \cdots \oplus V_r \qquad 
\mbox{where} \quad V_1,\ldots,V_r \in  \Irr(\bJ_k).\]
The compatibility in Theorem~\ref{bimod} now implies that we have an 
isomomorphism
\[ P^a \cong {^*V_1}\oplus \cdots \oplus {^*V_r} \qquad 
\mbox{(as left $\bH_{k,\xi}$-modules)}.\]
Since $M$ is the unique simple quotient of $P$, we deduce that $M$ also 
is the unique simple quotient of $P^{a}$. In particular, $P^a$ is 
indecomposable. Hence we must have $r=1$. Let us write $V_1=
E_{1,\spadesuit}^k$ where $E_1\in \Irr_{\QQ}(W)$. By 
Proposition~\ref{lem3}, we have $[E_1:M]\neq 0$ and so $\alpha_M
\leq \alpha_{E_1}=a$, using ($*$). 

Now consider $E$ such that $[E:M]\neq 0$ and $\alpha_M=\alpha_E$. By 
Proposition~\ref{lem3} and the definition of $a$, we have $\alpha_E\geq a$.
Furthermore, if $\alpha_E=a$, then we necessarily have $E\cong E_1$. 
Now, since $\alpha_M=\alpha_E$, we cannot have $\alpha_E>a$. Hence
$a=\alpha_E$ and so $E\cong E_1$. 

Finally, we do have $[E(M):M]\neq 0$ and $\alpha_M=\alpha_{E(M)}$. Hence, 
we can apply the previous argument to $E(M)$ and this yields $E(M) \cong
E_1\cong E$, as required. This completes the proof of Theorem~\ref{basic}.

\begin{exmp} \label{excep} Assume that $L$ is a positive multiple of
the length function. Then, as mentioned in (\ref{III}), the properties 
P1--P15 are known to hold.  Hence, the above results guarantee the 
existence of a canonical basic set for $\bH_{k,\xi}$ if the characteristic 
of $k$ is zero of a good prime for $W$. Explicit tables for these canonical 
basic sets for the exceptional types $G_2$, $F_4$, $E_6$, $E_7$, $E_8$ 
are given by Jacon \cite[\S 3.3]{Jac0}, using the known information on 
decomposition numbers (see \cite{GeLu}, \cite{mye6}, \cite{habil} 
and M\"uller \cite{muell}).
\end{exmp}

\begin{exmp} \label{canbn} Let $W_n$ be a Coxeter group of type $B_n$,
as in Example~\ref{schurbn}. Recall that a weight function $L$ on $W_n$
is specified by two integers $a,b\geq 0$. Let us denote the
corresponding specialized algebra simply by $H_n$.

Now, as in the case of the symmetric group (see Example~\ref{expsym}), there
is a theory of ``Specht modules'' for $H_n$; see Dipper--James--Murphy 
\cite{DJM3}. Thus, for any bipartition $\ulambda\vdash n$, we have a 
corresponding Specht module $S_k^{\ulambda}\in H_n\modA$. (Actually, we 
shall denote by $S_k^{\ulambda}$ the module that is labelled by 
$(\lambda_{(2)}^*, \lambda_{(1)}^*)$ in \cite{DJM3}, where the star denotes 
the conjugate partition.) As before, this module has the property that
\[ [E^\ulambda:M]=\mbox{multiplicity of $M$ as a composition factor of
$S_k^{\ulambda}$},\]
for any $M \in \Irr(H_n)$. Furthermore, there is an $H_n$-equivariant 
symmetric bilinear form on $S_k^{\ulambda}$ and, taking quotients by the 
radical, we obtain modules
\[ D^{\ulambda}:=S_k^{\ulambda}/ \mbox{rad}(S_k^{\ulambda}) \in
H_n\modA.\]
Let $\Lambda_{2,n}^{\circ}:=\{\ulambda\vdash n\mid D^\ulambda\neq \{0\}\}$. 
Then, by \cite[\S 6]{DJM3}, we have
\[\Irr(H_n)=\{D^{\ulambda}\mid \ulambda \in \Lambda_{2,n}^{\circ}\}.\]
Furthermore, the decomposition numbers $[S^\ulambda:D^\umu]$ satisfy
conditions similar to those in Example~\ref{expsym}, where we have
to consider the dominance order on bipartitions. 

Is it possible to interprete $\Lambda_{2,n}^\circ$ in terms of canonical
basic sets? The answer is ``yes'', but we need to choose the parameters
$a,b$ appropriately. So let us assume that $b>(n-1)a>0$, that is, we
are in the ``asymptotic case'' defined in Example~\ref{schurbn}. All primes 
are $L$-good in this case. Furthermore, as remarked in (\ref{III}), the 
properties P2--P8 and P15' hold. Hence, Theorem~\ref{basic} shows that we 
have a canonical basic set $\cB_{k,\xi}$. Using a compatibility of the
dominance order on bipartitions with the invariants $\alpha_E$ for this case
(see \cite[Cor.~5.5]{GeIa05}), one can use exactly the same arguments
as in Example~\ref{expsym} to show that
\[ \cB_{k,\xi}=\{E^\ulambda \mid \ulambda \in \Lambda_{2,n}^\circ\}
\qquad \mbox{(assuming that $b>(n-1)a>0$)}.\]
However, contrary to the situation in the symmetric group, a combinatorial 
description of $\Lambda_{2,n}^\circ$ is much harder to obtain. Dipper, 
James and Murphy already dealt with the following cases. As in 
Example~\ref{expsym}, let us set 
\[e:=\min \{i \geq 2 \mid 1+\xi^a+\xi^{2a}+\cdots + \xi^{(i-1)a}=0\}.\]
(If no such $i$ exists, we set $e=\infty$.)  Following Dipper--James 
\cite[4.4]{DJ2}, let 
\[ f_n(a,b):=\prod_{i=-(n-1)}^{n-1} (\xi^b+\xi^{ai}).\]
Now there is a distinction between the cases where $f_n(a,b)$ is zero or 
not.

\medskip
$\bullet$ {\bf The case $f_n(a,b)\neq 0$.} Then we have 
\[ \Lambda_{2,n}^\circ=\{\ulambda\vdash n \mid \mbox{$\lambda_{(1)}$
and $\lambda_{(2)}$ are $e$-regular}\}.\]
Indeed, by \cite[Theorem~6.9]{DJM3}, we have the inclusion ``$\subseteq$''.
On the other hand, Dipper--James \cite[4.17 and 5.3]{DJ2} showed that 
there is a Morita equivalence
\[H_n\modA \quad \cong \quad \Bigl(\bigoplus_{r=0}^n H_k(\Sym_r,\xi^a)
\otimes_k H_k(\Sym_{n-r},\xi^a)\Bigr)\modA;\]
In particular, we have a bijection
\[ \Irr(H_n) \stackrel{\sim}{\rightarrow} \coprod_{r=0}^n
\Bigl(\Irr(H_k(\Sym_r,\xi^a)) \times \Irr(H_k(\Sym_{n-r},\xi^a))\Bigr).\]
Consequently, there are as many simple $H_n$-modules as there are
bipartitions $\ulambda$ such that $\lambda_{(1)}$ and
$\lambda_{(2)}$ are $e$-regular. This yields the inclusion ``$\supseteq$''.

In \cite{GeJa}, we show that, for {\em any} values of $a,b\geq 0$ such that
$f_n(a,b)\neq 0$, the set $\Lambda_{2,n}^\circ$ is a canonical 
basic set in the sense of Definition~\ref{canbase}, provided the
characteristic of $k$ is not~$2$.

\medskip
$\bullet$ {\bf The case $\xi^a=1$.} If $\xi^b\neq -1$, then $f_n(a,b)\neq 0$
and we can apply the previous case. So let us now assume that $\xi^b=-1$. 
The special feature of this case is that, by \cite[Remark~5.4]{DJ2}, the
simple $H_n$-modules are obtained by extending (in a unique way) the 
simple modules of the parabolic subalgebra $k\Sym_n=\langle T_{s_1},
\ldots,T_{s_{n-1}} \rangle_k$ to $H_n$. By \cite[Theorem~7.3]{DJM3}, 
we have
\[ \Lambda_{2,n}^\circ=\{\ulambda \vdash n \mid \mbox{$\lambda_{(1)}$
is $e$-regular and $\lambda_{(2)}=\varnothing$}\}.\]
The general case where $f_n(a,b)=0$ will be discussed in Section~8. 
\end{exmp}

The above discussion shows that the set $\Lambda_{2,n}^\circ$ arising
from the theory of Specht modules can be interpreted as a {\em canonical 
basic set} assuming that $b$ is large with respect to~$a$. It would certainly 
be interesting to know the canonical basic sets for other choices of 
$a$ and $b$ as well, for example, the equal parameter case where $a=b>0$, 
or the case where $b=0$ and $a>0$ (which is relevant to type $D_n$). 
This will be discussed in Section~8.

\section{The Fock space and canonical bases} \label{MGsec6}

We now take an excursion to another area of representation theory: the 
theory of canonical bases for highest weight modules of quantized 
enveloping algebras, as developed by Kashiwara and Lusztig. Of course, it 
is not the place here to give any reasonably detailed introduction (we 
refer the reader to Kashiwara \cite{Kash} or Lusztig \cite{Lucan}), but 
what we can do is to explain some of the main ideas behind a deep 
combinatorial construction arising from that theory, namely, the 
{\em crystal graph} of a highest weight module. (In the following 
section, we will see applications to Iwahori--Hecke algebras.) 

The starting point is the following idea. Recall that the simple
$\QQ \Sym_n$-modules are parametrized by the set of all partitions 
of $n$. More generally, the simple modules of a Coxeter group of type 
$B_n$ are parametrized by the set of all pairs of partitions such that 
the total sum of all parts equals $n$. Even more generally, the simple 
modules of the semidirect product 
\[ G_{r,n}:=(\ZZ/r\ZZ)^n \rtimes \Sym_n \quad \mbox{(where $r\geq 1$)}\]
are parametrized by the set $\Pi_{r,n}$ of all $r$-tuples
$\ulambda=(\lambda_{(1)},\ldots,\lambda_{(r)})$ where each
$\lambda_{(i)}$ is a partition of some $a_i\geq 0$ and where
$n=a_1+\cdots+a_r$. We shall now fix $r$ and consider the simple
modules for all the groups $G_{r,n}$ ($n\geq 1$) at the same time. 
For this purpose, let $\fF_{r,n}$ be the ${\CC}$-vector space with 
basis $\Pi_{r,n}$. The ``Fock space'' is defined as the direct sum of 
all these spaces:
\[\fF^{(r)}:=\bigoplus_{n\geq 0} \fF_{r,n};\]
where $\fF_{0,n}$ is the $1$-dimensional space with basis $\Pi_{0,n}=\{
\uvar:=(\varnothing,\ldots,\varnothing)\}$. The point
of collecting all these spaces into one object is that there are natural 
operators sending $\fF_{r,n}$ into $\fF_{r,n+1}$, and vice versa. Indeed,
the {\em branching rule} for the induction of representations from
$G_{r,n}$ to $G_{r,n+1}$ gives rise to a linear map 
\[ \text{ind}\colon \fF^{(r)} \rightarrow \fF^{(r)}, \qquad 
\ulambda \mapsto \sum_{\umu} \umu 
\qquad (\ulambda \in \Pi_{r,n}),\]
where the sum runs over all $\umu \in \Pi_{r,n+1}$ which
can be obtained from $\ulambda$ by increasing exactly one
part by~$1$. Similarly, the restriction of representations from
$G_{r,n}$ to $G_{r,n-1}$ gives rise to a linear map 
\[ \text{res}\colon \fF^{(r)} \rightarrow \fF^{(r)}, \qquad 
\ulambda \mapsto \sum_{\umu} \umu 
\qquad (\ulambda \in \Pi_{r,n}),\]
where the sum runs over all $\umu \in \Pi_{r,n-1}$ which
can be obtained from $\ulambda$ by decreasing exactly one
part by~$1$. Now, fixing a positive integer $l\geq 2$ and a set
of parameters
\begin{equation*}
\bu=\{ u_1,\ldots,u_r\} \qquad \mbox{where} \qquad 
u_i\in \ZZ,
\end{equation*}
the operators $\mbox{ind}$ and $\mbox{res}$ can be refined into sums of
linear operators
\[ \mbox{ind}=\sum_{i=0}^{l-1} f_i \quad \mbox{and}\quad\mbox{res}=
\sum_{i=0}^{l-1} e_i\quad \mbox{where}\quad e_i,f_i\colon \fF^{(r)}
\rightarrow \fF^{(r)}.\]
The definition of these refined operators (which we will give further below)
also has an interpretation in terms of representations; see 
\cite[\S 12.1]{Ar3}. Now it is a remarkable fact that these operators 
turn out to satisfy the {\em Serre relations} for the affine 
Kac--Moody algebra $\hat{\fsl}_l$. This algebra is defined as
\[ \fg=\hat{\fsl}_l:=\bigl({\CC}[t,t^{-1}] \otimes_{\CC} 
{\fsl}_l\bigr) \oplus \CC c \oplus \CC d,\]
where $t$ is an indeterminate and ${\fsl}_l$ is the usual  Lie algebra 
of complex $l \times l$-matrices with trace zero and product 
$[X,Y]=XY-YX$. Thus, as a $\CC$-vector space, $\fg$ is spanned by 
$\{t^n \otimes X \mid n \in \ZZ, X \in \fsl_l\} \cup \{c,d\}$; the 
Lie product is given by 
\begin{gather*}
[t^n \otimes X,t^m \otimes Y]=t^{n+m} \otimes [X,Y]+ \mbox{trace}(XY)\,n
\, \delta_{0,n+m} \,c,  \\ [c,t^n\otimes X]=0, \quad [d,t^n \otimes X]=n
\,t^n \otimes X, \quad [c,d]=0.
\end{gather*}
In dealing with representations of $\fg$, it will be convenient to work 
with a presentation of the corresponding universal enveloping algebra 
$U(\fg)$.  We set
\begin{alignat*}{2}
e_0&:=t \otimes E_{l1}, \qquad &e_i &:=1 \otimes E_{i,i+1} \quad 
\mbox{($1\leq i\leq l-1$)},\\ f_0&:=t^{-1} \otimes E_{1l}, \qquad 
&f_i&:=1 \otimes E_{i+1,i} \quad \mbox{($1\leq i\leq l-1$)},\\
h_0&:=-\sum_{i=1}^{l-1} h_i +c,\qquad &h_i&:=1 \otimes 
(E_{i,i}-E_{i+1,i+1}) \quad \mbox{($1\leq i\leq l-1$)},
\end{alignat*}
where the $E_{i,j}$ are the usual matrix units in $\fsl_l$. Let 
\[ \fh:=\langle d,h_0,h_1,\ldots,h_{l-1}\rangle_{\CC}\subseteq \fg\]
be the Cartan subalgebra. We have $[h,h']=0$ for all $h,h'\in \fh$. 
For $0\leq i \leq l-1$, we define linear forms $\alpha_i\in \fh^*=
\Hom_{\CC}(\fh,\CC)$ by 
\[ \alpha_i(h_j)=\left\{\begin{array}{rl} 2 & \qquad \mbox{if $i=j$},\\
-1 & \qquad \mbox{if $i-j \equiv \pm 1 \bmod l$},\\ 0 & \qquad 
\mbox{otherwise}; \end{array}\right.\]
and $\alpha_i(d)=\delta_{i0}$. By the classification in Kac 
\cite{Kac}, the matrix 
\[A=\bigl(\alpha_i(h_j)\bigr)_{0\leq i,j\leq l-1}\]
is the Cartan matrix of affine type $A_{l-1}^{(1)}$. With this notation, 
$U(\fg)$ is the associative algebra (with $1$) generated by elements
\[\{e_i,f_i,h_i \mid  0 \leq i \leq l-1\} \cup\{d\},\]
and subject to the following defining relations:
\begin{gather*}
h_je_i-e_ih_j=\alpha_{i}(h_j) e_i,\qquad h_jf_i -f_ih_j= -\alpha_i(h_j) 
f_i,\\ de_i-e_id=\delta_{i0}e_i,\qquad df_i-f_id=-\delta_{i0}f_i,\\
e_if_j-f_j e_i=\delta_{ij}h_i, \qquad h_ih_j=h_jh_i,\qquad h_id=dh_i,\\ 
e_i e_j=e_je_i \quad \mbox{and}\quad f_i f_j=f_jf_i \qquad 
\mbox{(if $i-j \not\equiv \pm 1 \bmod l$)}.
\end{gather*}
Furthermore, if $l \geq 3$ and $i-j \equiv \pm 1 \bmod l$, then 
\[e_i^2 e_j-2e_i e_j e_i+e_j e_i^2=0 \quad \mbox{and}\quad  
f_i^2 f_j-2f_i f_j f_i+f_j f_i^2=0;\]
if $l=2$ and $i\neq j$, then
\begin{align*}
e_i^3  e_j-3(e_i^2  e_j  e_i+e_i  e_j  e_i^2) -e_j  e_i^3=0,\\
f_i^3  f_j-3(f_i^2  f_j  f_i+ f_i  f_j  f_i^2) -f_j  f_i^3=0.
\end{align*}
For a sketch of proof and further references, see Ariki \cite[\S 3.2]{Ar3}.
We will now define a linear action of $U(\fg)$ on $\fF^{(r)}$. For this 
purpose, we need some further notation. Let $\ulambda=
(\lambda_{(1)}, \ldots,\lambda_{(r)})\in \Pi_{r,n}$ and write
\[ \lambda_{(c)}=(\lambda_{(c),1} \geq \lambda_{(c),2} \geq \cdots 
\geq 0) \qquad \mbox{for} \qquad c=1,\ldots,r.\]
The diagram of $\ulambda$ is defined as the set 
\[ [\ulambda]:=\{(a,b,c) \mid 1 \leq c \leq r, \;1 \leq b 
\leq \lambda_{(c),a} \mbox{ for $a=1,2,\ldots$}\}.\]
For any ``node'' $\gamma=(a,b,c)\in [\ulambda]$, we set 
\[ \mbox{res}_l(\gamma):=(b-a+u_c)  \quad \bmod l\]
and call this the $l$-residue of $\gamma$ with respect to the 
parameters $\bu$. If $\mbox{res}_l(\gamma)=i$, we say that 
$\gamma$ is an $i$-node of $\ulambda$. Let 
\[ W_i(\ulambda) :=\mbox{number of $i$-nodes of 
$\ulambda$}.\]
Now suppose that $\ulambda \in \Pi_{r,n}$ and 
$\umu\in \Pi_{r,n+1}$ for some $n\geq 0$. We write
\[ \gamma=\umu/\ulambda \qquad \mbox{if}\qquad
[\ulambda] \subset [\umu] \quad \mbox{and}\quad
[\umu]=[\ulambda] \cup \{\gamma\};\]
Then we call $\gamma$ an addable node for $\ulambda$ or a 
removable node for $\umu$. Let 
\begin{align*}
A_i(\ulambda) &:=\mbox{set of addable $i$-nodes for 
$\ulambda$},\\
R_i(\umu) &:=\mbox{set of removable $i$-nodes for 
$\umu$}.
\end{align*}
Now let $0 \leq i \leq l-1$. We define a linear operation of $e_i$ on
$\fF^{(r)}$ by
\[ e_i.\ulambda=\sum_{\umu} \, \umu 
\qquad (\ulambda\in \Pi_{r,n})\]
where the sum runs over all $\umu \in \Pi_{r,n-1}$ such that 
$[\umu] \subset [\ulambda]$ and $\mbox{res}_l
(\ulambda/\umu) \equiv i\bmod l$.
Similarly, we define a linear operation of $f_i$ on $\fF^{(r)}$ by
\[ f_i.\ulambda=\sum_{\umu}\umu \qquad 
(\ulambda\in \Pi_{r,n})\]
where the sum runs over all $\umu \in \Pi_{r,n+1}$ such that 
$[\ulambda] \subset [\umu]$ and $\mbox{res}_l(
\umu/\ulambda)\equiv i\bmod l$. Note that, indeed,
we have $\mbox{res}=\sum_i e_i$ and $\mbox{ind}=\sum_i f_i$. 

Next, we define a linear operation of $h_i$ and $d$ on $\fF^{(r)}$ by
\[ h_i.\ulambda=N_i(\ulambda)\,\ulambda 
\qquad \mbox{and} \qquad d.\ulambda=-W_0(\ulambda)\,
\ulambda,\]
where $N_i(\ulambda)=|A_i(\ulambda)|-
|R_i(\ulambda)|$.
The following result (extended to a ``quantized'' version) is originally
due to Hayashi for $r=1$; see \cite[Chap.~10]{Ar3} for the proof and 
further historical remarks. 

\begin{prop}[See Ariki \protect{\cite[Lemma~13.35]{Ar3}}] \label{hash1}
Via the above maps, $\fF^{(r)}$ becomes an integrable $U(\fg)$-module. 
Let $M^{(r)}(\uvar)$ be the submodule generated by 
$\uvar\in \fF_{r,0}$. Then $M^{(r)}
(\uvar)$ is a highest weight module with highest 
weight given by $\sum_{i=1}^r \Lambda_{\gamma_i}$ where $\gamma_i \in 
\{0,1,\ldots,l-1\}$ is such that $\gamma_i\equiv u_i \bmod l$.
\end{prop}

Here, the fundamental weights $\Lambda_0,\Lambda_1, \ldots,\Lambda_{l-1}$ 
for $\fg$ are the elements of $\fh^*$ defined by $\Lambda_i(h_j)=\delta_{ij}$
and $\Lambda_i(d)=0$. We refer to Kac \cite{Kac} for the general theory
of integrable modules for Kac--Moody algebras.

In order to speak about the {\em canonical basis} of 
$M^{(r)}(\uvar)$, we need to consider a ``deformed'' 
version of $U(\fg)$, that is, the quantized universal enveloping algebra 
$U_v(\fg)$, where $v$ is an indeterminate. Following Ariki 
\cite[\S 3.3]{Ar3}, $U_v(\fg)$ is the associative ${\CC}(v)$-algebra
(with $1$), generated by elements
\[ \{E_i,F_i,K_i,K_i^{-1} \mid 0 \leq i \leq l-1\} \cup \{D,D^{-1}\},\]
subject to the following defining relations:
\begin{gather*}
K_jE_iK_j^{-1}=v^{\alpha_{i}(h_j)} E_i,\qquad K_jF_i
K_j^{-1}= v^{-\alpha_{i}(h_j)} F_i,\\ DE_iD^{-1}=
v^{\delta_{0i}}E_i,\qquad DF_iD^{-1}=v^{-\delta_{0i}}F_i,\\
E_i  F_j-F_j E_i=\delta_{ij}\frac{K_i-K_i^{-1}}{v-v^{-1}},
\quad K_iK_j=K_jK_i,\quad DK_i=K_iD, \\
E_i  E_j=E_j E_i \quad \mbox{and}\quad F_i  
F_j=F_j F_i \qquad \mbox{(if $i-j \not\equiv \pm 1 \bmod l$)}.
\end{gather*}
Furthermore, if $l \geq 3$ and $i-j \equiv \pm 1 \bmod l$, then 
\begin{align*}
E_i^2 E_j-(v+v^{-1})E_i E_j E_i+E_j E_i^2&=0,\\ 
F_i^2 F_j-(v+v^{-1})F_i F_j F_i+F_j F_i^2&=0;
\end{align*}
if $l=2$ and $i\neq j$, then
\begin{align*}
E_i^3  E_j-(v^2+1+v^{-2})(E_i^2  E_j  E_i+E_i  
E_j  E_i^2) -E_j  E_i^3=0,\\
F_i^3  F_j-(v^2+1+v^{-2})(F_i^2  F_j  F_i+ F_i  
F_j  F_i^2) -F_j  F_i^3=0.
\end{align*}
We can extend the action of $U(\fg)$ on $\fF^{(r)}$ to an action of 
$U_v(\fg)$ on 
\[\fF_v^{(r)}:={\CC}(v) \otimes_{\CC} \fF^{(r)}.\]
This will depend on the choice of a total order on nodes. 

\begin{defn}[Foda et al. \protect{\cite[p.~331]{FLOTW}}] \label{above} 
We say that the node $\gamma=(a,b,c)$ is ``above'' the node $\gamma'=(a',
b',c')$ if 
\[ b-a+u_c<b'-a'+u_{c'} \quad \mbox{or if} \quad b-a+u_c=b'-a'+u_{c'}
\mbox{ and } c'<c.\]
\end{defn}

Now let $0 \leq i \leq l-1$. We define a linear operation of $E_i$
on  $\fF_v^{(r)}$ by
\[ E_i.\ulambda=\sum_{\umu} v^{-N_i^a
(\ulambda/\umu)} \, \umu 
\qquad (\ulambda\in \Lambda_{r,n})\]
where the sum runs over all $\umu \in \Pi_{r,n-1}$ such that 
$[\ulambda]=[\umu]\cup\{\gamma\}$, 
$\mbox{res}_l(\gamma)\equiv i\bmod l$ and 
\[N_i^a(\ulambda/\umu)=
|\{\gamma'\in A_i(\umu) \mid \gamma' \mbox{ above } \gamma\}|-
|\{\gamma'\in R_i(\ulambda) \mid \gamma' \mbox{ above } 
\gamma\}|.\]
Similarly, we define a linear operation of $F_i$ on $\fF_v^{(r)}$ by
\[ F_i.\ulambda=\sum_{\umu} v^{N_i^b
(\umu/\ulambda)} \, \umu \qquad 
(\ulambda\in \Lambda_{r,n})\]
where the sum runs over all $\umu \in \Pi_{r,n+1}$ such that 
$[\umu]=[\ulambda]\cup\{\gamma\}$, $\mbox{res}_l
(\gamma)\equiv i \bmod l$ and 
\[N_i^b(\umu/\ulambda)=
| \{\gamma'\in A_i(\ulambda) \mid \gamma \mbox{ above } 
\gamma'\}| -|\{\gamma'\in R_i(\umu) \mid \gamma \mbox{ above }
\gamma'\}|.\]
Next, we define linear operations of $K_i$ and $D$ on $\fF_v^{(r)}$ by 
\[ K_i.\ulambda=v^{N_i^b(\ulambda)} \, 
\ulambda \qquad \mbox{and}\qquad D.\ulambda=
v^{-W_0(\ulambda)}\ulambda.\]

\begin{thm}[Jimbo et al. \protect{\cite{JMMO}}, Foda et al. 
\protect{\cite{FLOTW}}, Uglov \protect{\cite{Ug}}] \label{jmmo} 
Via the above maps, $\fF_v^{(r)}$ becomes an integrable 
$U_v(\fg)$-module. The submodule $M_v^{(r)}(\uvar)$ 
generated by $\uvar$ is a highest weight module with 
highest weight $\sum_{i=1}^r \Lambda_{\gamma_i}$.
\end{thm}

\begin{rem} \label{conven} The above definition of a $v$-deformed action
of $U(\fg)$ on $\fF^{(r)}$ is not the only possible one. Ariki 
\cite[Theorem~10.10]{Ar3} considers an action of $U_v(\fg)$ on 
$\fF_v^{(r)}$ which is given by exactly the same formulae as above, but
where the exponents $N_i^a$ and $N_i^b$ are computed with respect to the 
following order on nodes: let $\gamma=(a,b,c)$ and $\gamma=(a',b',c')$; 
then
\begin{equation*}
\gamma \prec \gamma' \qquad \stackrel{\text{def}}{\Leftrightarrow} 
\qquad c'<c \quad \mbox{or if} \quad c=c' \mbox{ and } a'<a.
\end{equation*}
Note that $\prec$ does not depend on the parameters $\bu$~! A relation 
between this order and the one in Definition~\ref{above} can be 
established as follows. Suppose we are only interested in
$\bigoplus_{0 \leq k \leq n} \fF_{r,k}$ for some fixed $n$. Then choose 
the integers $u_1,\ldots,u_r$ such that 
\[ u_1>u_2 > \cdots > u_r>0 \qquad \mbox{where}\qquad u_i-u_{i+1}>n-1
\mbox{ for all $i$}.\]
Let $\gamma=(a,b,c)$ and $\gamma'=(a',b',c')$ for some partitions
in $\Pi_{r,k}$ where $0 \leq k \leq n$. One easily checks that 
$\gamma$ lies above $\gamma'$ if and only if $\gamma\prec \gamma'$.
\end{rem}

Since $M_v^{(r)}(\uvar)$ is an integrable highest weight 
module, the general theory of Kashiwara and Lusztig provides us with a 
canonical basis of $M_v^{(r)}(\uvar)$. In Lusztig's 
setting, this basis is given by 
\[ \bB(\uvar)=\{b.\uvar \mid b 
\in \bB\}\setminus\{0\},\]
where $\bB$ is the {\em canonical basis}  of the subalgebra $U_v^-(\fg)
\subseteq U_v(\fg)$ generated by $F_0,F_1,\ldots,F_{l-1}$; see  
\cite[Theorem~7.3, Prop.~9.1]{Ar3} and the references there. Setting 
$v=1$, we obtain the canonical basis, denoted $\bB_1(\uvar)$, of 
the $U(\fg)$-module $M^{(r)}(\uvar)$. 

In order to make an efficient use of this result, we need a good
parametrization of the canonical basis. This is provided by Kashiwara's 
{\em crystal basis}, which is obtained by specializing $v$ to $0$. A 
modified action of the generators $E_i$ and $F_i$ sends the crystal basis 
to itself, and this gives the crystal basis the structure of a colored 
oriented graph. In this way, many properties of the representation
$M^{(r)}(\uvar)$ can be studied using combinatorial 
properties of that graph. Let us briefly discuss how all this works, 
where we follow Kashiwara \cite{Kash}. The starting point is the 
observation that the subalgebra
\[ U_{v,i}:=\langle E_i,F_i,K_i^{\pm 1}\rangle\subseteq U_v(\fg)
\qquad (0 \leq i \leq l-1)\]
is isomorphic to the quantized enveloping algebra of the $3$-dimensional 
Lie algebra $\fsl_2$. We have the relations:
\begin{gather*}
 K_iE_iK_i^{-1}=v^2E_i, \qquad K_iF_iK_i^{-1}=v^{-2}F_i,\\
E_iF_i-F_iE_i=\frac{K_i-K_i^{-1}}{v-v^{-1}}.
\end{gather*}
We will study properties of a module for $U_v(\fg)$ by restricting the 
action to $U_{v,i}$ and using known properties of representations of $U_{v,i}
\cong U_v(\fsl_2)$. Recall that, for any $m\geq 0$, there are precisely 
two simple $U_{v,i}$-modules of dimension $m+1$ (up to isomorphism), 
denoted by $V_m^{\pm}$. We have that $V_m^{\pm}$ is a cyclic module 
generated by a vector $u_0$ such that
\[ E_i.u_0=0 \qquad \mbox{and}\qquad K_i.u_0=\pm v^m u_0.\]
See \cite[\S 1.3]{Kash} for further details on the construction of
$V_m^{\pm}$. 

Now let $M$ be an integrable $U_v(\fg)$-module. We do not need to recall
the precise definition here, but what is important is the fact that,
regarded as a $U_{v,i}$-module, $M$ can be written as a direct sum of 
modules of the form $V_m^+$; see \cite[1.4.1 and 3.3.1]{Kash}. The
explicit knowledge of the action of $U_{v,i}$ on $V_m^+$ implies that
we have a direct sum decomposition
\[ M=\bigoplus_{n \in \ZZ} M_n \qquad \mbox{where}\qquad
M_n:=\{u \in M \mid K_i.u=v^nu\};\]
see \cite[\S 1.4]{Kash}. For a fixed $n\in \ZZ$, we have $\dim M_n<\infty$ 
and any $u\in M_n$ has a unique expression 
\[ u=\sum_{m \geq 0,-n} \frac{1}{[m]!}F_i^m.u_m \quad \mbox{where $u_m 
\in M_{n+2m}$ is such that $E_i.u_m=0$}.\]
Here, $[0]!=1$ and $[m]!=[m-1]!(v^m-v^{-m})/(v-v^{-1})$ for $m>0$. For
these facts, see \cite[\S 2.2]{Kash}. 

\begin{defn}[Kashiwara; see \protect{\cite[2.2.4]{Kash}}] \label{defkash}
For $0 \leq i \leq l-1$, we define linear maps $\tilde{E}_i,\tilde{F}_i
\colon M \rightarrow M$ by 
\begin{align*}
\tilde{E}_i.u&=\sum_{m \geq 1,-n}  \frac{1}{[m{-}1]!}F_i^{m-1}.u_m 
\qquad (u \in M_n),\\
\tilde{F}_i.u&=\sum_{m \geq 0,-n}  \frac{1}{[m{+}1]!}F_i^{m+1}.u_m 
\qquad (u \in M_n).
\end{align*}
\end{defn}

Next, let $R\subseteq {\CC}(v)$ be the local ring of all rational
functions with no pole at $v=0$. A {\em local basis} of $M$ is a
pair $(\cL,\cB)$ where $\cL$ is an $R$-submodule of $M$ such that
$M$ is a free as an $R$-module, ${\CC}(v)\otimes_R \cL=M$, and $\cB$ is 
a basis of the $\CC$-vector space $\cL/v\cL$. Note that a complete set 
of representatives of $\cB$ in $\cL$ forms a ${\CC}(v)$-basis of $M$. 
A {\em crystal basis} for $M$ is a local basis $(\cL,\cB)$ satisfying 
certain additional conditions. These conditions include the following 
(see \cite[Theorem~4.1.2]{Kash}):
\begin{itemize}
\item we have $\tilde{E}_i.\cL\subseteq \cL$ and $\tilde{F}_i.\cL
\subseteq \cL$ for $0 \leq i \leq l-1$; hence, we obtain actions
of $\tilde{E}_i$ and $\tilde{F}_i$ on $\cL/v\cL$ which we denote
by the same symbols;
\item we have $\tilde{E}_i.\cB \subseteq \cB\cup\{0\}$ and 
$\tilde{F}_i.\cB \subseteq \cB \cup\{0\}$ for $0 \leq i \leq l-1$;
\item for $b,b'\in \cB$ and $0 \leq i \leq n-1$, we have 
$b'=\tilde{F}_i.b \Leftrightarrow \tilde{E}_i.b'=b$.
\end{itemize}
Thus, we can define a graph with vertices indexed by the elements
of $\cB$; for $b\neq b'$ in $\cB$ and $0 \leq i \leq n-1$, we have 
a colored oriented edge 
\[ b \stackrel{i}{\longrightarrow} b' \qquad \mbox{if and only if}
\qquad b'=\tilde{F}_i.b.\]
This graph is the {\em crystal graph} of $\cB$.  By \cite[Chap.~9]{Ar3},
if $M$ is an integrable highest weight module, then $M$ has a crystal 
basis which is unique up to scalar multiples; furthermore, using Lusztig's 
canonical basis $\bB$ of $U_v^-(\fg)$, we obtain a crystal basis of $M$ 
by setting
\[ \cL:=\sum_{b \in \bB} Rb.u_0 \qquad \mbox{and}
\qquad \cB:=\{b.u_0+v\cL\mid b \in \bB\}\setminus
\{0\},\] 
where $u_0\in M$ is a highest weight vector.

Now the problem is to obtain an explicit description of the crystal graph 
for our highest weight module $M_v^{(r)}(\uvar)$. For
this purpose, we first describe a crystal basis of the integrable module
$\fF_v^{(r)}$. The corresponding graph will have vertices labelled by the 
set $\Pi_r:=\bigcup_{n\geq 0} \Pi_{r,n}$ and it will have the property that 
the connected component of $\uvar$ is the crystal graph 
of $M_v^{(r)}(\uvar)$.  

\begin{thm}[Jimbo et al. \protect{\cite{JMMO}}, Foda et al. 
\protect{\cite{FLOTW}}, Uglov \protect{\cite{Ug}}] \label{crys1}  
We set
\[ \cL_r:=\sum_{\ulambda \in \Pi_r} R \,
\ulambda \subseteq \fF_v^{(r)} \qquad \mbox{and}\qquad
\cB_r:=\{b_{\ulambda} \mid \ulambda 
\in \Pi_r\},\]
where $b_{\ulambda}$ denotes the image of $\ulambda$
in $v\cL_r$. Then $(\cL_r,\cB_r)$ is a crystal basis for $\fF_v^{(r)}$.
Given $\ulambda,\umu \in \Pi_r$ and $i\in \{0,1,
\ldots,l-1\}$, we have $b_{\ulambda} \stackrel{i}{\rightarrow} 
b_{\umu}$ if and only if $\umu$ is obtained from $\ulambda$ by adding 
a so-called ``good'' $i$-node.
\end{thm}

The ``good'' nodes are defined as follows. Let $\ulambda\in \Pi_r$ and 
let $\gamma$ be an $i$-node of $\ulambda$. We say that $\gamma$ is a 
{\em normal} node if, whenever $\gamma'$ is an $i$-node of $\ulambda$ 
below $\gamma$, there are strictly more removable $i$-nodes between 
$\gamma'$ and $\gamma$ than there are addable $i$-nodes between 
$\gamma'$ and $\gamma$. If $\gamma$ is a highest normal $i$-node of 
$\ulambda$, then $\gamma$ is called a {\em good} node. Note that these 
notions heavily depend on the definition of what it means for one node 
to be ``above'' another node. These definitions (for $r=1$) first appeared 
in the work of Kleshchev \cite{Klesh} on the modular branching rule for the 
symmetric group; see also the discussion of these results in 
\cite[\S 2]{LLT}.

For any $n\geq 0$, we define a subset $\Lambda_{r,n}^{(\bu)}
\subseteq \Pi_{r,n}$ recursively as follows. We set 
$\Lambda_{r,0}^{(\bu)}=\{\uvar\}$. For $n\geq 1$, the set 
$\Lambda_{r,n}^{(\bu)}$ is constructed as follows.
\begin{itemize}
\item[(1)] We have $\uvar \in \Lambda_{r,n}^{(\bu)}$; 
\item[(2)] Let $\ulambda\in \Pi_{r,n}$. Then $\ulambda$
belongs to $\Lambda_{r,n}^{(\bu)}$ if and only if 
$\ulambda/\umu=\gamma$ where $\umu
\in \Lambda_{r,n-1}^{(\bu)}$ and $\gamma$ is a good $i$-node
of $\ulambda$ for some $i\in \{0,1,\ldots,l-1\}$.
\end{itemize}
Thus, the set $\Lambda_r^{(\bu)}:=\bigcup_{n\geq 0} \Lambda_{r,n}^{(\bu)}$
labels the vertices in the connected component containing $b_{\uvar}$
of the crystal graph of $\fF_v^{(r)}$. Hence, by general results on crystal 
bases (see the ``unicity theorem'' in \cite[4.1.5]{Kash}), we have:

\begin{cor} \label{crys2} The crystal graph of $M_v^{(r)}(\uvar)$ has
vertices labelled by the elments in $\Lambda_r^{(\bu)}$. The edges
$b_{\ulambda} \stackrel{i}{\rightarrow} b_{\umu}$ are given as in
Theorem~\ref{crys1}.
\end{cor}

The above result provides a purely combinatorial description of the
crystal graph of $M_v^{(r)}(\uvar)$. An example is given in 
Table~\ref{tabcr}. In that example, we have 
\[ \Lambda_{2,3}^{(0,1)}=\{((3),\varnothing),\; ((2),(1)),\; ((1,(2)),\;
(\varnothing,(3))\}.\]
If we had used the ordering of nodes $\prec$ in Remark~\ref{conven},
we would obtain a quite different labelling of the vertices and edges in
the crystal graph; in that case, $\Lambda_{2,3}^{(0,1)}$ would have to
be replaced by the set (see \cite[\S 1.3]{Jac0}):
\[ \{((3),\varnothing),\; ((2,1),\varnothing),\; ((1),(2)),\; ((2),(1))\}.\]
In Section~8, we will see interpretations of these different labellings
in terms of modular representations of Iwahori--Hecke algebras of type
$B_n$. 

\begin{table}[htbp] \caption{Part of the crystal graph for $l=r=2$,
$\bu=(0,1)$} \label{tabcr} 
\begin{center}
\begin{picture}(290,190)
\put(  4,11){$\bigl($}
\put( 13,11){$0$}
\put( 23,11){$1$}
\put( 33,11){$0$}
\put( 10,10){\line(1,0){30}}
\put( 10,20){\line(1,0){30}}
\put( 10,10){\line(0,1){10}}
\put( 20,10){\line(0,1){10}}
\put( 30,10){\line(0,1){10}}
\put( 40,10){\line(0,1){10}}
\put( 43,11){$,$}
\put( 52,11){$\varnothing\bigr)$}
\put( 84,11){$\bigl($}
\put( 93,11){$0$}
\put(103,11){$1$}
\put(123,11){$1$}
\put( 90,10){\line(1,0){20}}
\put( 90,20){\line(1,0){20}}
\put( 90,10){\line(0,1){10}}
\put(100,10){\line(0,1){10}}
\put(110,10){\line(0,1){10}}
\put(113,11){$,$}
\put(120,10){\line(1,0){10}}
\put(120,20){\line(1,0){10}}
\put(120,10){\line(0,1){10}}
\put(130,10){\line(0,1){10}}
\put(132,11){$\bigr)$}
\put(154,11){$\bigl($}
\put(163,11){$0$}
\put(183,11){$1$}
\put(193,11){$0$}
\put(160,10){\line(1,0){10}}
\put(160,20){\line(1,0){10}}
\put(160,10){\line(0,1){10}}
\put(170,10){\line(0,1){10}}
\put(173,11){$,$}
\put(180,10){\line(1,0){20}}
\put(180,20){\line(1,0){20}}
\put(180,10){\line(0,1){10}}
\put(190,10){\line(0,1){10}}
\put(200,10){\line(0,1){10}}
\put(202,11){$\bigr)$}
\put(224,11){$\bigl($}
\put(253,11){$1$}
\put(263,11){$0$}
\put(273,11){$1$}
\put(230,11){$\varnothing$}
\put(241,11){$,$}
\put(250,10){\line(1,0){30}}
\put(250,20){\line(1,0){30}}
\put(250,10){\line(0,1){10}}
\put(260,10){\line(0,1){10}}
\put(270,10){\line(0,1){10}}
\put(280,10){\line(0,1){10}}
\put(282,11){$\bigr)$}

\put( 60,55){\vector(-1,-1){25}}
\put(103,55){\vector(+1,-3){8}}
\put(230,55){\vector(+1,-1){25}}
\put(182,55){\vector(-1,-3){8}}
\put( 39,45){$0$}
\put( 95,40){$1$}
\put(169,42){$0$}
\put(245,45){$1$}
\put( 64,61){$\bigl($}
\put( 73,61){$0$}
\put( 83,61){$1$}
\put( 70,60){\line(1,0){20}}
\put( 70,70){\line(1,0){20}}
\put( 70,60){\line(0,1){10}}
\put( 80,60){\line(0,1){10}}
\put( 90,60){\line(0,1){10}}
\put( 93,61){$,$}
\put(100,61){$\varnothing$}
\put(111,61){$\bigr)$}

\put(174,61){$\bigl($}
\put(180,61){$\varnothing$}
\put(191,61){$,$}
\put(203,61){$1$}
\put(213,61){$0$}
\put(200,60){\line(1,0){20}}
\put(200,70){\line(1,0){20}}
\put(200,60){\line(0,1){10}}
\put(210,60){\line(0,1){10}}
\put(220,60){\line(0,1){10}}
\put(222,61){$\bigr)$}

\put(100,115){\vector(-1,-3){12}}
\put(185,115){\vector(+1,-3){12}}
\put( 87,100){$1$}
\put(194,100){$0$}
\put( 84,121){$\bigl($}
\put( 93,121){$0$}
\put( 90,120){\line(1,0){10}}
\put( 90,130){\line(1,0){10}}
\put( 90,120){\line(0,1){10}}
\put(100,120){\line(0,1){10}}
\put(103,121){$,$}
\put(110,121){$\varnothing$}
\put(120,121){$\bigr)$}
\put(163,121){$\bigl($}
\put(170,121){$\varnothing$}
\put(182,121){$,$}
\put(193,121){$1$}
\put(190,120){\line(1,0){10}}
\put(190,130){\line(1,0){10}}
\put(190,120){\line(0,1){10}}
\put(200,120){\line(0,1){10}}
\put(202,121){$\bigr)$}

\put(130,175){\vector(-2,-3){25}}
\put(150,175){\vector(+2,-3){25}}
\put(111,160){$0$}
\put(164,160){$1$}
\put(124,181){$\bigl($}
\put(130,181){$\varnothing,\varnothing$}
\put(154,181){$\bigr)$}
\end{picture}
\end{center}
(The numbers inscribed in the boxes of the diagrams are the $l$-residues.)
\end{table}

\begin{Par} {\bf FLOTW-partitions.} \label{crys3}  Assume that the
parameters in $\bu$ satisfy the condition
\[ 0 \leq u_1 \leq u_2 \leq \cdots \leq u_r\leq l-1.\]
Then it is shown in Foda et al. \cite[2.11]{FLOTW} that $\ulambda
\in \Pi_{r,n}$ belongs to $\Lambda_{r,n}^{(\bu)}$ if and only if the
following conditions are satisfied:
\begin{itemize}
\item[(a)] For all $1\leq j \leq r-1$ and $i=1,2,\ldots$, we have:
\[ \lambda_{(j+1),i} \geq \lambda_{(j),i+u_{j+1}-u_j} \qquad \mbox{and}
\qquad \lambda_{(1),i} \geq \lambda_{(r),i+l+u_1-u_r};\]
\item[(b)] for all $k>0$, among the residues appearing at the right ends
of the rows of $[\ulambda]$ of length $k$, at least one element of
$\{0,1,\ldots,l-1\}$ does not occur.
\end{itemize}
Note that this provides a non-recursive description of the elements
of $\Lambda_r^{(\bu)}$.
\end{Par}

\section{The theorems of Ariki and Jacon} \label{MGsec7}

The aim of this section is to explain the applications of the results
on the canonical basis of the Fock space $\fF^{(r)}$ to the problem
of parametrizing the simple modules of non-semisimple Iwahori--Hecke 
algebras.  These results actually hold for a wider class of ``Hecke
algebras'', which we now introduce.

\begin{Par} {\bf Ariki--Koike algebras.} \label{arko} Let $k$ be an 
algebraically closed field and let $\zeta_l$ be an element of order 
$l \geq 2$ in $k^\times$. Let $r,n\geq 1$ and  fix parameters
\[ \bu=(u_1,\ldots,u_r) \qquad \mbox{where} \qquad  u_i\in \ZZ.\]
Having fixed these data, we let $\bH_{r,n}^{(\bu)}$ be the associative 
$k$-algebra (with $1$), with generators $S_0,S_1,\ldots, S_{n-1}$ and 
defining relations as follows:
\begin{gather*} S_0S_1S_0S_1=S_1S_0S_1S_0 \quad \mbox{and}\quad 
S_0S_i=S_iS_0 \quad \mbox{(for $i>1$)},\\ S_iS_j=S_jS_i \quad 
\mbox{(if $|i-j|>1$)}, \\ S_iS_{i+1} S_i=S_{i+1}S_iS_{i+1} \quad 
\mbox{(for $1 \leq i \leq n-2$)},\\ (S_0-\zeta_l^{u_1})(S_0-
\zeta_l^{u_2})\cdots (S_0-\zeta_l^{u_r})=0,\\ (S_i-\zeta_l)(S_i+1)=0 
\quad \mbox{for $1 \leq i \leq n-1$)}.
\end{gather*}
This algebra can be seen as an Iwahori--Hecke algebra associated with 
the group $G_{r,n}$ (which is a Coxeter group only for $r=1,2$); see
Ariki \cite[Chap.~13]{Ar3} and Brou\'e--Malle \cite{BM2} for further 
details and motivations for studying this class of algebras.

Here, we will be mostly interested in the case $r=2$, when 
$\bH_{2,n}^{(\bu)}$ can be identified with an Iwahori--Hecke algebra of 
type $B_n$. Indeed, let $\zeta_{2l}\in k^\times$ be a square root of 
$\zeta_l$ such that
\[ \zeta_{2l}^2=\zeta_l \qquad \mbox{and}\qquad \zeta_{2l}^l=-1.\]
Then, setting $T:=\zeta_{2l}^{l-2u_2}S_0$, we obtain 
\begin{align*}
T^2&=\zeta_{2l}^{f} +(\zeta_{2l}^{f}-1)T \qquad \mbox{where}\qquad
f=l+2(u_1-u_2),\\ S_i^2&=\zeta_{2l}^2+(\zeta_{2l}^2-1)S_i 
\qquad \mbox{for $1\leq i\leq n-1$}.
\end{align*}
Thus, the map $T\mapsto T_t$, $S_i \mapsto T_{s_i}$ defines an
isomorphism of $k$-algebras
\[ \bH_{2,n}^{(\bu)} \cong \bH_{k,\xi}=k 
\otimes\bH_A(W_n,L), \qquad \xi=\zeta_{2l},\]
where $W_n$ is defined as in Example~\ref{schurbn} and $L$ is a weight
function such that 
\[ L(t) \equiv l+2(u_2-u_1) \bmod 2l \quad \mbox{and}\quad L(s_i) 
\equiv 2 \bmod 2l.\]
\end{Par}

As mentioned at the beginning of the previous section, we have a 
labelling of the simple $\CC G_{r,n}$-modules by $\Pi_{r,n}$. We write 
this as
\[ \Irr_{\CC}(G_{r,n})=\{E^{\ulambda} \mid \ulambda
\in \Pi_{r,n}\}.\]
Now Dipper--James--Mathas \cite[\S 3]{DJMa} have generalized the theory 
of {\em Specht modules} for Iwahori--Hecke algebra of the symmetric 
group to $\bH_{r,n}^{(\bu)}$; see Example~\ref{canbn} for the case $r=2$
and also Graham--Lehrer \cite{GrLe}.
Thus, for any $\ulambda\in \Pi_{r,n}$, there is a {\em Specht module} 
$S^\ulambda \in \bH_{r,n}^{(\bu)}$. Each $S^\ulambda$ carries a 
symmetric bilinear form and, taking quotients by the radical, we obtain 
a collection of modules $D^{\ulambda}$. As before, we set 
$\Lambda_{r,n}^\circ:=\{\ulambda \vdash n\mid D^\ulambda \neq \{0\}\}$. 
Then, by \cite[Theorem~3.30]{DJMa},  we have 
\[ \Irr(\bH_{r,n}^{(\bu)})=\{D^\ulambda \mid \ulambda \in 
\Lambda_{r,n}^\circ\}.\]
Furthermore, the entries of the decomposition matrix
\[ D=\bigl([S^\ulambda :D^{\umu}]\bigr)_{\ulambda \in \Pi_{r,n},\umu
\in \Lambda_{r,n}^\circ}\]
satisfy conditions analogous to those in Example~\ref{canbn}, where
one has to consider the dominance order on $r$-tuples of partitions
defined in \cite[3.11]{DJMa}. 

We are now ready to state Ariki's theorem which establishes the link 
to the Fock space $\fF^{(r)}$. For any $\umu \in \Lambda_{r,n}^\circ$, 
we define an element of $\fF^{(r)}$ by 
\[ P_\umu:=\sum_{\ulambda \vdash n} [S^\ulambda:D^\umu]\, \ulambda
\in \fF_{r,n}.\]
Then we consider the subspace of $\fF^{(r)}$ generated by all these
elements:
\[\cM^{(r)}:=\langle P_\umu \mid \umu \in \Lambda_{r,n}^\circ 
\mbox{ for some $n\geq 0$}\rangle_{\CC} \subseteq \fF^{(r)},\]
where $P_{\uvar}=\uvar \in \fF_{r,0}$.

\begin{thm}[Ariki \protect{\cite[Theorem~12.5]{Ar3}}] \label{mainar}
Assume that the characteristic of $k$ is zero. Then we have 
\[ \cM^{(r)}=M^{(r)}(\uvar) \quad \mbox{and}\quad \bB_1(\uvar)=
\{P_\umu \mid \umu \in \Lambda_{r,n}^\circ \mbox{ for some $n\geq 0$}\},\]
where $M^{(r)}(\uvar)$ is the highest weight module for $U(\hat{\fsl}_l)$
as in Proposition~\ref{hash1} and $\bB_1(\uvar)$ is the Kashiwara--Lusztig 
canonical basis of $M^{(r)}(\uvar)$ (at $v=1$).
\end{thm}

The importance of this result for the modular representation theory
of Iwahori--Hecke algebras (and for finite groups of Lie type, via
the results in Section~1 and~2) can hardly be overestimated. There are 
efficient and purely combinatorial algorithms for computing the 
canonical basis; see Lascoux--Leclerc--Thibon \cite{LLT} (the 
``LLT algorithm'') for the case $r=1$ and Jacon \cite{Jac3} for $r\geq 2$. 
Hence, if $k$ has characteristic zero, the above result shows that 
these algorithms compute the decomposition numbers of $\bH_{r,n}^{(\bu)}$. 

If the characteristic of $k$ is not zero, the elements $P_\umu$ will
no longer coincide with the elements of the canonical basis. But the
first part of the above statement remains valid:

\begin{thm}[Ariki--Mathas \protect{\cite{ArMa}}] \label{arama} We have
$\cM^{(r)}=M^{(r)}(\uvar)$ and 
\[|\Irr(\bH_{r,n}^{(\bu)})|=\dim_{\CC} (\fF_{r,n} \cap M^{(r)}(\uvar)).\]
In particular, the number of simple modules of $\bH_{r,n}^{(\bu)}$
only depends on $l$ and on the congruence classes modulo $l$ of the
numbers $u_1,\ldots,u_r$.
\end{thm}

Using the above results, Ariki identified the indexing set
$\Lambda_{r,n}^\circ$ arising from the theory of Specht modules: 

\begin{thm}[Ariki \protect{\cite[Theorem~4.3]{Ar2}}] \label{arkl} 
We have $\Lambda_{r,n}^\circ=\Lambda_{r,n}^{(\bu)}$, where $\bu=(u_1,
\ldots,u_r)$ is chosen such that $u_i-u_{i+1}>n-1$ for all $i$.
\end{thm}

The condition $u_i-u_{i+1}>n-1$ really means that we are working with the
order $\prec$ on nodes defined in Remark~\ref{conven}. Then the set
$\Lambda_{r,n}^{(\bu)}$ is precisely the set of so-called {\em Kleshchev
bipartations}; see also Ariki \cite[\S 12.2]{Ar3}.

Now the algebra $\bH_{r,n}^{(\bu)}$ also is symmetric, with a trace
function $\tau\colon \bH_{r,n}^{(\bu)} \rightarrow k$ satisfying
properties analogous to those for Iwahori--Hecke algebras of finite
Coxeter groups; see Bremke--Malle \cite{BrMa}, Malle--Mathas \cite{MaMa}.
Working with a suitable generic version of $\bH_{r,n}^{(\bu)}$, we have 
corresponding Laurent polynomials $\bc_{E^\ulambda}$ for any $\ulambda
\in \Pi_{r,n}$. Explicit combinatorial formulae for $\bc_{E^\ulambda}$,
generalizing those in Example~\ref{schurbn} for type $B_n$, are obtained 
by Geck--Iancu--Malle \cite{GIM} and, independently, by Mathas \cite{Mat}
(proving a conjecture of Malle \cite{MaCon}). 

Thus, it makes sense to consider the existence of ``canonical basic 
sets'' for $\bH_{r,n}^{(\bu)}$, as in Definition~\ref{canbase}. Recently,
Jacon \cite{Jac0}, \cite{Jac4} has shown that such canonical basic sets 
indeed exist for $\bH_{r,n}^{(\bu)}$. (His results still rely on
Ariki's Theorem~\ref{mainar}.) These are labelled by sets 
$\Lambda_{r,n}^{(\bu)}$ for various choices of $\bu$ where the 
non-recursive description of Foda et al. in (\ref{crys3}) applies. 

Let us discuss some applications of Jacon's results to the case $r=2$. 
Then $\bH_{2,n}^{(\bu)}$ can be identified with an Iwahori--Hecke algebra
of type $B_n$, as explained above. So let $W_n$ 
be a Coxeter group of type $B_n$, with weight function $L \colon W_n 
\rightarrow \NN$ determined by two integers $a,b\geq 0$ as in 
Example~\ref{schurbn}:
\begin{center}
\begin{picture}(250,40)
\put( 10,25){$B_n$}
\put( 10,05){$L$ :}
\put( 65,05){$b$}
\put( 95,05){$a$}
\put(125,05){$a$}
\put(215,05){$a$}
\put( 61,33){$t$}
\put( 91,33){$s_1$}
\put(121,33){$s_2$}
\put(205,33){$s_{n-1}$}
\put( 65,25){\circle*{5}}
\put( 95,25){\circle*{5}}
\put(125,25){\circle*{5}}
\put(215,25){\circle*{5}}
\put( 65,27){\line(1,0){30}}
\put( 65,23){\line(1,0){30}}
\put( 95,25){\line(1,0){50}}
\put(160,25){\circle*{1}}
\put(170,25){\circle*{1}}
\put(180,25){\circle*{1}}
\put(195,25){\line(1,0){20}}
\end{picture}
\end{center}
Let $\theta \colon A \rightarrow k$ be a specialization into a field $k$ 
of characteristic $\neq 2$; let $\xi=\theta(v^2)$. In Example~\ref{canbn},
we have already seen a convenient parametrization of the simple modules
of $\bH_{k,\xi}$ in the cases where $f_n(a,b)\neq 0$ or $\xi^a=1$.
So let us now assume that
\[ f_n(a,b)=0 \qquad \mbox{and}\qquad \xi^a\neq 1.\]
Then we have 
\[ \xi^{b+ad}=-1 \qquad \mbox{for some $d\in \ZZ$ such that $|d|<n-1$}.\]
In particular, $\xi$ is a root of unity of even order. Let $l\geq 2$ be
the multiplicative order of $\xi^a$ and set $\zeta_l:=\xi^a$.
Then $\xi^b=-\xi^{-ad}=-\zeta_l^{-d}$.  Let $\zeta_{2l}\in k$ be 
a square root of $\zeta_l$ such that $\zeta_{2l}^l =-1$. Thus, we have the 
quadratic relations
\[T_t^2=\zeta_{2l}^{l-2d} +(\zeta_{2l}^{l-2d}-1)T_t \qquad 
\mbox{and} \qquad T_{s_i}=\zeta_{2l}^2+(\zeta_{2l}^2-1)T_{s_i} \]
for $1\leq i\leq n-1$. Consequently, we have 
\[ \bH_{k,\xi}\cong \bH_{2,n}^{(\bu)} \qquad \mbox{where} \qquad
\bu=(u_1,u_2),\quad u_2-u_1 \equiv d \bmod l.\]

\begin{thm}[The equal parameter case; Jacon 
\protect{\cite[Theorem~3.2.3]{Jac0}}] \label{jac01} In the above setting, 
assume that $a=b=1$, $f_n(1,1)=0$ and $\xi\neq 1$.  Then $\xi$ has even 
order $l\geq 2$ and 
\[ \cB_{k,\xi}:=\{E^\ulambda \mid\ulambda\in\Lambda_{2,n}^{(\bu)}\}\]
is a canonical basic set for $\bH_{k,\xi}$, where $\bu=(1,l/2)$.
See (\ref{crys3}) for an explicit description of $\Lambda_{2,n}^{(\bu)}$.
\end{thm}

(Note that, in this case, we have $d+1 \equiv l/2 \bmod l$.)

\begin{thm}[The case $b=0$; Jacon \protect{\cite[Theorem~3.2.5]{Jac0}}] 
\label{jac02} In the above setting, assume that $b=0$, $a=1$, 
$f_n(1,0)=0$ and $\xi\neq 1$. Then $\xi$ has even order $l\geq 2$ and 
\[ \cB_{k,\xi}:=\{E^\ulambda \mid\ulambda\in\Lambda_{2,n}^{(\bu)}\}\]
is a canonical basic set for $\bH_{k,\xi}$, where $\bu=(0,l/2)$.
See (\ref{crys3}) for an explicit description of $\Lambda_{2,n}^{(\bu)}$.
\end{thm}

(Note that, in this case, we have $d \equiv l/2 \bmod l$.)

The above result yields a description of the canonical basic set for
Iwahori--Hecke algebras of type $D_n$ ($n\geq 2$). Indeed, let $W_n'$ 
be the subgroup of $W_n$ generated by $s_0,s_1,\ldots,s_{n-1}$ where 
$s_0=ts_1t$.  As already pointed out in Example~\ref{schurbn}, $W_n'$ 
is a Coxeter group of type $D_n$. Furthermore, let $L'$ be the 
restriction of the weight function $L$ in Theorem~\ref{jac02} to $W_n'$; 
then $L'(s_i)=1$ for $0\leq i \leq n-1$. We denote the corresponding
Iwahori--Hecke algebra by $\bH'$. Let $\cB_{k,\xi}'\subseteq 
\Irr_{\QQ}(W_n')$ be the canonical basic set for $\bH_{k,\xi}'$. 
By \cite[Theorem~5.5]{my00}, we have 
\[ \cB_{k,\xi}'=\{E' \in \Irr_{\QQ}W_n') \mid \mbox{$E'$ occurs in the
restriction of some $E\in \cB_{k,\xi}$}\},\]
where $\cB_{k,\xi}$ is as in Theorem~\ref{jac02}. Using the information
on the restriction of modules from $W_n$ to $W_n'$ in Example~\ref{schurbn},
this yields:

\begin{thm}[Jacon \protect{\cite[Theorem~3.2.7]{Jac0} and \cite{Jac1}}] 
\label{mudn} Let $W_n'$ be a Coxeter group of type $D_n$ ($n\geq 2$) and 
$L'(s_i)=1$ for $0\leq i \leq n-1$. Assume that $\xi$ has even order, 
$l\geq 2$ say. Then 
\begin{align*}
\bB_{k,\xi}'&=\{E^{[\lambda,\mu]} \mid (\lambda,\mu)\in 
\Lambda_{2,n}^{(\bu)},\lambda \neq \mu\}\\&\qquad \cup 
\{E^{[\lambda,\pm]} \mid \mbox{ $n$ even and $\lambda \vdash n/2$ is 
$l/2$-regular}\}
\end{align*}
is a canonical basic set for $\bH_{k,\xi}'$, where $\bu=(0,l/2)$.
\end{thm}

\begin{table}[htbp] \label{tabb3} \caption{Decomposition numbers for 
$B_3$ with $\xi^b=1$, $\xi=-1$.}
\begin{center}
\renewcommand{\arraystretch}{1.1}
$\begin{array}{lr} \begin{array}{|c|c|cccc|} 
\hline \ulambda & \dim E^{\ulambda} & 
\multicolumn{4}{c|}{[E:M]} \\ \hline
(3,\varnothing)    & 1 & 1 & . & . & . \\
(\varnothing,3 )   & 1 & . & 1 & . & . \\
(111,\varnothing)  & 1 & 1 & . & . & . \\
(\varnothing,111)  & 1 & . & 1 & . & . \\ 
(21,\varnothing)   & 2 & . & . & 1 & . \\
(\varnothing,21)   & 2 & . & . & . & 1 \\
(1,2)              & 3 & 1 & . & . & 1 \\
(2,1)              & 3 & . & 1 & 1 & . \\
(11,1)             & 3 & . & 1 & 1 & . \\
(1,11)             & 3 & 1 & . & . & 1 \\
\hline\end{array}& \quad 
\begin{array}{|cc|c|cccc|} 
\hline \multicolumn{2}{|c|}{b=0}& \alpha_{\ulambda} & 
\multicolumn{4}{c|}{[E:M]} \\ \hline
\rightarrow&(3,\varnothing)    & 0 & 1 & . & . & . \\
\rightarrow&(\varnothing,3 )   & 0 & . & 1 & . & . \\
\rightarrow&(1,2)              & 1 & 1 & . & 1 & . \\
\rightarrow&(2,1)              & 1 & . & 1 & . & 1 \\
&(21,\varnothing)   & 2 & . & . & . & 1 \\
&(\varnothing,21)   & 2 & . & . & 1 & . \\
&(11,1)             & 3 & . & 1 & . & 1 \\
&(1,11)             & 3 & 1 & . & 1 & . \\
&(111,\varnothing)  & 6 & 1 & . & . & . \\
&(\varnothing,111)  & 6 & . & 1 & . & . \\ 
\hline\end{array}\\ \\  
\begin{array}{|cc|r|cccc|} 
\hline \multicolumn{2}{|c|}{b=2}& \alpha_{\ulambda} & 
\multicolumn{4}{c|}{[E:M]} \\ \hline
\rightarrow &(3,\varnothing)    & 0 & 1 & . & . & . \\
\rightarrow &(21,\varnothing)   & 1 & . & 1 & . & . \\
\rightarrow &(2,1)              & 2 & . & 1 & 1 & . \\
&(11,1)             & 3 & . & 1 & 1 & . \\
&(111,\varnothing)  & 3 & 1 & . & . & . \\
&(\varnothing,3 )   & 3 & . & . & 1 & . \\
\rightarrow &(1,2)              & 3 & 1 & . & . & 1 \\
&(1,11)             & 6 & 1 & . & . & 1 \\
&(\varnothing,21)   & 7 & . & . & . & 1 \\
&(\varnothing,111)  & 12& . & . & 1 & . \\ 
\hline\end{array}\quad&
\begin{array}{|cc|r|cccc|} 
\hline \multicolumn{2}{|c|}{b=4}& \alpha_{\ulambda} & 
\multicolumn{4}{c|}{[E:M]} \\ \hline
\rightarrow&(3,\varnothing)    & 0 & 1 & . & . & . \\
\rightarrow&(21,\varnothing)   & 1 & . & 1 & . & . \\
&(111,\varnothing)  & 3 & 1 & . & . & . \\
\rightarrow&(2,1)              & 4 & . & 1 & 1 & . \\
&(11,1)             & 5 & . & 1 & 1 & . \\
\rightarrow&(1,2)              & 7 & 1 & . & . & 1 \\
&(\varnothing,3 )   & 9 & . & . & 1 & . \\
&(1,11)             & 10& 1 & . & . & 1 \\
&(\varnothing,21)   & 13& . & . & . & 1 \\
&(\varnothing,111)  & 18& . & . & 1 & . \\ 
\hline\end{array}\end{array}$
\end{center}
\end{table}

In Table~3, we consider the case where $n=3$, $\xi=-1$, $a=1$ and $b\geq 0$
is even. Note that, for any such $b$, we have $\xi^b=1$ and so we
are always dealing with the same algebra $\bH_{k,\xi}$. In the first
matrix, the rows are just ordered according to increasing value of 
$\dim E^\ulambda$. Now, depending on which value of $b$ we take, we 
obtain a different canonical basic set (indicated by ``$\rightarrow$''). 
The one for $b=0$ is given by Theorem~\ref{jac02}, and it yields a 
canonical basic set for type $D_3$; the one for $b=4$ corresponds to 
the ``asymptotic case'' in Example~\ref{canbn}, and it yields the indexing
set $\Lambda_{2,3}^\circ$ (arising from the Specht module theory). It 
turns out that this is also obtained for $b=2$.

In the above two results, we only considered the cases where $a=1$ and 
$b\in \{0,1\}$. But similar arguments apply to all choices of $a,b\geq 0$
where $f_n(a,b)=0$, and this yields explicit combinatorial descriptions
of canonical basic sets in terms of the sets $\Lambda_{r,n}^{(\bu)}$
arising from crystal graphs, for suitable values of $\bu$. For the details,
see \cite{GeJa}. 


\end{document}